\documentclass[english,12pt]{article}%
\usepackage{bbm}
\usepackage{graphicx}
\usepackage{subcaption}
\usepackage{makeidx}
\usepackage{amssymb}
\usepackage[all]{xy}
\usepackage{float}
\usepackage[latin1]{inputenc}
\usepackage{amsfonts}
\usepackage{amsmath}
\usepackage{amssymb}
\usepackage{url}
\usepackage{hyperref}
\usepackage{mathabx}
\usepackage{wrapfig}
\usepackage{dsfont} 
\usepackage{resizegather}
\usepackage{yfonts}
\usepackage{authblk}
\usepackage{geometry}
\usepackage[dvipsnames]{xcolor}

\usepackage{amsthm}
\usepackage{url}
\usepackage{hyperref}
\providecommand{\U}[1]{\protect\rule{.1in}{.1in}}

\usepackage{xcolor}

\newtheorem{prop}{Proposition}[section]

\newtheorem{defi}[prop]{Definition}

\newtheorem{rmk}[prop]{Remark}
\newtheorem{lem}[prop]{Lemma}

\newtheorem{theo}[prop]{Theorem}
\newtheorem{examp}[prop]{Example}

\newcommand{\tr}{\mbox{\rm Tr}}

\newcommand{\EE}{\mathbb{E}}

\newcommand{\HH}{\mathbb{H}}

\newcommand{\LL}{\mathbb{L}}

\newcommand{\NN}{\mathbb{N}}
\newcommand{\PP}{\mathbb{P}}

\newcommand{\RR}{\mathbb{R}}

\newcommand{\TT}{\mathbb{T}}
\newcommand{\UU}{\mathbb{U}}

\newcommand{\XX}{\mathbb{X}}

\newcommand{\Aa}{ {\cal A }}
\newcommand{\Ba}{ {\cal B }}
\newcommand{\Ca}{ {\cal C }}
\newcommand{\Da}{ {\cal D }}

\newcommand{\Ka}{ {\cal K }}

\newcommand{\Ea}{ {\cal E }}
\newcommand{\Sa}{ {\cal S }}

\newcommand{\Va}{ {\cal V }}

\newcommand{\Qa}{ {\cal Q }}
\newcommand{\Oa}{ {\cal O }}
\newcommand{\Ia}{ {\cal I }}
\newcommand{\Xa}{ {\cal X }}
\newcommand{\Ma}{ {\cal M }}
\newcommand{\Ta}{ {\cal T}}

\newcommand{\Ja}{ {\cal J }}
\newcommand{\Pa}{ {\cal P }}

\newcommand{\Wa}{ {\cal W }}

\newcommand{\point}{\mbox{\LARGE .}}
\newcommand{\cqfd}{\hfill\blbx \\}
\def\blbx{\hbox{\vrule height 5pt width 5pt depth 0pt}\medskip}

\def \PP{\mathbb{P}}
\def \RR{\mathbb{R}}

\def \EE{\mathbb{E}}

\def \LL{\mathbb{L}}

\def \BB{\mathbb{B}}
\newcommand{\cchi}{\protect\raisebox{2pt}{$\chi$}}

\newcommand{\vertiii}[1]{{\left\vert\kern-0.25ex\left\vert\kern-0.25ex\left\vert #1 
    \right\vert\kern-0.25ex\right\vert\kern-0.25ex\right\vert}}

\begin{document}

  \title{A Lyapunov approach to stability of positive semigroups: 
  An overview with illustrations}
  \author{Marc Arnaudon$^{1}$, Pierre Del Moral$^{2}$ \& El Maati Ouhabaz$^{1}$}
\affil{\footnotesize $^{1}$University of Bordeaux, Institut de Math\'ematiques de Bordeaux, France. {\footnotesize E-Mail:\,} \texttt{\footnotesize marc.arnaudon@math.u-bordeaux.fr, Elmaati.Ouhabaz@math.u-bordeaux.fr}}
\affil{\footnotesize $^{2}$Centre de Recherche Inria Bordeaux Sud-Ouest, Talence, 33405, FR. {\footnotesize E-Mail:\,} \texttt{\footnotesize pierre.del-moral@inria.fr}}
\maketitle
  \begin{abstract}
The stability analysis of possibly time varying positive  semigroups on non necessarily compact state spaces, including Neumann and Dirichlet boundary conditions is a notoriously difficult subject. These crucial questions arise in a variety of areas of applied mathematics, including nonlinear filtering, rare event analysis, branching processes, physics and molecular chemistry.  This article presents an overview of some recent Lyapunov-based approaches, focusing principally on practical and powerful tools for designing Lyapunov functions. These techniques include semigroup comparisons as well as conjugacy principles on non necessarily bounded manifolds with locally Lipschitz boundaries. All  the Lyapunov methodologies discussed in the article are 
illustrated in a variety of situations, ranging from conventional Markov semigroups on general state spaces to more sophisticated conditional stochastic processes possibly restricted to some non necessarily bounded domains, including locally Lipschitz and smooth hypersurface boundaries, Langevin diffusions as well as coupled harmonic oscillators.\\



\end{abstract}
  {\footnotesize
\setcounter{tocdepth}{2}
\newpage\tableofcontents
}

\newpage
\section{Introduction}

This review article outlines some of the main points of the stability theory of possibly time varying positive semigroups on non necessarily compact state spaces, including Neumann and Dirichlet boundary conditions. We present an overview of some recent Lyapunov-based approaches, focusing principally on practical and powerful tools for designing Lyapunov functions.

Foster-Lyapunov criterion dates back to the 1950s with the seminal articles~\cite{foster,harris}. These criteria are nowadays an essential tool to analyze the stability properties of Markov semigroups on general state spaces~\cite{bremaud,douc-04,hairer-mattingly,hasminski,meyn-tweedie,meyn-tweedie-2}.
There is also a vast literature on subgeometric convergence rates for Markov chains, starting with the foundational articles~\cite{nummelin,tuominen} based on sequences of Lyapunov-type functions defined in terms of some well chosen subgeometrical rate, followed by 
the control of modulated moments of the return-time to some regular set. More practical Foster-Lyapunov conditions are presented in~\cite{dfms,fort,jarner,veret-1,veret-2}.
Subgeometric convergence rates for continuous time Markov processes are also discussed in \cite{dfg-09,fort-r,hairer-mattingly}, see also the more recent article~\cite{bernou}.

 The use of Foster-Lyapunov criteria in the context of positive semigroup
arising in discrete time nonlinear filtering goes back to the pioneering articles~\cite{douc-moulines-ritov,whiteley}, based on coupling techniques developed in~\cite{klepsyna-2,klepsyna-3}. The extension of Foster-Lyapunov criterion to discrete or continuous time varying positive semigroups and their normalized versions
on general state spaces were further developed in~\cite{dhj-21}, extending Dobrushin's ergodic coefficient techniques introduced
 in~\cite{dg-ihp,dg-cras} and further developed in \cite{dm-2000,dm-sch-2,dm-miclo-ledoux,dm-04,dm-penev-2017} to unbounded state space models.

Recall that the Dobrushin's ergodic coefficient of a Markov semigroup  is the operator norm of the Markov integral operator  acting on probability measures equipped with the total variation norm (see for instance~\cite{dm-miclo-ledoux} and references therein). In the same vein, the $V$-Dobrushin's ergodic coefficient of a Markov transition introduced in~\cite{dm-penev-2017} is defined as the operator norm of the Markov integral operator acting on probability measures equipped with the $V$-norm. In this framework, the contraction w.r.t. $V$-norms is deduced by coupling
the Foster-Lyapunov criterion with a local contraction on a sufficiently large compact sub-level set of the Lyapunov function.

This operator-theoretic framework is discussed in Section~\ref{brief-review-sec} in the context of discrete time and homogeneous Markov semigroups. Section~\ref{ref-sec-discrete} is dedicated to $V$-norm contraction coefficients and the exponential convergence of Markov semigroups.  This rather elementary operator-theoretic framework is further extended in Section~\ref{subgeo-sec} to derive in a rather simple way subexponential convergence rates, stripped of all analytical superstructure, and probabilistic irrelevancies. 

The extension of this framework to more general classes of time varying Markov  semigroups with possibly continuous time indices is discussed in Section~\ref{sec-Markov-V-norm} as well as in Section~\ref{sec-Markov-lyap} in the context of diffusion semigroups.
Exponential stability theorems for more general classes of positive semigroups are discussed in Section~\ref{norm-sg-V-norm}.

To take the discussion one step further and underline the role of Lyapunov conditions, we emphasize that local contraction principles (a.k.a. local minorization conditions)
on the compact sub-level sets  
 of a prescribed Lyapunov function are generally  easily verifiable conditions. This property is often deduced from a Doeblin type local minorization property of integral operators on the
 compact sub-level sets of the Lyapunov function.
  For instance, this  local minorization condition is satisfied as soon as the semigroup is lower bounded by an absolutely continuous integral operator 
  (a.k.a.  transition kernel operator). This class of models includes hypo-elliptic diffusion semigroups as well as some regular  jump processes
  on non necessarily bounded domains. 
  
We also underline that for diffusion semigroups  with smooth densities on bounded manifolds with entrance boundaries (i.e. boundary states that cannot be reached from the inside), the existence of a sufficiently strong Lyapunov function is essential 
  to ensure the stability of the semigroup. 
  In this context, the transition densities are null on entrance boundary states so that
  the local minorization condition alone applied to some exhausting sequence of
  compact subsets is not sufficient to ensure the stability of the process.  
 The  exhausting sequence of
  compact subsets  needs to be equivalent to  the sub-level sets  
 of some sufficiently strong Lyapunov function near entrance boundaries.    For a more thorough discussion on this subject we refer to Section~\ref{brief-review-sec} and the article~\cite{dhj-21}, see also the series of Riccati-type diffusions discussed in Section~\ref{ricc-sec}.
 
 The general problem of constructing Lyapunov functions for positive semigroups, including for Markov semigroups often requires to have some good intuition about a candidate for a Lyapunov function on some particular class of model. As for deterministic dynamical systems, the design of Lyapunov functions for sub-Markov semigroups associated with a non-absorbed stochastic process requires to use some physical insight on the stability and the  behavior of the free evolution stochastic process near possible absorbing boundaries.
 
 Constructing Lyapunov functions for general classes of positive semigroups is well known as a very hard problem in system theory as well as in applied probability literature.
  The main subject of this article is to find {\em practical ways} to design these Lyapunov functions for various classes of positive semigroups that have been discussed in the literature, including conditional diffusions on manifolds with Neumann and Dirichlet boundaries. We did our best to cover the subject as broadly as possible, we also refer to the article~\cite{dhj-21} for additional historical and reference pointers. Due to the vast literature on this subject we apologize for possible omissions of some important contributions due to the lack of knowledge.

The remainder of this article is structured as follows: 

In Section~\ref{brief-review-sec}, we begin with 
a brief review on the stability of Markov semigroup.  The extension of these results to time varying positive semigroups are discussed in Section~\ref{brief-review-sec-pos}.
Section~\ref{norm-sg-V-norm} is dedicated to exponential stability theorems for normalized semigroups. In Section~\ref{time-homogenous-sec}, we present some consequences of these results in the context of time homogenous models, including  existence of ground states and quasi-invariant measures. Section~\ref{sub-Markov-sec-intro} presents different tools to design Lyapunov functions for continuous time Markov semigroups and sub-Markov semigroups. We also illustrate these results  through different examples of  semigroups arising in physics and applied probability, including overdamped Langevin diffusions, Langevin and hypo-elliptic diffusions, as well as typical examples of solvable one-dimensional sub-Markov semigroups such as the harmonic oscillator, the half-harmonic oscillator and the Dirichlet heat kernel.
General comparison and conjugacy principles to construct Lyapunov functions for positive semigroups are provided in Section~\ref{lyapunov-principles-sec}. Boundary problems are discussed in some detail in Section~\ref{boundary-sec}.
We then turn in Section~\ref{ricc-sec} to the design of Lyapunov functions for Riccati type processes, including positive definite matrix valued diffusions, logistic and multivariate birth and death processes arising respectively in Ensemble Kalman-Bucy filter theory and population dynamic analysis.

In Section~\ref{conditional-diff-sec} we illustrate the power of the Lyapunov approach in the context of
 multivariate conditional diffusions. Section~\ref{hyper-surface-bound-sec} is dedicated to illustrations with explicit computations of  geometrical objects for the Lyapunov functions discussed in Section~\ref{sec-smooth-boundary} in the context of hypersurface Dirichlet boundaries. 

 \subsection{Some basic notation}\label{notation-sec}

Let $\Ba(E)$ be the algebra  of locally bounded  measurable functions on 
a locally compact Polish space $E$.
We denote by $\Ba_b(E)\subset \Ba(E)$ the sub-algebra  of bounded measurable functions endowed with the supremum  norm $\|\point \|$.

For a given uniformly positive function $V\in \Ba(E)$, we let 
$\Ba_V(E)\subset \Ba(E)$ be the sub-space of functions $f\in \Ba(E)$ equipped with the norm
$\Vert f\Vert_V:=\Vert f/V\Vert$. 

We also let $\Ba_{\infty}(E)\subset \Ba(E)$ be the subalgebra of  locally bounded and uniformly positive functions $V$ that grow at infinity; that is, $\sup_KV<\infty$ for any compact set $K\subset E$, and for any $r\geq V_{\star}:=\inf_E V>0$  the $r$-sub-level set $\Va(r):=\{V\leq r\}\subset E$ is a non empty compact subset. 
We denote by $\Ba_0(E):=\{1/V~:~V\in \Ba_{\infty}\}\subset \Ba_b(E)$ the sub-algebra of (bounded) positive functions, locally lower bounded and that vanish at infinity.
For a given $V\in \Ba_{\infty}(E)$,  consider the subspace
$$
\Ba_{0,V}(E):=\left\{f~\in \Ba(E)~:~\vert f\vert /V\in \Ba_0(E)\right\}.
$$

 We denote by  $\Ca(E)\subset \Ba(E)$ the sub-algebra  of continuous functions and by  $\Ca_b(E)\subset \Ca(E)$ the sub-algebra  of bounded continuous functions. 
 
 We also set  $\Ca_V(E):=\Ba_V(E)\cap \Ca(E)$,  $\Ca_{0}(E):=\Ba_{0}(E)\cap \Ca(E)$ and $\Ca_{\infty}(E):=\Ba_{\infty}(E)\cap \Ca(E)$ and $\Ca_{0,V}(E):=\Ba_{0,V}(E)\cap \Ca(E)$.
Note that none of the sub-algebras $\Ba_{0}(E)$ and $\Ba_{\infty}(E)$ have an unit unless $E$ is compact, the null function $0\not\in \Ba_0(E)$ but the unit function $1\in \Ca_{0,V}(E)$ as soon as $V\in\Ba_{\infty}(E)$. 

Let $\Ma_b(E)$ be the set of bounded signed measures $\mu$ on $E$ equipped with the total variation norm $\Vert\mu\Vert_{tv}:=\vert \mu\vert(E)/2$,  where $
\vert\mu \vert:=\mu_++\mu_-$ stands for  the total variation measure associated with  a Hahn-Jordan decomposition $\mu=\mu_+-\mu_-$ of the measure. 
Also let $\Pa(E)\subset \Ma_b(E)$ be the subset of probability measures on $E$.  Recall that for any $\mu_1,\mu_2\in \Pa(E)$ and $\epsilon\in ]0,1]$ we have
\begin{equation}\label{ref-coupling-tv}
\Vert \mu_1-\mu_2\Vert_{\tiny tv}\leq 1-\epsilon\Longleftrightarrow\left( \exists \nu\in \Pa(E)~:~\mu_1\geq \epsilon~\nu\quad\mbox{\rm and}\quad \mu_2\geq \epsilon~\nu\right).
\end{equation}

With a slight abuse of notation, we denote by $0$ and $1$ the null and unit scalars as well as the null and unit function on $E$.

 Let  $Q_{s,t}$,   be a semigroup of positive integral operators on $\Ba_b(E)$ indexed by  continuous time indices $s, t\in \Ta=\RR_+:=[0,\infty[$ or by a discrete time index set
 $\Ta=\NN$, with $s\leq t$. 
The action of $Q_{s,t}$ on $\Ba_b(E)$ is given for any $f\in \Ba_b(E)$  by the formulae
\begin{equation}\label{right-actions}
Q_{s,t}(f)(x):=\int Q_{s,t}(x,dy)~f(y).
\end{equation}

The  left action of $Q_{s,t}$ on $\Ma_b(E)$  is given for any  $\eta\in \Ma_b(E)$ by the formulae
\begin{equation}\label{right-actions}
(\eta\, Q_{s,t})(dy):=\int\eta(dx)~Q_{s,t}(x,dy).
\end{equation}

In this notation, the semigroup property takes the following form
\begin{equation}\label{def-Q-s-t-intro}
Q_{s,u} Q_{u,t}=Q_{s,t}\quad \mbox{\rm with}\quad Q_{s,s}=I,\quad \mbox{\rm the identity operator}.
\end{equation}
In the above display, $Q_{s,u} Q_{u,t}$ is a shorthand notation for the composition $Q_{s,u} \circ Q_{u,t}$ of the left or right-action operators. Unless otherwise stated, all the semigroups discussed in this article are  indexed by conformal indices $s\leq t$ in the set $\Ta$. To avoid repetition, we often write $Q_{s,t}$ 
without specifying the order $s\leq t$ of the indices $s,t\in\Ta$.

We  denote by $
 \Ma_V(E)
$ be the space of measures $\mu\in  \Ma_b(E)$   equipped with the operator $V$-norm
$\vertiii{\mu}_{V}:=\vert \mu\vert(V)$, and  by $
\Pa_V(E)\subset \Ma_V(E)
$ be the convex set of probability measures. Whenever $V\geq 1$, for any $\rho>0$ we have 
the norm equivalence formulae 
\begin{equation}\label{norm-equivalence}
\rho \vertiii{\mu}_V\leq  \vertiii{\mu}_{1+\rho V}\leq (1+\rho)\vertiii{\mu}_{V}
\end{equation}

We associate with a function $h\in \Ba_{0,V}(E)$  the Boltzmann-Gibbs transformation
 \begin{equation}\label{def-Psi-H}
 \Psi_{h}~:~\mu\in \Pa_{V}(E)\mapsto \Psi_{h}(\mu)\in \Pa_{V^h}(E)
 \end{equation}
with the probability measure
$$
 \Psi_{h}(\mu)(dx):=\frac{h(x)}{\mu(h)}~\mu(dx) \quad \mbox{and}\quad
V^h:=V/h\in \Ba_{\infty}(E).
$$

We also denote by   $\vertiii{Q}_{V}$ the operator norm of  a bounded linear operator $Q:f\in \Ba_{V}(E)\mapsto Q(f)\in \Ba_{V}(E)$; that is
\begin{equation}\label{def-op-norm}
\vertiii{Q}_{V}:=\sup\{\Vert Q(f)\Vert_V~:~f\in  \Ba_V(E)\quad\mbox{\rm such that}\quad \Vert f\Vert_V\leq 1\}.
\end{equation}
In terms of the $V$-conjugate semigroup 
$$ 
f\in \Ba_b(E)\mapsto Q^V(f):=Q(Vf)/V\in \Ba_b(E)
$$
we have
$$
\vertiii{Q}_{V}=\Vert Q^V(1)\Vert=\vertiii{Q^V}:=\sup\{\Vert Q^V(f)\Vert~:~f\in  \Ba_b(E)\quad\mbox{\rm such that}\quad \Vert f\Vert\leq 1\}.
$$
For a given measurable function $f$ and a given measurable subset, we use the shorthand notation
$$
-\infty\leq \inf_A f:=\inf_{x\in A}f(x)\leq \sup_A f:=\sup_{x\in A}f(x)\leq+\infty.
$$
For a given  $s\in \Ta$ and $\tau\in \Ta$ with $\tau>0$, we consider the time mesh
$$
 [s,\infty[_{\tau}:=\{s+n\tau \in [s,\infty[~:~n\in\NN\}.
$$

Throughout, unless otherwise is stated we write $c$ for some positive constants whose values may vary from line to line, and we write $c_{\alpha}$, as well as $c(\beta)$ and $c_{\alpha}(\beta)$ when their values may depend on  some parameters $\alpha,\beta$ defined on some parameter sets. We also set $a\wedge b=\min(a,b)$, $a\vee b=\max(a,b)$, and $a_+=a\vee 0$ for $a,b\in \RR$.

 \subsection{$V$-positive semigroups}\label{reg-cond-sec}

 We say that $Q_{s,t}$ is a
 $V$-positive semigroup on $\Ba_V(E)$ for some Lyapunov function $V\in \Ba_{\infty}(E)$  as soon as there exists some $\tau>0$ and some  function
 $\Theta_{\tau}\in \Ba_0(E)$ such that for any  $0<f\in\Ba_V(E)$ and $ s<t$ we have
 $ 0<Q_{s,t}(f)\in \Ba_{0,V}(E)$ as well as
 \begin{equation}\label{ref-V-theta-intro}
 Q_{s,s+\tau}(V)/V\leq \Theta_{\tau}\quad \mbox{and}\quad  \sup_{\vert t-s\vert\leq \tau}\left(\vertiii{Q_{s,t}}\vee\vertiii{Q_{s,t}}_V\right)<\infty. 
\end{equation}

As shown in Section~\ref{foster-lyap}, the l.h.s. criterion in (\ref{ref-V-theta-intro}) can be seen as a uniform Foster-Lyapunov condition (a.k.a. drift condition).

The irreducibility condition $f>0\Longrightarrow Q_{s,t}(f)>0$ is satisfied if and only if we have  $Q_{s,t}(1)>0$.  We check this claim by contradiction. Assume that  $Q_{s,t}(1)>0$ and consider a function $f>0$ and some $x\in E$ such that  $Q_{s,t}(f)(x)=0$. In this case, for any $\epsilon>0$ we would have $$
\epsilon~Q_{s,t}\left(1_{f\geq \epsilon}\right)(x)\leq Q_{s,t}(f)(x)=0
$$  by Fatou's lemma  we would find the contradiction
$$
\liminf_{\epsilon\rightarrow 0}Q_{s,t}\left(1_{f\geq \epsilon})\right(x)=0\geq Q_{s,t}(1)(x)\Longrightarrow Q_{s,t}(1)(x)=0.
$$
Without further mention, all semigroups $Q_{s,t}$ considered in this article are assumed to be semigroups of positive integral operators $Q_{s,t}$ on $\Ba_b(E)$ satisfying the irreducibility condition $Q_{s,t}(1)>0$ for any $s\leq t$. 
Notice  that the condition 
$$
0<f\in\Ba_V(E)\quad\Longrightarrow\quad  \forall s<t\qquad 0<Q_{s,t}(f)\in \Ba_{0,V}(E)
$$
is met as soon as $Q_{s,t}$ is a strong $V$-Feller semigroup (i.e. for any $s<t$ we have $Q_{s,t}(\Ba_V(E))\subset\Ca_{V}(E)$ and when we have $Q_{s,t}(V)/V\in \Ba_0(E)$). To check this claim, observe that for any positive function $f\in \Ba_V(E)$ and $s<t$ the function $Q_{s,t}(f)$ is positive and  continuous;  and thus locally lower bounded. In this situation, whenever $\Vert f\Vert_V\leq 1$, for any $s<t$ we have the comparison property
$$
Q_{s,t}(f) /V\leq Q_{s,t}(V)/V \in \Ba_0(E)\Longrightarrow Q_{s,t}(f)/V\in \Ba_{0}(E)\Longleftrightarrow Q_{s,t}(f)\in \Ca_{0,V}(E).
$$
In summary,  a strong $V$-Feller semigroup $Q_{s,t}$ is $V$-positive on $\Ba_V(E)$ as soon as  there exists some $\tau>0$ and some  function
 $\Theta_{\tau}\in \Ba_0(E)$ such that the l.h.s. condition in (\ref{ref-V-theta-intro}) is met and for any $s<t$ we have
$$
Q_{s,t}(V)/V\in \Ba_0(E)\quad \mbox{and}\quad Q_{s,s+\tau}(V)/V\leq \Theta_{\tau}\in \Ba_0(E). 
$$

When $V\in\Ca_{\infty}(E)$, we say that  $Q_{s,t}$ is a
 $V$-positive semigroup on $\Ca_V(E)$  as soon as $Q_{s,t}(\Ca_V(E))\subset \Ca_{0,V}(E)$ for any $s<t$ and  condition  (\ref{ref-V-theta-intro}) is met.
 
A $V$-Feller semigroup $Q_{s,t}$ for some $V\in\Ca_{\infty}(E)$, in the sense that for any $s<t$ we have $Q_{s,t}(\Ca_V(E))\subset\Ca_{V}(E)$, is also said to be $V$-positive on $\Ca_V(E)$
as soon as  there exists some $\tau>0$ and some  function
 $\Theta_{\tau}\in \Ba_0(E)$ such that the l.h.s. condition in (\ref{ref-V-theta-intro}) is met  and for any $s<t$ we have
$$
Q_{s,t}(V)/V\in \Ca_0(E)\quad \mbox{and}\quad Q_{s,s+\tau}(V)/V\leq \Theta_{\tau}\in \Ba_0(E). 
$$

Last but not least, observe that positive semigroups $\Qa_{s,t}$ with continuous time indices $s\leq t\in \RR_+$ can be turned into discrete time models by setting $Q_{p,n}=\Qa_{p\tau,n\tau}$ for any $p\leq n\in \NN$ and some parameter $\tau>0$. Up to a time rescaling, the parameter $\tau>0$ arising in the definition of a discrete time $V$-positive semigroups $Q_{p,n}$ can be chosen as the unit time parameter. In this context, the r.h.s. condition in (\ref{ref-V-theta-intro}) is automatically satisfied.

\section{A brief review on Markov semigroups}\label{brief-review-sec}

The stability analysis of positive semigroups presented in this article is mainly based on discrete time operator-type contraction techniques combining Lyapunov inequalities with local minorization conditions.  This section presents a brief overview of this 
operator-theoretic framework.
Our presentation is nearly self-contained and follows that of Section 8 in the book~\cite{dm-penev-2017}  (see also Section 2 in~\cite{dhj-21}).

\subsection{$V$-norm contraction coefficients}\label{ref-sec-discrete}

In this section we are mainly interested in the contraction properties of discrete time Markov integral operators. We only consider time homogeneous Markov semigroups $P:=P_{t,t+1}$ on $\Ba_b(E)$, so that $P_t:=P_{0,t}=P P_{1,t}$.  
In what follows $E$ is assumed to be a Polish space. One key mathematical object is the $V$-norm contraction coefficient.

We further assume that there exists a Lyapunov function $V\in \Ba_{\infty}(E)$ and parameters $\epsilon \in ]0,1[$ and $c<\infty$ such that
\begin{equation}\label{PVC}
 P(V)\leq \epsilon~V+c.
 \end{equation}
 
Note that the class of Markov semigroups considered in this section is more general than the one discussed in (\ref{ref-V-theta-intro}). Indeed, in our context the Lyapunov condition stated in the l.h.s. of $(\ref{ref-V-theta-intro})$ takes the form
 $$
 P(V)/V\leq \Theta\in \Ba_0(E).
 $$
 This condition ensures that {\em for any} $0<\epsilon< \Vert\Theta\Vert$, the set $K_{\epsilon}:=\{\Theta\geq \epsilon\}$ is a non empty compact subset and we have
 $$
 P(V)\leq \epsilon~ 1_{E-K_{\epsilon}} V+
 1_{K_{\epsilon}}(\Theta~V)\leq  \epsilon~ V+c_{\epsilon}\quad \mbox{\rm with}\quad c_{\epsilon}:=\Vert\Theta\Vert
 \sup_{K_{\epsilon}}V.
 $$

Replacing $V$ by $2^{-1}(1+\epsilon V/c)$ there is no loss of generality to assume that $c=1/2$ and $V\geq 1/2$.
Also assume there exists some  $r_0\geq 1$  and some function $\alpha : r\in [r_0,\infty[\,\mapsto\alpha(r)\in \, ]0,1]$,  such that for any $r\geq r_0$ we have
\begin{equation}\label{loc-min}
\sup_{(x,y)\in \Va(r)^2}\Vert\delta_xP-\delta_yP\Vert_{\tiny tv}\leq 1-\alpha(r)\quad \mbox{\rm with}\quad \Va(r):=\{V\leq r\}.
 \end{equation}
The $V$-Dobrushin coefficient $\beta_{V}(P)$ of  $P$  is defined by the  $V$-norm  operator 
\begin{equation}\label{defi-beta-V}
 \beta_{V}(P):=\sup_{\mu,\eta\in \Pa_{V}(E)}{\vertiii{(\mu-\eta)P}_{V}}/{\vertiii{ \mu-\eta}_{V}}.
\end{equation}
As show in Section 8 in~\cite{dm-penev-2017} (see also Section 2.3 in~\cite{dhj-21}), the supremum in (\ref{defi-beta-V}) is attained on Dirac masses $(\mu,\eta)=(\delta_x,\delta_y)$; that is,  we have
$$
 \beta_{V}(P)=\sup_{(x,y)\in E^2}\frac{\Vert \delta_xP-\delta_yP\Vert_V}{V(x)+V(y)}.
$$ 
The terminology $V$-Dobrushin coefficient comes from the fact that we recover the standard Dobrushin coefficient $\beta(P):= \beta_{1/2}(P)$ by choosing the constant function $V=1/2$.
Theorem 8.2.21 in~\cite{dm-penev-2017}  (see also Lemma 2.3 in~\cite{dhj-21}) shows that the Lyapunov inequality (\ref{PVC}) combined with the local minorization condition (\ref{loc-min}) yield a $V$-norm contraction estimate for some well chosen Lyapunov function. 

\begin{lem}[\cite{dm-penev-2017}]\label{ref-beta-V-est}
Assume (\ref{PVC}) and (\ref{loc-min}). In this situation, for any $r\geq r_0\vee r_{\epsilon}$  with $r_{\epsilon}:=1/(1-\epsilon)$ we have
\begin{equation}\label{betaV-est}
  \beta_{V_{\epsilon,r}}(P)\leq 
1-\alpha_{\epsilon}(r)
 \end{equation}
with the rescaled Lyapunov function
$$
V_{\epsilon,r}:=\frac{1}{2}\left(1+\frac{1}{1+\epsilon}~\frac{\alpha(r)}{r}~V\right)\quad \mbox{and}\quad
\alpha_{\epsilon}(r):=\frac{\alpha(r)}{2}~\frac{(1-\epsilon)}{(1+\epsilon)+ \frac{\alpha(r)}{2}}\left(1-\frac{r_{\epsilon}}{r}\right)>0.
$$
\end{lem}
For the convenience of the reader a proof of the $V$-norm estimate (\ref{betaV-est}) is provided in the Appendix, on page~\pageref{betaV-est-proof}.
As an aside, note that whenever (\ref{PVC}) is met with $c=1/2$ we have
$$
P(V_{\epsilon,r})\leq \epsilon ~V_{\epsilon,r}+c_{\epsilon,r}\quad \mbox{\rm with}\quad 
c_{\epsilon,r}:=\frac{1}{2}\left((1-\epsilon)+
\frac{1}{1+\epsilon}~\frac{\alpha(r)}{2r}\right)
$$

The equivalence of the $V$-norm and the $V_{\epsilon,r}$-norm yields without further work the following contraction theorem.
\begin{theo}\label{theo-intro-Markov-hom} 
For any $t\in\NN$ and any $\mu,\eta\in \Pa_V(E)$ we have
\begin{equation}\label{lipschitz-inq-M-hom}
\vertiii{(\mu-\eta) P_t}_{V}\leq c_{\epsilon,r}~\beta_{V_{\epsilon,r}}(P)^t~\vertiii{\mu-\eta}_{V}\quad
\mbox{with}\quad c_{\epsilon,r}:=1+2r(1+\epsilon)/{\alpha(r)}.
\end{equation}
\end{theo}
The contraction estimate (\ref{lipschitz-inq-M-hom}) ensures the existence of a single invariant probability measure $\mu_{\infty}=\mu_{\infty}P_{t}\in \Pa_V(E)$. Similar approaches are presented in the article~\cite{hairer-mattingly}, simplifying the Foster-Lyapunov methodologies and the small-sets return times estimation techniques developed in~\cite{meyn-tweedie}. While the geometric convergence of discrete time Markov semigroups towards their invariant measure $\mu_{\infty}$ (a.k.a. Harris-type theorem) are well known, the $V$-norm contraction coefficients techniques developed in~\cite{dm-penev-2017} provides a very direct and short proof of Theorem~\ref{theo-intro-Markov-hom}. In Section~\ref{sec-Markov-V-norm}, this natural operator-theoretic framework is extended without difficulties 
to time varying and continuous time indices. 
For instance, for any collection of Lyapunov functions $V_t\in \Ba_{\infty}(E)$  and any Markov transitions $M_t$ indexed by $t\in\NN$ from 
$\Ba_{V_t}(E)$ into $\Ba_{V_{t-1}}(E)$ we have
$$
P_t:=P_{t-1}M_t=M_1\ldots M_t\Longrightarrow
\vertiii{(\mu-\eta)P_t}_{V_t}\leq \left(\prod_{1\leq s\leq t} \beta_{V_{s-1},V_s}(M_s)\right)~\vertiii{\mu-\eta}_{V_0}
$$
with the $V$-norm contraction coefficients
$$
 \beta_{V_{t-1},V_t}(M_t):=\sup_{\mu,\eta\in \Pa_{V}(E)}{\vertiii{(\mu-\eta)M_t}_{V_{t}}}/{\vertiii{ \mu-\eta}_{V_{t-1}}}=\sup_{(x,y)\in E^2}\frac{\Vert \delta_xM_t-\delta_yM_t\Vert_{V_t}}{V_{t-1}(x)+V_{t-1}(y)}.
$$
Whenever all the Markov transitions $M_t$ satisfy (\ref{PVC}) anf (\ref{loc-min}) with the same  Lyapunov function $V$ and the same function $\alpha(r)$, choosing the function $V_t:=V_{\epsilon,r}$ defined in Lemma~\ref{ref-beta-V-est} we have $ \beta_{V_{t-1},V_t}(M_t)\leq (1-\alpha_{\epsilon}(r))$, with the parameter $\alpha_{\epsilon}(r)$ defined in Lemma~\ref{ref-beta-V-est}. More general time-inhomogeneous models can probably be handle extending the analysis developed in~\cite{dg-ihp,dm-2000} as well as in Theorem 4.18 and Theorem 4.20 in~\cite{seneta} and in Proposition 1 in~\cite{witheley2} using the standard Dobrushin coefficient to non necessarily compact spaces in terms of $V$-norms contraction coefficients.
\begin{rmk}\label{rmk-unif-tv}
Whenever $P(V)/V\leq \Theta$ for some $\Theta\in \Ba_0(E)$ a direct application of (\ref{lipschitz-inq-M-hom}) yields for any $t\geq 1$ the estimate
$$
\vertiii{(\mu-\eta) P_{t+1}}_{V}\leq c_{\epsilon,r}~\beta_{V_{\epsilon,r}}(P)^t~\vertiii{\mu-\eta}_{V_{\Theta}}\quad \mbox{\rm with}\quad V_{\Theta}:=V~\Theta 
$$
In some situations (cf. for instance Section~\ref{ricc-sec}), we can choose $\Theta=1/V$. In this context, for any any $\mu,\eta\in \Pa(E)$ we have the  uniform total variation norm estimates
$$
2~\vertiii{(\mu-\eta) P_{t+1}}_{\tiny tv}\leq \vertiii{(\mu-\eta) P_{t+1}}_{V}\leq 2c_{\epsilon,r}~\beta_{V_{\epsilon,r}}(P)^t~\vertiii{\mu-\eta}_{\tiny tv}\leq 2c_{\epsilon,r}~\beta_{V_{\epsilon,r}}(P)^t~
$$
In this situation, besides the state space $E$ may not be compact, for sufficiently large $t>0$ the standard Dobrushin contraction coefficient $\beta(P_t)<1$ of the Markov transition $P_t$ yields  exponential decays.
\end{rmk}

The extension to more general positive semigroups is slightly 
more involved and relies on the stability properties of triangular arrays of Markov operators (cf.~\cite{dhj-21}). In Section~\ref{brief-review-sec-pos} we present an overview on the stability properties of this class of models.

Combining (\ref{ref-coupling-tv}) with (\ref{loc-min}), for any $r\geq r_0$ there exists some probability measure $\nu_r$ such that for any $\mu=(\mu_+-\mu_-)\in \Ma_V(E)$ with $\mu(1)=0$ and any bounded positive function $f\geq 0$ we have
$$
\mu_{\pm}(P(f))\geq\mu_{\pm}(1_{\Va(r)})~\alpha(r)~\nu_r(f).$$
This implies that
$$
\frac{\Vert\mu P\Vert_{\tiny tv}}{\Vert \mu\Vert_{\tiny tv}}=\left\Vert\left(\frac{\mu_+}{\mu_+(1)}-\frac{\mu_-}{\mu_-(1)}\right)P\right\Vert_{\tiny}\leq 1-\alpha(r)\left(\frac{\mu_+(1_{\Va(r)})}{\mu_+(1)}\wedge\frac{\mu_-(1_{\Va(r)})}{\mu_-(1)}\right).
$$
On the other hand, by Markov inequality we have
$$
1-\frac{1}{r}\left(\frac{\mu_+(V)}{\mu_+(1)}\vee\frac{\mu_-(V)}{\mu_-(1)}\right)\leq \frac{\mu_+(1_{\Va(r)})}{\mu_+(1)}\wedge\frac{\mu_-(1_{\Va(r)})}{\mu_-(1)}.
$$
This implies that
$$
\frac{\Vert\mu P\Vert_{\tiny tv}}{\Vert \mu\Vert_{\tiny tv}}\leq 
1-\alpha(r)\left(1-\frac{1}{r}~\frac{\vertiii{\mu}_V}{\Vert\mu\Vert_{\tiny tv}}\right).
$$
We summarize the above discussion in the following lemma.
\begin{lem}\label{lem-alpha-min}
Assume that (\ref{loc-min}) is met. In this case for any $r\geq r_0$ and $\mu\in \Ma_V(E)$ with $\mu(1)=0$ we have
\begin{equation}
\Vert\mu P\Vert_{\tiny tv}\leq (1-\alpha(r))~\Vert \mu\Vert_{\tiny tv}+\frac{\alpha(r)}{r}~\vertiii{\mu}_V.
\end{equation}
\end{lem}
Whenever (\ref{PVC}) is met with $c=1/2$,  for any $\mu\in \Ma_V(E)$ with $\mu(1)=0$, recalling that $\Vert \mu\Vert_{\tiny tv}=\vertiii{\mu }_{1/2}=\mu_+(1)=\mu_-(1)$  we have
$$
\vertiii{\mu P}_V\leq \epsilon~\vertiii{\mu}_V+\Vert \mu\Vert_{\tiny tv}.
 $$
 This yields for any $\rho>0$ the estimate
 $$
 \vertiii{\mu P}_{1/2+\rho V}\leq 
 ((1-\alpha(r))+\rho)~\Vert \mu\Vert_{\tiny tv}+\left(\frac{\alpha(r)}{\rho r}+ \epsilon\right)~\rho ~\vertiii{\mu}_V
 $$
 from which we readily check the following lemma.
 \begin{lem} 
 Assume  (\ref{PVC}) and (\ref{loc-min})  are met with $c=1/2$. In this situation,
 for any  $(r,\rho)$ such that ${r_{\epsilon}}/{r}<\rho/\alpha(r)<1$ we have
 $$
\beta_{1/2+\rho V}(P)\leq 1-
 \alpha^{\rho}_{\epsilon}(r) $$
 with the parameter $r_{\epsilon}$ defined in Lemma~\ref{ref-beta-V-est} and
 $ \alpha^{\rho}_{\epsilon}(r)$ defined by
 $$
 \alpha^{\rho}_{\epsilon}(r):=\left( \alpha(r)-\rho\right)\wedge \left((1-\epsilon)\left(1-\frac{\alpha(r)}{\rho}\frac{r_{\epsilon}}{r}\right)\right)>0.
 $$
 \end{lem}
 For instance, for any $\delta>0$ we have
 $$
 \begin{array}{l}
\displaystyle r>(1+\delta)~r_{\epsilon}\quad \mbox{\rm and}\quad
 \rho:=(1+\delta)~\alpha(r)~\frac{r_{\epsilon}}{r}\\
 \\
\displaystyle \Longrightarrow
  \alpha^{\rho}_{\epsilon}(r)\geq \left(\alpha(r)~\left(1-(1+\delta)~\frac{r_{\epsilon}}{r}\right)\right)\wedge\left((1-\epsilon)~\frac{\delta}{1+\delta}\right)>0.
 \end{array}$$
 
 We end this section with a weaker version of another popular condition ensuring the non expansive property of Markov semigroups w.r.t. $V$-norms (see for instance~\cite{andrieu,connor,dfms,dfg-09} and references therein). 
 Instead of (\ref{PVC}), we further assume there exist an increasing function $\varphi~:~v\in [1,\infty[\mapsto \varphi(v)\in [0,\infty[$  some constant $c\geq 0$ and some function $V\in \Ba_{\infty}(E)$ such that $V\geq 1$ and
\begin{equation}\label{P-varphi}
P(V)\leq V-\varphi(V)+c\quad \mbox{\rm and}\quad \varphi(V)/V\in \Ba_0(E).
\end{equation}
Note that when $\varphi(V)=(1-\epsilon)V$ the r.h.s. condition in the above display is not met (unless $E$ is compact) but the l.h.s. inequality coincides with the Lyapunov condition (\ref{PVC}). In contrast with the Lyapunov condition introduced in~\cite{dfms}, the function $\varphi$ is not required to be concave.

 Whenever (\ref{loc-min}) is satisfied, Lemma~\ref{lem-alpha-min}
ensures that for any $r\geq r_1:=\varphi(r_0)$ and $\mu\in \Ma_{\varphi(V)}(E)$ with $\mu(1)=0$ we have
\begin{equation}\label{ref-r1-alpha1}
\Vert\mu P\Vert_{\tiny tv}\leq (1-\alpha_{1}(r))~\Vert \mu\Vert_{\tiny tv}+\frac{\alpha_{1}(r)}{r}~\vertiii{\mu}_{\varphi(V)}\quad\mbox{\rm with}\quad
 \alpha_{1}(r):=\alpha\left(\varphi^{-1}(r)\right).
\end{equation}
On the other hand, by (\ref{P-varphi}) $\mu\in \Ma_{V}(E)$ with $\mu(1)=0$  we have
$$
\vertiii{\mu P}_V\leq \vertiii{\mu}_V-\vertiii{\mu}_{\varphi(V)}+2c\Vert \mu\Vert_{\tiny tv}.
$$
This yields for any $\rho>0$ the inequality
\begin{equation}\label{ref-lyp-est}
\vertiii{\mu P}_{1+\rho V}\leq \vertiii{\mu}_{1+\rho V}-\left(
2\left(\alpha_{1}(r)-\rho c\right)~\Vert \mu\Vert_{\tiny tv}+\left(\rho-
\frac{2\alpha_{1}(r)}{r}\right)~\vertiii{\mu}_{\varphi(V)}\right)
\end{equation}
from which we readily check the following uniform estimates.
\begin{lem}
 Assume (\ref{loc-min})  and (\ref{P-varphi}) are met. In this situation,   (\ref{ref-r1-alpha1}) is met 
for some parameters $(r_1,\alpha_1)$. In addition,
for any $r\geq r_1$  and parameter $\rho >0$  such that
$
\rho c\leq \alpha_{1}(r)$ and $\rho \geq {2\alpha_{1}(r)}/{r}
$
we have
\begin{equation}\label{ref-unif-est}
\sup_{t\geq 0}\beta_{1+\rho V}(P_t)\leq 1\quad \mbox{or equivalently}\quad
\sup_{t\geq 0}\vertiii{\mu P_t}_{1+\rho V}\leq \vertiii{\mu}_{1+\rho V}
\end{equation}
for any $\mu\in \Ma_{V}(E)$ with $\mu(1)=0$.
\end{lem}
\subsection{Subgeometric convergence}\label{subgeo-sec}

Consider the discrete time and homogeneous Markov semigroups discussed in Section~\ref{ref-sec-discrete}. Without further mention, we shall assume that  conditions (\ref{loc-min})  and (\ref{P-varphi}) are met for some increasing function $\varphi$ and some Lyapunov function $1\leq V\in\Ba_{\infty}(E)$.

Consider a concave increasing differentiable  function $\varphi_1~:~[1,\infty[\mapsto [1,\infty[$ such that  $\varphi_1(V)/V\in \Ba_0(E)$. We further assume there exists some parameters $\chi>0$ and $\kappa_2>0$ such that
 \begin{equation}\label{ref-hyp-varphi}
\partial \varphi_1(V)~\varphi(V)\geq \varphi_2(V)\quad \mbox{\rm with}\quad \varphi_2(V):=\kappa_2~V
\left(\varphi_1(V)/V\right)^{1+\chi}\in \Ba_{\infty}(E)
\end{equation}
Observe that 
$$
\kappa_2~\varphi_2(V)/V\leq\partial \varphi_1(V)~\varphi(V)/V\leq \partial \varphi_1(1)~\varphi(V)/V\in \Ba_0(E)
$$
Thus, whenever the function $\varphi_1(V)/V$ is locally lower bounded and upper semi-continuous the above estimate ensures that
$$
\varphi(V)/V\in \Ba_0(E)\Longrightarrow \varphi_1(V)/V, \varphi_2(V)/V\in \Ba_0(E)
$$
The prototype of model we have in mind is the case
\begin{equation}\label{varphi-proto}
\varphi(v):=\kappa_0~v^{\delta}\quad\mbox{\rm and}\quad
\varphi_1(v):=\kappa_1~v^{1-\upsilon\delta}
\end{equation}
for some parameters $\upsilon,\delta\in ]0,1[$ and $\kappa_1\geq 1$. In this context, we have
$$
\partial \varphi_1(V)~\varphi(V)=\kappa_0\kappa_1(1-\upsilon\delta)~V^{\delta(1-\upsilon)}= \varphi_2(V):=\kappa_2~V~\left(\varphi_1(V)/V\right) ^{1+\chi}\in \Ba_{\infty}(E)
$$
with the parameters
$$
 \kappa_2:=\kappa_0\kappa_1^{-\frac{1-\delta}{\upsilon\delta}}(1-\upsilon\delta)\quad
 \mbox{\rm and}\quad\chi:=\frac{1-\delta}{\upsilon\delta}.
$$
Applying Jensen's inequality and using (\ref{P-varphi}) we prove that
\begin{equation}\label{varphi-P2}
P(\varphi_1(V))\leq   \varphi_1(V)-\varphi_2(V)+c_1\quad
\mbox{\rm with}\quad c_1=c\partial\varphi_1(1).
\end{equation}
We check this estimate using the fact that $\partial\varphi_1$ is decreasing. Thus, for any $0\leq u\leq v$ we have
$$
\varphi_1(v-u)\leq \varphi_1(v)-\partial \varphi_1(v)~u.
$$
In the same vein, for any $0\leq  v\leq v-u$ we have
$$
\forall w\in [v,v-u]\quad    \partial\varphi_1(v)\geq \partial\varphi_1(w)\quad \mbox{\rm and therefore}\quad
\varphi_1(v-u)\leq \varphi_1(v)-\partial \varphi_1(v)~u.
$$
The Lyapunov inequality (\ref{varphi-P2}) applied to (\ref{varphi-proto}) is closely related to 
 Lemma 3.5 in~\cite{jarner}
We further assume there exists some  $r_2\geq 1$  and some function $\alpha_2 : r\in [r_2,\infty[\,\mapsto\alpha_2(r)\in \, ]0,1]$,  such that for any $r\geq r_2$ we have
\begin{equation}\label{loc-min-2}
\sup_{(x,y)\in \{\varphi_2(V)\leq r\}^2}\Vert\delta_xP-\delta_yP\Vert_{\tiny tv}\leq 1-\alpha_2(r).
 \end{equation}
The above condition is automatically met with $r_2=\varphi_2(r_0)$ and  $\alpha_{2}(r):=\alpha\left(\varphi^{-1}_2(r)\right)$ as soon as (\ref{loc-min}) is satisfied and $\varphi_2$ is increasing.
In this situation, arguing as above, Lemma~\ref{lem-alpha-min}
ensures that for any $r\geq r_2$ and $\mu\in \Ma_{\varphi_2(V)}(E)$ with $\mu(1)=0$ we have
$$
\Vert\mu P\Vert_{\tiny tv}\leq (1-\alpha_{2}(r))~\Vert \mu\Vert_{\tiny tv}+\frac{\alpha_{2}(r)}{r}~\vertiii{\mu}_{\varphi_2(V)}.
$$
On the other hand, by (\ref{varphi-P2}) we have
$$
\vertiii{\mu P}_{\varphi_1(V)}\leq \vertiii{\mu}_{\varphi_1(V)}-\vertiii{\mu}_{\varphi_2(V)}+2c_1\Vert \mu\Vert_{\tiny tv}.
$$
Applying Jensen's inequality for any $\mu\in \Ma_{V}(E)$ with $\mu(1)=0$  we find that
\begin{equation}\label{ref-2-Jensen}
\vertiii{\mu}_{\varphi_2(V)}=\kappa_2~\frac{\vert\mu\vert\left( V~\left(\varphi_1(V)/V\right)^{1+\chi}\right)~}{\vert \mu\vert(V)}~\vert \mu\vert(V)\geq \kappa_2
~\vertiii{\mu}_{\varphi_1(V)}^{1+\chi}/\vertiii{\mu}_V^{\chi}
\end{equation}
Following word-for-word the proof of (\ref{ref-lyp-est}) we readily check the following lemma.
\begin{lem}\label{ref-lem-sub-geo}
For any $r\geq r_2\vee r_1$ and $\rho >0$ such that
$$
\rho (c_1\wedge c)\leq \alpha_{1}(r)\wedge \alpha_2(r)\quad \mbox{and}\quad \delta_{\rho}(r):=\kappa_2\left(\rho-2~\frac{\alpha_{2}(r)\vee \alpha_1(r)}{r}\right)>0
$$
and for any $\mu\in \Ma_{V}(E)$ with $\mu(1)=0$ we have
\begin{eqnarray*}
\vertiii{\mu P}_{1+\rho \varphi_1(V)}
&\leq& \vertiii{\mu}_{1+\rho \varphi_1(V)}-\delta_{\rho}(r)
~\vertiii{\mu}_{\varphi_1(V)}^{1+\chi}/\vertiii{\mu}_V^{\chi}.
\end{eqnarray*}
\end{lem}
To take the final step, recall that
$
 \vertiii{\mu}_{1+\rho \varphi_1(V)}\leq (1+\rho)\vertiii{\mu}_{\varphi_1(V)}
$. Combining this estimate with the uniform estimate (\ref{ref-unif-est}) for any $t\in\Ta=\NN$ we check that
$$
\vertiii{\mu P_{t+1}}_{1+\rho \varphi_1(V)}- \vertiii{\mu P_t}_{1+\rho  \varphi_1(V)}\leq -\varsigma_{\mu}\left(\vertiii{\mu P_t}_{1+\rho\varphi_1(V)}\right)
$$
with the function
$
\varsigma_{\mu}(u):=\omega_{\mu}(\rho)
~u^{1+\chi}$
and the parameters
\begin{equation}\label{def-omega-rho}
\omega_{\mu}(\rho):=\omega(\rho)~\vertiii{\mu}_{1+\rho V}^{-\chi}
\quad
\mbox{\rm and}\quad
\omega(\rho):=\rho^{\chi}~\delta_{\rho}(r)/(1+\rho)^{1+\chi}
\end{equation}

\begin{lem}[\cite{butkovsky}]\label{lem-butkovsky}
Consider a decreasing sequence of positive numbers $u_t$ such that for any $t\in\NN$ we have
$$
u_{t+1}-u_t\leq -\varsigma(u_t)
$$
for some continuous increasing function $\varsigma$ from $]0,u_0]$ into $]0,\infty[$. In this situation, for any $t\geq 1$ we have
$$
u_t\leq I_{\varsigma}^{-1}(t)\quad \mbox{with}\quad
I_{\varsigma}(u):=\int_u^{u_0}\frac{dv}{\varsigma(v)}.
$$
\end{lem}
\proof
Since  $\varsigma$ is increasing we have $v\in [u_{n+1},u_n]\Longrightarrow \varsigma(v)\leq \varsigma(u_n)$.
On the other hand, by the mean value theorem there exists some $v\in [u_{n+1},u_n]$ such that
$$
I_{\varsigma}(u_{n+1})-I_{\varsigma}(u_n)=-\frac{u_{n+1}-u_{n}}{\varsigma(v)}\geq \frac{\varsigma(u_n)}{\varsigma(v)}\geq 1.
$$
We conclude that $I_{\varsigma}(u_n)\geq n+I_{\varsigma}(u_0)$. This ends the proof of the lemma.
\cqfd
Observe that for any $t\geq 1$ we have
\begin{eqnarray*}
\varsigma(v)=~v^{1+\chi}&\Longrightarrow&
I_{\varsigma}(u)=\frac{1}{\chi}~(u^{-\chi}-u_0^{-\chi})\\
&\Longrightarrow& I_{\varsigma}^{-1}\left(\omega t\right)=
\left(\frac{1}{u_0^{-\chi}+t\omega~\chi }\right)^{1/\chi}
\leq \left(\chi\omega\right)^{-1/\chi}~t^{-1/\chi}.
\end{eqnarray*}

The above lemma readily yields the following polynomial convergence theorem.
\begin{theo}\label{theo-subgeo}
Assume conditions   (\ref{ref-hyp-varphi}) and  (\ref{loc-min-2}) are satisfied for some function $\varphi_1$ and some parameter $\chi>0$. In this situation, for any $t\geq 1$ and any $\mu\in \Ma_{V}(E)$ with $\mu(1)=0$  we have the polynomial convergence estimates
$$
\vertiii{\mu P_t}_{1+\rho \varphi_1(V)}
\leq~c_{\rho}^{\chi}~t^{-1/\chi}~\vertiii{\mu}_{1+\rho V}
\quad\mbox{with}\quad c_{\rho}^{\chi}:=\left(\chi\omega(\rho)\right)^{-1/\chi}
$$
In the above display, $\rho$ and $\omega(\rho)$ stands for the parameters defined in Lemma~\ref{ref-lem-sub-geo} and in (\ref{def-omega-rho}).
\end{theo}
\begin{rmk}
Without further work, we recover the sequence of polynomial rates of convergence discussed in~\cite{jarner} by choosing 
for any $1<i<n$ the parameters
$$
\begin{array}{l}
\delta:=(n-1)/n\qquad \upsilon:=(i-1)/(n-1)\\
\\
\Longrightarrow 1-\upsilon\delta=1-(i-1)/n\qquad
\delta(1-\upsilon)=1-i/n\qquad 1/\chi=i-1.
\end{array}
$$
\end{rmk}
As expected, the operator-theoretic framework described above can be extended easily to situations where the function $\varphi_2$ in (\ref{ref-hyp-varphi}) has the following form
\begin{equation}\label{weaker-subgeo-cond}
\varphi_2(V):=V~\psi\left(\varphi_1(V)/V\right)\in \Ba_{\infty}(E)
\end{equation} for some convex increasing function $\psi~:~[0,\infty[\mapsto [0,\infty[$ such that $\psi(0)=0$. In this situation, arguing as in (\ref{ref-2-Jensen}) for any $\mu\in \Ma_{V}(E)$ with $\mu(1)=0$ we have
$$
\vertiii{\mu}_{\varphi_2(V)}\geq \vertiii{\mu}_V~
\psi\left({\vertiii{\mu}_{\varphi_1(V)}}/{\vertiii{\mu}_V}\right)
$$
Using  the fact that 
$\psi(\lambda v)\leq \lambda \psi(v)$, for any $\lambda\in [0,1]$ and $v\geq 0$,
for any $t\geq 0$ we check that
\begin{eqnarray*}
\vertiii{\mu P_t}_{\varphi_2(V)}&\geq& \sup_{s\geq 0}\vertiii{\mu P_s}_V~\frac{\vertiii{\mu P_t}_V}{\sup_{s\geq 0}\vertiii{\mu P_s}_V}~
\psi\left({\vertiii{\mu P_t}_{\varphi_1(V)}}/{\vertiii{\mu P_t}_V}\right)\\
&\geq & \sup_{s\geq 0}\vertiii{\mu P_s}_V~\psi\left({\vertiii{\mu P_t}_{\varphi_1(V)}}/{\sup_{s\geq 0}\vertiii{\mu P_s}_V}\right)\
\end{eqnarray*}
On the other hand, by  (\ref{ref-unif-est}) we have $\rho\sup_{s\geq 0}\vertiii{\mu P_s}_V\leq (1+\rho) \vertiii{\mu}_V$.  This yields the rather crude estimate
\begin{eqnarray*}
\vertiii{\mu P_t}_{\varphi_2(V)}
&\geq & \vertiii{\mu}_V~\psi\left(\frac{\rho}{1+\rho}
\frac{\vertiii{\mu P_t}_{\varphi_1(V)}}{\vertiii{\mu}_{V}}\right).
\end{eqnarray*}
Thus, recalling the norm equivalence formulae (\ref{norm-equivalence})  we prove that
\begin{eqnarray*}
\vertiii{\mu P_t}_{\varphi_2(V)}&\geq&  \vertiii{\mu}_{1+\rho V}~\psi_{\rho}\left(
{\vertiii{\mu P_t}_{1+\rho\varphi_1(V)}}/{\vertiii{\mu}_{1+\rho V}}\right)
\end{eqnarray*}
with the rescaled convex function
\begin{equation}\label{rescaled-psi}
\psi_{\rho}(v):=\frac{1}{1+\rho}~\psi\left(\frac{\rho^2}{(1+\rho)^2}~v\right).
\end{equation}
This shows that
$$
v_0:=\vertiii{\mu}_{1+\rho V}\quad\mbox{\rm and}\quad
u_t:=\vertiii{\mu P_{t}}_{1+\rho \varphi_1(V)}\Longrightarrow
u_{t+1}-u_t\leq -\delta_{\rho}(r)
~\varsigma_{v_0}\left(u_t\right)
$$
with the same $\delta_{\rho}(r)$ as in Lemma~\ref{ref-lem-sub-geo} (with $\kappa_2=1$) and the function
$\varsigma_{v_0}(u):=v_0~\psi_{\rho}\left(u/v_0\right)$.
In this context, Lemma~\ref{lem-butkovsky} readily yields the following convergence theorem.
\begin{theo}\label{ref-theo-subgeo-psi}
Assume there exists some function $\varphi_1$ satisfying (\ref{ref-hyp-varphi}), with the  function $\varphi_2$ defined in (\ref{weaker-subgeo-cond}) in terms of some convex increasing function $\psi$. In this situation, for any $t\geq 0$ and any $\mu\in \Ma_{V}(E)$ with $\mu(1)=0$  we have
$$
\vertiii{\mu P_{t}}_{1+\rho\, \varphi_1(V)}\leq J_{\psi_{\rho}}^{-1}(t)~\vertiii{\mu}_{1+\rho V}
$$
with the parameter $\rho$ defined in Lemma~\ref{ref-lem-sub-geo}, the function $\psi_{\rho}$ defined in (\ref{rescaled-psi}) and 
$$
J_{\psi_{\rho}}(u):=\int_{u}^{\iota}\frac{dv}{\psi_{\rho}\left(v\right)}\quad \mbox{with}\quad
\iota:=1\vee \Vert  \varphi_1(V)/V\Vert.$$
\end{theo}
\proof
Observe that
$$
\begin{array}{l}
\displaystyle I_{\varsigma_{v_0}}(u):=\int_u^{u_0}\frac{dv}{\varsigma_{v_0}(v)}=
\int_{u/v_0}^{u_0/v_0}\frac{dv}{\psi_{\rho}\left(v\right)}=:I_{\psi_{\rho}}\left(u/v_0\right).
\end{array}$$
On the other hand, we have
$
u_0\leq \iota~v_0$ and $
I_{\psi_{\rho}}(u)\leq J_{\psi_{\rho}}(u)
$.
Since $I_{\psi_{\rho}}$ and $J_{\psi_{\rho}}$ are decreasing their inverse are also decreasing and choosing $u=I_{\psi_{\rho}}^{-1}(t)$ in the above display we have
$$
\begin{array}{l}
t=I_{\psi_{\rho}}(I_{\psi_{\rho}}^{-1}(t))\leq J_{\psi_{\rho}}(I_{\psi_{\rho}}^{-1}(t))\\
\\
\Longrightarrow J_{\psi_{\rho}}^{-1}(t)\geq 
I_{\psi_{\rho}}^{-1}(t)\Longrightarrow
I_{\varsigma_{v_0}}^{-1}(t)=v_0~I_{\psi_{\rho}}^{-1}(t)\leq v_0~J_{\psi_{\rho}}^{-1}(t).
\end{array}$$
The end of the proof is now a direct consequence of Lemma~\ref{lem-butkovsky}.
\cqfd

\subsection{Some $V$-norm stability theorems}\label{sec-Markov-V-norm}

Let $P_{s,t}$ be  a semigroup of Markov integral operators $P_{s,t}$ on $\Ba_b(E)$
indexed by  continuous time indices $s, t\in \Ta=\RR_+:=[0,\infty[$ or by discrete time indices
 $s,t\in \Ta=\NN$. We further assume that
\begin{equation}\label{V-geo-drift}
 P_{s,s+\tau}(V)\leq \epsilon_{\tau}~V+c_{\tau}
\end{equation}
 for some parameter $\epsilon_{\tau} \in ]0,1[$ and some finite constant $c_{\tau}<\infty$. 
The geometric drift condition (\ref{V-geo-drift}) ensures that  the sequence $\vertiii{P_{s,s+n\tau}}_V$ indexed by $s\geq 0$ and $n\geq 1$ is uniformly bounded.
In this context,  (\ref{ref-V-theta-intro}) applied to $Q_{s,t}=P_{s,t}$ ensures that the operator norms of $P_{s,t}$ are uniformly bounded w.r.t. any time horizon. More precisely, whenever (\ref{V-geo-drift}) is met we have the equivalence 
\begin{equation}\label{discrete-2-continuous-intro-K-M}
 \sup_{s\geq 0 }\sup_{t\geq s}\vertiii{P_{s,t}}_V<\infty  \Longleftrightarrow \sup_{\vert t-s\vert\leq \tau}\vertiii{P_{s,t}}_V<\infty.
\end{equation}
Note that the condition (\ref{discrete-2-continuous-intro-K-M}) is automatically satisfied  whenever (\ref{V-geo-drift}) is met for any $\tau>0$ with $\sup_{\tau\in [0,1]}c_{\tau}<\infty$.
For instance, consider  the Markov transition semigroup $P_{s,t}$ of a continuous time stochastic flow $X_{s,t}(x)$ on some locally compact normed vector space $(E,\Vert\point\Vert)$ with generator $L_t$ defined on some common domain $\Da(L)\subset\Ba(E)$. In this context,  for any non negative function $V\in \Da(L)$ and any parameters $a>0$, $c<\infty$ and $\tau>0$ we have 
\begin{equation}\label{V-L-ex}
\begin{array}{l}
\forall t\in\Ta=\RR_+\qquad L_{t}(V)\leq -a~ V+c
\\
\\
\Longrightarrow (\ref{V-geo-drift})~ \mbox{\rm and}~(\ref{discrete-2-continuous-intro-K-M})\quad \mbox{\rm with}\quad \epsilon_{\tau}=(1+a\tau)^{-1}<1\quad \mbox{\rm and}\quad c_{\tau}=c\tau.
\end{array}
\end{equation}
The above estimate is rather well known,  a detailed proof is provided in the appendix on page~\pageref{V-L-ex-proof}.
Further examples of Markov diffusion semigroups on $\RR^n$ satisfying (\ref{V-geo-drift}) are discussed in Section~\ref{sec-Markov-lyap}.
We further assume there exists some  $r_0\geq 1$  and some function $\alpha_{\tau} : r\in [r_0,\infty[\,\mapsto\alpha_{\tau}(r)\in \, ]0,1]$,  such that for any $r\geq r_0$ we have
\begin{equation}\label{V-dob}
\sup_{(x,y)\in \Va(r)^2}\Vert\delta_xP_{s,s+\tau}-\delta_yP_{s,s+\tau}\Vert_{\tiny tv}\leq 1-\alpha_{\tau}(r),
\end{equation}
with the compact level sets $\Va(r)$ introduced in (\ref{loc-min}).
By Theorem~\ref{theo-intro-Markov-hom}, conditions (\ref{V-geo-drift}), (\ref{discrete-2-continuous-intro-K-M}) and (\ref{V-dob})  ensure the existence of some parameter $\overline{\tau}>0$  such that \begin{equation}
\sup_{\vert t-s\vert\leq \overline{\tau}}{ \beta_{V}(P_{s,t})}<\infty\quad\mbox{\rm and}\quad
\sup_{s\geq 0} \beta_{V}(P_{s,s+\overline{\tau}})<1. \label{beta-V-norm}
 \end{equation}
 In the above display, $\beta_{V}(P_{s,t})$ stands for the $V$-Dobrushin coefficient  of the Markov operator  $P_{s,t}$ introduced in (\ref{defi-beta-V}).
The next  exponential contraction theorem is a direct consequence of the operator norm estimates (\ref{beta-V-norm}) and it is valid on abstract measurable spaces as well as for any function $V\geq 1$.
\begin{theo}\label{theo-intro-Markov}
Let $P_{s,t}$ be a semigroup of Markov integral operators  $P_{s,t}$ on some measurable state space $E$ satisfying  condition (\ref{beta-V-norm})  for some function $V\geq 1$ and some parameter $\overline{\tau}>0$.
In this situation, there exists a parameter $b>0$ and some finite constant $c<\infty$ such that for any $s\leq t$ and $\mu,\eta\in \Pa_V(E)$ we have the exponential estimate
\begin{equation}\label{lipschitz-inq-M}
\vertiii{(\mu-\eta) P_{s,t}}_{\tiny V}\leq~c~e^{-b(t-s)}~\vertiii{ \mu-\eta}_{\tiny V}.
\end{equation}
In particular, the above exponential Lipschitz estimates are met as soon as conditions (\ref{V-geo-drift}), (\ref{discrete-2-continuous-intro-K-M}) and (\ref{V-dob}) are satisfied. The estimates (\ref{lipschitz-inq-M}) also hold for any $s\geq 0$ and $t\in [s,\infty[_{\tau}$ 
as soon as  (\ref{V-geo-drift}) and (\ref{V-dob}) are satisfied for some $\tau>0$ and  $\epsilon_{\tau} \in ]0,1[$. 
\end{theo}

Theorem~\ref{theo-intro-Markov} can be seen as an extension
 of Harris'  theorem  to time varying Markov semigroup. 
The proof of Theorem~\ref{theo-intro-Markov} is based on the {\em discrete time} $V$-norm operator contraction techniques presented in Section~\ref{ref-sec-discrete}. The r.h.s. condition in (\ref{discrete-2-continuous-intro-K-M}) is a technical condition only made for continuous time semigroups to ensure that (\ref{lipschitz-inq-M}) also holds for continuous time indices.
Note that the strength of conditions (\ref{V-geo-drift}) and (\ref{V-dob}) depends on the strength of the function $V$: when the function $V$ is bounded, the  geometric drift condition (\ref{V-geo-drift}) and 
the uniform norm condition (\ref{discrete-2-continuous-intro-K-M}) are trivially met but in this case condition (\ref{V-dob}) is a uniform contraction condition on the state $E$.
In the reverse angle, when $V\in \Ba_{\infty}(E)$ is a function with compact sub-level sets, the geometric drift condition  (\ref{V-geo-drift}) combined with (\ref{discrete-2-continuous-intro-K-M}) ensures that $\mu P_{s,t}$ is a tight collection  of probability measures indexed by $s\leq t$. For time homogenous models $P_{s,s+t}=P_t$, following Remark~\ref{rmk-unif-tv}, uniform total variation estimates can be derived for any $t\geq \tau$ as soon $\Vert P_{\tau}(V)\Vert<\infty$. Some examples satisfying this condition are discussed in Section~\ref{ricc-sec}.

Using (\ref{ref-coupling-tv}) we readily check that the local contraction condition (\ref{V-dob}) is met if and only if
for any $s\geq 0$ and any $(x,y)\in \Va(r)^2$ there exists some probability measure $\mu$ on $E$ (that may depends on the parameters $(\tau,r,s,x,y)$) such that
$$
\forall z\in\{x,y\}\qquad
\delta_zP_{s,s+\tau}(dy)\geq \alpha_{\tau}(r)~\mu(dy).
$$
For instance, the above condition is met as soon as
\begin{equation}\label{ref-P-min}
P_{s,s+\tau}(x,dy)\geq p_{s,s+\tau}(x,y)~\nu_{\tau}(dy)
\end{equation}
for some  Radon positive measure $\nu_{\tau}$ on $E$ and some density function $p_{s,s+\tau}$, satisfying for any $r\geq r_0$ the local minorization condition  
 \begin{equation}\label{min-p}
  0<\inf_{s\in \Ta}\inf_{\Va(r)^2}p_{s,s+\tau}\quad \mbox{\rm and}\quad 0<\nu_{\tau}(\Va(r))<\infty.
  \end{equation}
For locally compact Polish spaces condition $0<\nu_{\tau}(\Va(r))<\infty$ is met as soon as
$V$ has compact sub-levels sets $\Va(r)$ with non empty interior and $\nu_{\tau}$ is a
Radon measure of full support; that is $\nu_{\tau}$ is finite on compact sets and strictly
positive on non-empty open sets. For time homogeneous models, also note that the l.h.s. minorization condition (\ref{min-p}) is satisfied as soon as $(x,y)\in (E^{\circ})^2\mapsto p_{\tau}(x,y)$ is a continuous positive function on the interior $E^{\circ}$ of the set  $E$.

Several illustrations of Theorem~\ref{theo-intro-Markov} are discussed in Section~\ref{sec-Markov-lyap} in the context of diffusion processes on Euclidean spaces as well as in Section~\ref{ricc-sec} in the context of Riccati-type diffusion on positive definite matrix spaces and multivariate birth and death jump type processes on countable state spaces. The stability of Markov semigroups on manifolds with entrance boundaries can also be analyzed using the Lyapunov techniques developed in Section~\ref{boundary-sec}. For instance, as shown in Section~\ref{bounded-dom-sec}, any absolutely continuous Markov semigroup $P_{s,t}$ on a bounded connected subset $E\subset \RR^n$ with locally Lipschitz boundary 
$\partial E$ satisfies the conditions of Theorem~\ref{theo-intro-Markov} with the (non unique) Lyapunov function $V(x)=1/\sqrt{d(x,\partial E)}$ and
 the  distance to the boundary defined for any $x\in E$ by
$$
d(x,\partial E):=\inf{\left\{\Vert x-y\Vert~:~y\in \partial E\right\}}.
$$
We illustrate the above discussion 
with some elementary one dimensional examples.
\begin{examp}
Consider a one dimensional 
Brownian motion  $X_t(x)$ starting at $X_0(x)=x\in E:=[0,1]$ and reflected at the boundaries $\partial E=\{0,1\}$. 
We recall that the Markov transition of the process $t\in\Ta:=\RR_+\mapsto X_t(x)\in E$ is symmetric and absolutely continuous; that is we have
$$
P_t(x,dy):=\PP(X_t(x)\in dy)=p_t(x,y)~\nu(dy)\quad 
\mbox{\rm with}\quad \nu(dy):=1_{[0,1]}(y)dy
$$
and the density $p_t(x,y)$ is given by the spectral decomposition
$$
p_t(x,y)=1+2 \sum_{n\geq 1}~e^{-(n\pi)^2t/2}~\cos{(n\pi x)}~\cos{(n\pi y)}
$$
In this situation,  $P_t=P_{0,t}$ coincides with the Neumann heat semigroup on $[0,1]$.
Since the Neumann heat kernel $p_t(x,y)$ is smooth as well as bounded and strictly positive on the compact interval $[0,1]$, the conditions of Theorem~\ref{theo-intro-Markov} are satisfied with the unit Lyapunov function $V(x)=1$, as well as for any of the Lyapunov functions $V(x)=1/\sqrt{x}$, $V(x)=1/\sqrt{1-x}$ or $V(x)=1/\sqrt{x}+1/\sqrt{1-x}$.  Indeed, note that the minorization condition (\ref{min-p}) holds for any of the Lyapunov 
functions $V$ discussed above. Since $\nu(V)\leq \infty$,  we  have $\Vert P_{\tau}(V)\Vert<\infty$ for any $\tau>0$, so that the Lyapunov condition (\ref{V-geo-drift}) is also met.

The same reasoning applies to the one dimensional positive Riccati-type diffusions  with an entrance boundary at the origin discussed in Section~\ref{ricc-sec}. Reflecting this class of positive diffusions at $x=1$, 
the conditions of Theorem~\ref{theo-intro-Markov} are satisfied on $E=]0,1]$ with  the Lyapunov functions $V(x)=1/\sqrt{x}$ as well as for $V(x)=1/\sqrt{x}+1/\sqrt{1-x}$.
\end{examp}

In the same vein, assume there exists some  increasing  differentiable concave function $ \varphi~:~v\in [1,\infty[\mapsto \varphi(v)\in [0,\infty[$  with bounded differential $\Vert\partial  \varphi\Vert<\infty$  and some function $V\in \Ba_{\infty}(E)$ such that $V\geq 1$ and $V,\varphi(V)\in\Da(L)$. In addition, there exists some finite constant $c>\varphi(1)$ such that for any $t\in\Ta=\RR_+$ we have
\begin{equation}\label{L-sub-geo}
L_t(V)\leq - \varphi(V)+c\quad \mbox{\rm and}\quad \varphi(V)/V\in \Ba_0(E)
\end{equation}
We set
\begin{equation}\label{ref-c-tau-sub}
V_{\tau}:=V+\tau \, \varphi(V)\quad \mbox{and}\quad 
\varphi_{\tau}(v)=
\varphi\left(v/\beta_{\tau}\right)\quad \mbox{with}\quad 
 \beta_{\tau}:=1+\tau \Vert  \varphi(V)/V\Vert.
\end{equation}
Observe that
$$
\varphi_{\tau}\left(V_{\tau}\right)=\varphi\left(V\left(1+\tau  \varphi(V)/V\right)/\beta_{\tau}\right)\leq 
\varphi\left( V\right)
\quad\mbox{\rm with}\quad \beta_{\tau}:=1+\tau \Vert  \varphi(V)/V\Vert
$$
and
$$
\varphi_{\tau}\left(V_{\tau}\right)/V_{\tau}\leq \varphi\left( V\right)/V\in \Ba_0(E)
$$
This estimate ensures that $\varphi_{\tau}\left(V_{\tau}\right)/V_{\tau}\in \Ba_0(E)$ as soon as $\varphi_{\tau}\left(V_{\tau}\right)/V_{\tau}$ is locally lower bounded and upper semi-continuous.
\begin{lem}\label{lem-211}
Assume condition (\ref{L-sub-geo}) is satisfied for some functions $(\varphi,V)$ and some constant $c>\varphi(1)$. In this situation,
for any $s\in\Ta$ and $\tau>0$ we have
$$
 \forall s\in\Ta\qquad
P_{s,s+\tau}(V_{\tau})\leq V_{\tau}-\tau ~\varphi_{\tau}(V_{\tau})+c_{\tau}
$$
with the function $V_{\tau}$ defined in (\ref{ref-c-tau-sub}) and $
c_{\tau}:=c~\tau~\left(1+\Vert \partial  \varphi\Vert
~{\tau}/{2}\right).
$
\end{lem}
The proof of the above lemma is provided in the Appendix, on page~\pageref{lem-211-proof}.
Next theorem is the continuous time version of the polynomial convergence theorem, Theorem~\ref{theo-subgeo} presented in Section~\ref{subgeo-sec}. The continuous time version of Theorem~\ref{ref-theo-subgeo-psi} can be obtained using the same lines of arguments, thus it is left to the reader.
\begin{theo}
Assume  (\ref{L-sub-geo}) is met for some function $(\varphi,V)$ and the function $\varphi_{\tau}$ defined in (\ref{ref-c-tau-sub}) satisfies  (\ref{ref-hyp-varphi}) for some parameters $\tau,\chi>0$ and some  functions $(\varphi_{\tau,1},\varphi_{\tau,2})$ such that $\varphi_{\tau}(V_{\tau})/V_{\tau},\varphi_{\tau,1}(V_{\tau})/V_{\tau}\in \Ba_0(E)$ and $\varphi_{\tau,2}(V_{\tau})\in \Ba_{\infty}(E)$. We also assume that (\ref{V-dob}) is satisfied with the compact level sets of the function $V_{\tau}$.
In this situation, there exists some constant $c<\infty$ (that may depends on $\tau$)  such that for any $s\geq 0$ and any $t\in [s,\infty[_{\tau}$ and  $\mu,\eta\in \Pa_V(E)$we have
$$
\vertiii{(\mu-\eta) P_{s,t}}_{\tiny tv}\leq c~(t-s)^{-1/\chi}
~\vertiii{\mu-\eta}_{V}
$$
The above polynomial convergence estimates also holds for any continuous time indices as soon as $\sup_{\vert s-t\vert\leq \tau}\vertiii{P_{s,t}}_{V}<\infty$.
\end{theo}
\proof
Applying Theorem~\ref{theo-subgeo}, there exists some parameter $\rho_{\tau}>0$ and some finite constant $c_{\tau}$ such that for any $n\in \NN$ we have the polynomial convergence estimates
$$
\vertiii{\mu P_{s,s+n\tau}}_{1+\rho_{\tau} \varphi_{\tau, 1}(V_{\tau})}
\leq~c_{\tau}~n^{-1/\chi}~\vertiii{\mu}_{1+\rho_{\tau} V_{\tau}}\leq c_{\tau}~(1+\rho_{\tau})~n^{-1/\chi}~\vertiii{\mu}_{V_{\tau}}
$$
This implies that
$$
\vertiii{\mu P_{s,s+n\tau}}_{\tiny tv}\leq c~n^{-1/\chi}~\vertiii{\mu}_{V}\quad \mbox{\rm with}\quad
c=c_{\tau}~\beta_{\tau}(1+\rho_{\tau})
$$
This ends the proof of the theorem.\cqfd
\subsection{Diffusion semigroups}\label{sec-Markov-lyap}
This section is mainly concerned with the design of Lyapunov functions for continuous time Markov semigroups. To simplify notation, we only consider time homogeneous models. 
 All the semigroups discussed in this section satisfy condition (\ref{discrete-2-continuous-intro-K-M}). Thus, by (\ref{ref-P-min}) the contraction theorem, Theorem~\ref{theo-intro-Markov}, applies to all the Markov semigroups discussed in this section as soon as the transition semigroups have a continuous density with respect to the Lebesgue measure. 

Section~\ref{sec-gen-Markov} presents some elementary principles based on spectral conditions on the drift function and a simple way to design Lyapunov functions in terms of the generator of diffusion process. These generator-type techniques are illustrated in Section~\ref{over-damped-sec} in the context of overdamped Langevin diffusions. The design of Lyapunov functions for hypo-elliptic diffusions and Langevin diffusions are discussed respectively in Section~\ref{hypo-elliptic-sec} and Section~\ref{Langevin-example}.
 
\subsubsection{Some general principles}\label{sec-gen-Markov}
Consider the Markov semigroup $P_t$ of a diffusion flow $X_{t}(x) $ on $E=\RR^n$ defined  by
\begin{equation}\label{ref-X-diff-stable}
dX_{t}(x)=b(X_{t}(x))\,dt+\sigma(X_{t}(x))\,dB_t.
\end{equation}
In the above display, $B_t$ is a $n_1$-dimensional Brownian motion starting at the origin for some $n\geq 1$, $b$ is a differentiable drift function from $\RR^n$ into itself with gradient-matrix $\nabla b=(\partial_{x_i}b^j)_{1\leq i,j\leq n}$, and $\sigma$ stands for some diffusion function from $\RR^n$ into $\RR^{n\times n_1}$. We set $\Sigma^2:=\sigma\sigma^{\prime}$, where $\sigma^{\prime}(x):=\sigma(x)^{\prime}$ stands for the transposition of the matrix $\sigma(x)$, so that $\Sigma^2(x):=\sigma(x)\sigma^{\prime}(x)$. The absolutely continuity of the transition semigroup
$P_t(x,dy)=\PP(X_t(x)\in dy)=p_t(x,y)\nu(dy)$ for some continuous transition densities $p_t(x,y)$ (w.r.t. the Lebesgue measure $\nu(dy)$) is ensured as soon as
$(b,\sigma)$ are globally Lipschitz continuous and the diffusion matrix is invertible or more generally satisfying a parabolic H\" ormander condition (see for instance~\cite{nualart,olivera-tudor,sanz} and references therein).
The generator $L$  of the diffusion flow $X_t(x)$ and its carr\'e du champ operator $\Gamma_L$ are given respectively by the formula
\begin{equation}
L(f):=b^{\prime}\nabla f+\frac{1}{2}~\tr\left(\Sigma^2\nabla^2f\right)\quad \mbox{\rm and}\quad
\Gamma_{L}(f,g):= (\nabla f)^{\prime}\Sigma^2 \nabla g.
\label{ref-gen-La} 
\end{equation}

The next proposition provides a rather elementary way to design a Lyapunov function.
\begin{prop}\label{prop-ref-tangent}
Assume that $\sigma(x)=\sigma_0$ for some $\sigma_0\in\RR^{n\times n_1}$ and we have
\begin{equation}\label{ref-tangent-X-stable}
\nabla b+(\nabla b)^{\prime}\leq -2\lambda~I\quad\mbox{for some}\quad \lambda>0.
\end{equation}
Then for any $v>0$ and $t>0$ there exists some  $\delta_t>0$ such that
\begin{equation}
V(x):=\exp{\left(v\Vert x\Vert\right)}\Longrightarrow P_{t}(V)/V\leq c_t/V^{\delta_t}.
\end{equation}
\end{prop}
The proof of the above proposition is rather technical, thus it is provided in the appendix on page~\pageref{prop-ref-tangent-proof}.

The next proposition is a slight extension of Theorem 2.6~\cite{kusuoka-stroock} on reversible semigroups to stochastic flows in Euclidean spaces. It provides a rather simple way to design Lyapunov functions in terms of generators.
\begin{prop}\label{condition-ks-prop-0}
Assume there exists some function $W\geq 0$ as well as some parameters
$\alpha>0$, $\beta\in\RR$ and $0<\epsilon< 1$ such that
\begin{equation}\label{condition-ks-prop}
\alpha~W +\beta+L(W)\leq -\epsilon~\Gamma_{L}(W,W).
\end{equation}
In this situation, for any $t> 0$ we have
\begin{equation}\label{expo-estimate-LW-Gamma}
V:=\exp{\left(2\epsilon W\right)}\Longrightarrow
P_t\left(V\right)/V
\leq v_t/V^{\delta_t}
\end{equation}
with the parameters
$$
v_t=\exp{\left(-2\beta\epsilon~(1-e^{-\alpha t})/{\alpha}\right)}\quad \mbox{and}\quad \delta_t:=(1-e^{-\alpha t}).
$$
\end{prop}
The proof of the above proposition follows word-for-word the proof of Theorem 2.6 in~\cite{kusuoka-stroock}, thus it is provided in the appendix on page~\pageref{condition-ks-prop-proof}.

We further assume that $P_t$ satisfies for any $t>0$ the sub-Gaussian estimate 
\begin{equation}\label{def-OU-A-Sigma-lem}
P_{t}(x,dy)\leq c_{t}~\exp{\left(-\frac{1}{2\sigma^2_{t}}\,\Vert y-m_{t}(x)\Vert^2\right)}~dy
\end{equation}
for some parameters $\sigma_{t}>0$ and some 
 some function $m_{t}$ on $\RR^n$ such that  
$$
\Vert m_t(x)\Vert\leq c_t~(1+\Vert x\Vert).
$$ 
In this situation, for any $n\geq 1$ and $t\geq 0$ we have
$$
V(x):=1+\Vert x\Vert^n \Longrightarrow \Vert P_t(V)/V\Vert<\infty.
$$ 
More refined estimates can be found when the function $m_t$ is such that
\begin{equation}\label{def-OU-A-Sigma-lem-epsilon}
 \vert m_{t}(x)\vert\leq \epsilon_{t}~\vert x\vert\quad \mbox{with}\quad \epsilon_{t}\in ]0,1[
\end{equation}
for some norm $\vert\point\vert$ on $\RR^n$.
In this situation, observe that any $v\geq 0$ and any centered Gaussian random variable $Y$ on $\RR^n$ with identity covariance matrix $I_n$ we have
$$
e^{-v\vert x\vert }~\EE\left(e^{v\vert m_t(x)+\sigma_t Y\vert}\right)\leq c_t~
e^{-v (1-\epsilon_t) \vert x\vert}.
$$
This yields the following lemma.
\begin{lem}\label{ref-lem-OU-A-Sigma}
Consider a Markov semigroup $P_{t}$ satisfying the sub-Gaussian estimate (\ref{def-OU-A-Sigma-lem}) as well as (\ref{def-OU-A-Sigma-lem-epsilon}) for some norm $\vert\point\vert$ on $\RR^n$.
Then for any $v\geq 0$ and $t>0$ there also exists some finite constant $\delta_t>0$ such that
$$
V(x):=\exp{\left(v\vert x\vert\right)}\Longrightarrow
P_{t}(V)/V\leq c_{t}/V^{\delta_{t}}. 
$$

\end{lem}

\subsubsection{Overdamped Langevin diffusion}\label{over-damped-sec}
Let $W(x)$ be some twice differentiable potential function from $\RR^n$ into $\RR$. The overdamped Langevin diffusion is defined by choosing in (\ref{ref-X-diff-stable}) the drift function
$$
b(x):=-\gamma~\nabla W(x)\quad\mbox{\rm and}\quad (n_1,\sigma(x))=(n,\rho~I)
\quad\mbox{\rm for some}\quad\gamma,\rho>0.
$$
 In this context, we have
 $$
 (\ref{ref-tangent-X-stable})\Longleftrightarrow \nabla^2 W \geq (\lambda/\gamma)~I\quad\mbox{for some}\quad \lambda>0.
 $$
 Also observe that
 $$
(\ref{condition-ks-prop})\Longleftrightarrow
\alpha~W +\beta+\frac{\rho^2}{2}~ \tr(\nabla^2W) \leq \left(\gamma-\epsilon~\rho^2\right)~\Vert\nabla W\Vert^2.
 $$
The above condition is clearly met when $W$ behaves as $\Vert x\Vert^{m}$ with $m\geq 1$ at infinity; that is, there exists some
 sufficiently large radius $r$ such that for any $\Vert x\Vert\geq r$ we have
$$
 \left\vert \tr(\nabla^2W(x))\right\vert\leq c_1~\|x\|^{(m-2)_+}
\quad\mbox{\rm
and}\quad
\left\Vert\nabla W(x) \right\Vert^2\geq c_2~\|x\|^{2(m-1)}. 
$$

\subsubsection{Hypo-elliptic diffusions}\label{hypo-elliptic-sec}
Consider the $\RR^n$-valued diffusion (\ref{ref-X-diff-stable}) with $(b(x),\sigma(x))=(A x,\Sigma)$, 
for some matrices $(A,\Sigma)$ with appropriate dimensions. We assume that $A$ is stable (a.k.a. Hurwitz); that is its spectral abscissa $ \varsigma(A)$ defined below is negative 
\begin{equation}\label{def-spectral-abscissa}
 \varsigma(A):=\sup\left\{
\mbox{\rm Re}\left(\lambda(A)\right)~:~\lambda(A)\in\mbox{\rm Spec}(A)\right\}<0.
\end{equation}
In the above display $\mbox{\rm Spec}(A)$ denotes  the spectrum of the matrix $A$, and $\mbox{\rm Re}\left(\lambda(A)\right)$ the real part of  $\lambda(A)$.
We also assume  that $R:=\Sigma\Sigma^{\prime}$ is positive semi-definite 
and the pair of matrices $(A,R^{1/2})$ are controllable, 
 in the sense that the $(n\times n^2)$-matrix
\begin{equation}\label{def-contr-obs}
\left[R^{1/2},AR^{1/2}\ldots, A^{r-1}R^{1/2}\right] \quad\mbox{has rank $n$}.
\end{equation}
Whenever  $\varsigma(A)<0$ we have
\begin{equation}\label{def-Ma-OU}
 P_{t}(x,dy)=\frac{1}{\sqrt{\mbox{\rm det}(2\pi C_t)}}
 \exp{\left(-\frac{1}{2}\left( y-m_{t}(x)\right)^{\prime}C^{-1}_{t}\left( y-m_{t}(x)\right)\right)}~dy
\end{equation}
 with the mean value function
 $$
x\mapsto m_t(x):=e^{tA}x\longrightarrow_{t\rightarrow\infty} 0
 $$
 and the covariance matrices $C_t$ defined for any $t>0$ by
 $$
0< C_t:=\int_0^t~e^{s A}Re^{s A^{\prime}}~ds\longrightarrow_{t\rightarrow\infty}
C_{\infty}:=\int_0^{\infty}~e^{s A}Re^{s A^{\prime}}~ds.
 $$
 Since $A$ is stable, there exists some norm $\vert \point\vert$ on $\RR^n$ such that
the corresponding  operator norm satisfies $\vert e^{tA}\vert\leq e^{l(A) t}$ for some log-norm parameter $l(A)<0$. This implies that
\begin{equation}\label{stable-exp}
\vert m_t(x)\vert=\vert e^{tA}x\vert\leq e^{l(A) t}~\vert x\vert.
\end{equation}
This clearly shows that the semigroup $P_t$ of the hypo-elliptic Ornstein-Ulhenbeck diffusion satisfies (\ref{def-OU-A-Sigma-lem}) and (\ref{def-OU-A-Sigma-lem-epsilon}), and thus the conditions of Lemma~\ref{ref-lem-OU-A-Sigma} are met.

Let $\Pa_t$ be the Markov semigroup of the $\RR^n$-valued linear diffusion
\begin{equation}\label{def-OU-b}
d\Xa_t(x)=\left(A\Xa_t(x)+a(\Xa_t(x))\right)dt+\Sigma~dB_t
\end{equation}
with some bounded drift function $a$ on $\RR^n$, an $(n\times n)$-matrix $A$ satisfying (\ref{def-spectral-abscissa}), some $n_1$-valued Brownian motion $B_t$ starting at the origin and some $(n\times n_1)$-matrix $\Sigma$ satisfying the
rank condition (\ref{def-contr-obs}). 

Using the stochastic interpolation formula (cf. Theorem 1.2 in~\cite{dss-19}) given by
$$
\Xa_t(x)-X_t(x)=\int_0^t~e^{(t-s)A^{\prime}}~a\left(\Xa_s(x)\right)~ds
$$
we check the almost sure estimate
$$
\vert  X_t(x)-\Xa_t(x)\vert\leq c
\quad\mbox{\rm for some finite constant $c<\infty$}.
$$
This yields the following proposition.
\begin{prop}\label{prop-ref-lang}
For any $v>0$ and $t>0$ there exists some  $\delta_{t}>0$  such that
$$
 V(x):=\exp{\left(v \vert x\vert\right)}\Longrightarrow \Pa_{t}(V)/V\leq c_{t}/
V^{\delta_{t}}.
$$
\end{prop}

\subsubsection{Langevin diffusion}\label{Langevin-example}
Consider the Langevin diffusion diffusion flow $\Xa_t(z)=(X_t(z),Y_t(z))\in (\RR^{r}\times\RR^r)$ starting at
$z=(x,y)\in (\RR^{r}\times\RR^r)$ and given by 
\begin{eqnarray*}
dX_t(z)&=&Y_t(z)/m~dt\\
dY_t(z)&=&\left(b(X_t(z))-\beta Y_t(z)/m\right)~dt+\sigma~dB_t.
\end{eqnarray*}
In the above display, $B_t$ stands for an $r$-dimensional Brownian motion starting at the origin, $\sigma,\beta,m>0$  some parameters and  $b$ a function of the form
$$
b(x):=-\gamma~x+a(x)\quad \mbox{with}\quad \gamma>0 \quad \mbox{and}\quad \Vert a\Vert<\infty.
$$
In statistical physics, the above diffusion represents the evolution of $N$ particles $X_t(z)=(X^i_t(z))_{1\leq i\leq N}\in \RR^{3N}$ with mass $m>0$, position $X_t(z)\in \RR^{3N}$ and momenta $Y_t(z)$. In this context, $\gamma>0$ stands for some friction parameter, and the diffusion parameter $\sigma>0$ is related to the Boltzmann constant and the temperature of the system. In this context, the function  $b(x)=-\nabla W(x)$ is often described by the gradient of some potential function $W$. For instance, for a quadratic confinement we have
$$
\begin{array}{l}
W(x):=\gamma \Vert x\Vert^2/2+w(x)\quad \mbox{with}\quad \Vert \nabla w\Vert<\infty\\
\\
\Longrightarrow
b(x)=-\nabla W(x):=-\gamma~x+a(x)\quad \mbox{and}\quad a(x)=\nabla w(x).
\end{array}$$
Notice that $\Xa_t(z)$ can be rewritten in vector form 
as in (\ref{def-OU-b}) with $n=2r$,  $a(x,y)=\left(
\begin{array}{c}
0\\
a(x)
\end{array}
\right)$ and the matrices
\begin{equation}\label{ref-A-Sigma}
A=\left(
\begin{array}{cc}
0&m^{-1}~I_{n\times n}\\
-\gamma~I_{n\times n} &-\beta m^{-1}~I_{n\times n}
\end{array}
\right)\quad \mbox{and}\quad \Sigma:=\left(
\begin{array}{cc}
0&0\\
0&\sigma I_{n\times n}
\end{array}
\right).
 \end{equation}
 It is a simple exercise to check that $A$ satisfies 
(\ref{def-spectral-abscissa}) and  (\ref{def-contr-obs}).

Consider the $\RR^2$-valued stochastic process $X_t=(q_t,p_t)$ defined by
\begin{equation}\label{Langevin-2d}
\left\{
\begin{array}{rcl}
\displaystyle dq_t&=&\displaystyle\beta~\frac{p_t}{m}~dt\\
\displaystyle dp_t&=&\displaystyle-\beta~\left(\frac{\partial W}{\partial q}(q_t)+\frac{\sigma^2}{2}~\frac{p_t}{m}\right)~dt+\sigma~dB_t
\end{array}
\right.
\end{equation}
with some positive constants $\beta,m,\sigma$, a Brownian motion $B_t$, and a smooth positive function $W$ on $\RR$ such that
for sufficiently large $r$ we have
$$
\forall ~\vert q\vert\geq r\qquad q\frac{\partial W}{\partial q}(q)\geq \delta~\left(W(q)+q^2\right)
$$
for some positive constant $\delta$. This condition is clearly met when $W$ behaves as $q^{2 l}$   for certain $l\geq 1$ at infinity. We let $V(q,p)$ be the function on $\RR^2$ defined by
$$
V(q,p)=1+\frac{1}{2m}~p^2+W(q)+\frac{\epsilon}{2}~\left(\frac{\sigma^2}{2}~q^2+2 pq\right)
\quad
\mbox{\rm with}\quad \epsilon<\frac{\sigma^2}{2m} .$$
In this situation, there exists some $a>0$ and $c<\infty$ such that
\begin{equation}\label{langevin-Lac}
L(V)\leq  -aV~+c.
\end{equation}
The proof of the above estimate is rather technical, thus it is provided in the appendix on page~\pageref{langevin-Lac-proof}.

\section{Stability of $V$-positive semigroups}\label{brief-review-sec-pos}

\subsection{Normalized semigroups}\label{norm-sg-V-norm}

For non necessarily Markov $V$-positive semigroups $Q_{s,t}$ one natural idea is to normalize the semigroups.
For any probability measure  $\eta\in \Pa_V(E)$  we let   $\Phi_{s,t}(\eta)\in  \Pa_V(E)$ be the normalized distribution defined for any $f\in\Ba_{V}(E)$ by the formula
\begin{equation}\label{def-Phi-s-t-intro-0}
\Phi_{s,t}(\eta)(f):=\frac{\eta Q_{s,t}(f)}{\eta Q_{s,t}(1)}\quad \mbox{\rm and we set}\quad
\overline{Q}_{s,t}(f)(x):=\frac{Q_{s,t}(f)(x)}{Q_{s,t}(1)(x)}=\Phi_{s,t}(\delta_x)(f).
\end{equation}
The mapping $\Phi_{s,t}$ is a well defined semigroup on $\Pa_V(E)$. The 
denormalisation formula connecting these semigroups is given for any $t\in [s,+\infty[_{\tau}$ by
\begin{equation}\label{product}
\mu Q_{s,t}(f)=\Phi_{s,t}(\mu)(f)~\prod_{u\in [s,t[_{\tau}}~\Phi_{s,u}(\mu)(Q_{u,u+\tau}(1)).
\end{equation}
with
$$
 [s,t[_{\tau}:=\{s+n\tau \in [s,t[~:~n\in\NN\}.
$$
To check this claim, observe that for any $t:=s+n\tau$ we have
$$
\begin{array}{l}
\Phi_{s,s+p\tau}(\mu)(Q_{s+p\tau,s+(p+1)\tau}(1))={\mu Q_{s,s+(p+1)\tau}(1)}/{\mu Q_{s,s+p\tau}(1)}\\
\\
\mbox{\rm and therefore}\quad \prod_{0\leq p<n}~\Phi_{s,s+p\tau}(\mu)(Q_{s+p\tau,s+(p+1)\tau}(1))=\mu Q_{s,s+n\tau}(1).
\end{array}
$$
The above formula coincides with the product  formula relating the unnormalized operators
$Q_{s,t}$
with the normalized semigroup $\Phi_{s,t}$ discussed in \cite[Section 1.3.2]{dm-2000}, see also \cite[Proposition 2.3.1]{dm-04} and \cite[Section 12.2.1]{dm-13}.

Observe that for any  $u\in [s,t[_{\tau}$ and $f\in\Ba_{b}(E)$ we have
$$
\overline{Q}_{s,t}=R^{(t)}_{s,u}\overline{Q}_{u,t}\quad \mbox{\rm with}\quad
R^{(t)}_{s,u}(f):={Q_{s,u}(f~Q_{u,t}(1)
 )}/{Q_{s,u}(Q_{u,t}(1))}
$$
Note that $\{R^{(t)}_{u,u+\tau},~u\in [s,t[_{\tau}\}$ forms a triangular array of time varying Markov semigroups. In terms of the standard Dobrushin coefficient discussed in Section~\ref{ref-sec-discrete} assume there exists some $\delta_{\tau} \in ]0,1[$ such that for any $t\in [s,+\infty[_{\tau}$ and $u\in [s,t[_{\tau}$ we have the strong contraction condition  $\beta(R^{(t)}_{u,u+\tau})\leq (1-\delta_{\tau})$. In this situation, for any $t\in [s,+\infty[_{\tau}$ we have the operator norm contraction
$$
\beta(\overline{Q}_{s,t})\leq\prod_{u\in [s,t[_{\tau}}
\beta(R^{(t)}_{u,u+\tau})\leq (1-\delta_{\tau} )^{(t-s)/\tau}
$$
For a more thorough discussion on these strong contraction conditions we refer the reader to Section 4 in~\cite{dhj-21} and references therein. Unfortunately, the strong contraction condition discussed above is rarely satisfied for non-compact state spaces and we need to resort to Lyapunov conditions. In this connection, note that the Lyapunov condition in (\ref{ref-V-theta-intro})
is stronger than the one discussed in (\ref{V-geo-drift}) in the context of Markov semigroups.
We also strengthen (\ref{ref-P-min}) and assume that  for any $s\geq 0$ and $\tau>0$, the integral operator $Q_{t,t+\tau}$ has
a density $ q_{s,s+\tau}$ with respect to some Radon positive measure $\nu_\tau$ on $E$; that is we have the formula
\begin{equation}\label{ref-Q-chi}
Q_{s,s+\tau}(x,dy)=q_{s,s+\tau}(x,y)~\nu_{\tau}(dy).
\end{equation}
We also assume there exists some $r_0>1$ such that for any $r\geq r_0$ we have
 \begin{equation}\label{min-max-q-vr}
  0<\iota_r(\tau):=\inf_{s\in \Ta}\inf_{\Va(r)^2}q_{s,s+\tau}\leq 
\sup_{s\in \Ta}\sup_{\Va(r)^2}q_{s,s+\tau}<\infty \quad \mbox{\rm and}\quad \nu_{\tau}(\Va(r))>0.
\end{equation}
In this situation,  for any $r\geq r_0$ and $\overline{r}\geq r$ we have the uniform estimate
$$
\inf_{\Va(r)}Q_{s,s+\tau}(1)\geq \jmath_{r,\overline{r}}(\tau):=
\inf_{\Va(r)}Q_{s,s+\tau}(1_{\Va(\overline{r})})\geq \iota_r(\tau)~\nu_{\tau}(\Va(r))>0.
$$
We associate with a given $\mu\in \Pa_V(E)$ and  some function $H\in \Ba_{0,V}(E)$  the finite rank (and hence compact) operator
$$
f\in \Ba_V(E)\mapsto T^{\mu,H}_{s,t}(f):=\frac{Q_{s,t}(H)}{\mu_s(Q_{s,t}(1))}~\mu_t(f)\in \Ca_{V}(E)
$$ 
with the flow of measures $\mu_t=\Phi_{s,t}(\mu_s)$ starting at  $\mu_0=\mu$. With this notation at hand, one has the following theorem.
\begin{theo}[\cite{dhj-21}]\label{theo-intro}
Consider a  $V$-positive semigroups $Q_{s,t}$
with a density (\ref{ref-Q-chi}) satisfying (\ref{min-max-q-vr}) for some parameter $\tau>0$ and some $r_0>1$. In this situation, there exists a parameter $b>0$ such that for any $\mu,\eta\in \Pa_V(E)$  and any $s\geq 0$ and $t\in [s,\infty[_{\tau}$ we have the local Lipschitz estimate
\begin{equation}\label{lipschitz-inq}
\vertiii{ \Phi_{s,t}(\mu)-\Phi_{s,t}(\eta)}_{\tiny V}\leq c(\mu,\eta)~e^{-b(t-s)}~\vertiii{ \mu-\eta}_{\tiny V}.
\end{equation}
For any $(\mu,H)\in (\Pa_V(E)\times \Ba_{0,V}(E))$ there exists some finite constant $c_H(\mu)<\infty$  such that for any $s\geq 0$ and $t\in [s,\infty[_{\tau}$ we have
\begin{equation}\label{kr-inq}
\vertiii{\frac{Q_{s,t}}{\mu_sQ_{s,t}(1)}-T_{s,t}^{\mu,H}}_V\leq ~c_H(\mu)~ e^{-b (t-s)}.
\end{equation}
For continuous time semigroups, the above estimates also hold for any continuous time indices $s\leq t$ as soon as for any $r\geq r_0$ there exists some $\overline{r}\geq r$ such that $\inf_{\delta\in [0,\tau]}\jmath_{r,\overline{r}}(\delta)>0$.

\end{theo}

The proof of Theorem~\ref{theo-intro} is based on discrete time type $V$-norm operator contraction techniques combining the geometric drift condition stated in the l.h.s. of (\ref{ref-V-theta-intro}) with the local minorization condition stated in (\ref{min-max-q-vr}). The condition $\inf_{\delta\in [0,\tau]}\jmath_{r,\overline{r}}(\delta)>0$ is a technical condition only made for continuous time semigroups to ensure that (\ref{lipschitz-inq}) and
(\ref{kr-inq})  also hold for continuous time indices.

Theses regularity conditions  are rather flexible as we will now explain.

Absolutely continuous integral operators arise in a natural way in discrete time settings~\cite{dm-2000,dm-04,douc-moulines-ritov,whiteley} and in the analysis of continuous time elliptic diffusion absorption models~\cite{aronson,ferre-phd,ferre-ldp,stroock}. In connection to this, two-sided estimates for stable-like processes are provided in~\cite{bogdan,kopnova,song,wang}.
Two sided Gaussian estimates can also be obtained  for some classes of degenerate diffusion processes of rank 2, that is when the Poisson brackets of the first order span the whole space~\cite{konakov}. 
  This class of diffusions includes frictionless Hamiltonian kinetic models.  
  
Diffusion density estimates can be extended to sub-Markovian semigroups using the multiplicative functional methodology developed in~\cite{dm-sch-2}. Whenever the trajectories of these diffusion flows, say  $t\mapsto X_t(x)$, where $x\in E$ is the initial position, are absorbed on the smooth boundary $\partial E$ of a open connected domain $E$, for any $\tau>0$ the densities $q_{\tau}(x,y)$ of the sub-Markovian semigroup $Q_{\tau}$  (with respect to the trace of the Lebesgue measure on $E$)  associated with the non absorption event are null at the boundary. Nevertheless, whenever these densities are positive and continuous on the open set $E^2$ for some $\tau>0$, they are uniformly positive and bounded on any compact subset of $E$;  thus condition (\ref{min-max-q-vr}) is satisfied. 
 
 In this context, whenever $T(x)$ stands for first exit time from $E$ and $T_r(x)$ the first exit time from the compact level set $\Va(r)\subset E$ starting from $x\in \Va(r)$, for any $\delta\in [0,\tau]$ and $r_+>r$ we have the estimate
\begin{eqnarray*}
Q_{\delta}(1_{\Va(r_+)})(x)&:=&\EE\left(1_{\Va(r_+)}(X_{\delta}(x))~1_{T(x)>\delta}\right)\geq \PP\left(T_{r_+}(x)>\delta\right)\geq \PP\left(T_{r_+}(x)>\tau\right). 
\end{eqnarray*}
In this context, we have
\begin{equation}\label{ref-over-r}
\inf_{x\in \Va(r)}\PP\left(T_{r_+}(x)>\tau\right)>0\Longrightarrow \inf_{\delta\in [0,\tau]}\inf_{\Va(r)}Q_{\delta}(1)\geq \inf_{\delta\in [0,\tau]}\jmath_{r,r_+}(\delta)>0.
 \end{equation}
Whenever the interior $E_+:=\Va(r_+)^{\circ}$ is a connected domain,   the l.h.s. estimate in (\ref{ref-over-r}) is met as soon as the sub-Markovian semigroup $Q^{+}_{\tau}$  associated with the non absorption event at the boundary $\partial E_+$ has a strictly positive continuous density $(x,y)\in E_+^2\mapsto q^{+}_{\tau}(x,y)$. To check this claim, observe that for any $x\in \Va(r)$  we have
\begin{eqnarray*}
 \PP\left(T_{r_+}(x)>\tau\right)=Q^{+}_{\tau}(1)(x)\geq Q^{+}_{\tau}(1_{\Va(r)})(x)&=&\int~q^{+}_{\tau}(x,y)~1_{\Va(r)}(y)~\nu_{\tau}(dy)\\&\geq &\nu_{\tau}(\Va(r))~\inf_{\Va(r)^2}q^{+}_{\tau}>0.
\end{eqnarray*}
It is out of the scope of this article to review the different classes of absolutely continuous operators and related two-sided Gaussian estimates arising in  the analysis of continuous time elliptic diffusion and particle absorption models. For a more thorough discussion on this topic we refer to the series of reference pointers presented above.

Needless to say that the design of Lyapunov functions is a crucial and challenging problem in the stability analysis of positive semigroups.
We have chosen to concentrate our review on presenting practical and general principles for designing Lyapunov functions.

\subsection{Time homogenous models}\label{time-homogenous-sec}
For time homogeneous models we use the notation
$$
(\Phi_t,Q_t,\overline{Q}_t):=(\Phi_{0,t},Q_{0,t},\overline{Q}_{0,t}).
$$
As expected for time homogeneous semigroups
a variety of results follow almost immediately from the estimates obtained in Theorem~\ref{theo-intro}. Following~\cite{dhj-21} (cf. for instance Section 4.1 and Section 4.3), these results include the existence of  an unique leading eigen-triple  \begin{equation}\label{def-eigen-triple}
 (\rho,\eta_{\infty},h)\in (\RR\times\Pa_V(E)\times \Ba_{0,V}(E))\quad \mbox{\rm with}\quad\eta_{\infty}(h)=1
 \end{equation}
 in the sense that for any $t\in \Ta$ we have
 \begin{equation}\label{intro-rho-h}
Q_t(h)=e^{\rho t}~h\quad \mbox{\rm and}\quad\eta_{\infty}Q_t=e^{\rho t}~\eta_{\infty}
\quad \mbox{\rm or equivalently}\quad \Phi_t(\eta_{\infty})=\eta_{\infty}.
 \end{equation}
The eigenfunction $h$ is sometimes called the ground state and the fixed point measure
$\eta_{\infty}$ the quasi-invariant measure. For any $x\in E$ we also have the product series formulation
$$
0<h(x):=\prod_{n\geq 0}~\left\{1+e^{-\rho \tau}
\left[\Phi_{n\tau}(\delta_x)({Q}_{\tau}(1))-\Phi_{n\tau}(\eta_{\infty})({Q}_{\tau}(1))\right]\right\}.
 $$

In this context, choosing $(\mu,H)=(\eta_{\infty},h)$ in (\ref{kr-inq}), we readily check that
$$
T^{\eta_{\infty},h}_{s,s+t}(f)= T(f):=\frac{h}{\eta_{\infty}(h)}~\eta_{\infty}(f)
\quad \mbox{\rm and}\quad
\vertiii{e^{-\rho t}~Q_{t}-T}_{V}\leq 
c_h(\eta_{\infty}) ~e^{-b t}.
 $$
In terms of the Boltzmann-Gibbs transformation $\Psi_{h}$ introduced in (\ref{def-Psi-H}), for any $\eta\in \Pa_V(E)$ we have the conjugate formulae
\begin{equation}\label{link-h-phi}
\Psi_{h}(\Phi_{t}(\eta))=
\Psi_{h}(\eta)P^h_t
\end{equation}
 with the Doob $h$-transform of $Q_t$ defined by the Markov semigroup
 $$
 P^h_t~:~f\in \Ba_{V^h}(E)\mapsto P^h_t(f):=e^{-\rho t}\frac{1}{h}~Q_{t}( hf)\in \Ba_{V^h}(E).
 $$

Observe that
$$
\eta_{\infty}=\Phi_t(\eta_{\infty})
\Longleftrightarrow\quad
 \eta^h_{\infty}:=\Psi_h(\eta_{\infty})=\eta^h_{\infty}P^h_t.
$$
The Markov semigroup $P^h_t$ is sometimes called the transition semigroup of the $h$-process, a.k.a. the process evolving in the ground state.

 We further assume that $Q_t$ is a sub-Markov semigroup of {\em self-adjoint operators} on $\LL_2(\nu)$ with respect to some locally finite measure $\nu$ on $E$. In addition, there exists an orthonormal basis $(\varphi_n)_{n\geq 1}$ associated with a decreasing sequence of eigenvalues $\rho_n\leq 0$ such that
\begin{equation}\label{spectral-Q}
Q_t(x,dy)=\sum_{n\geq 1}~e^{\rho_n t}~\varphi_n(x)~\varphi_n(y)~\nu(dy).
\end{equation}
In this context, the formulae (\ref{intro-rho-h}) are satisfied with the parameters
$$(\rho,h)=(\rho_1,\varphi_1)
\quad \mbox{\rm and}\quad
\eta_{\infty}(dx)=\Psi_h(\nu)(dx):=\frac{1}{\nu(h)}~h(x)~\nu(dx).
$$ 
Note that in this case $h$ has unit norm $\nu(h^2)=1$. 
The spectral resolution (\ref{spectral-Q}) yields 
for any $t\geq 0$ and $f\in \LL_2(\nu)$  the following decomposition
\begin{equation}\label{lambda-0-Ka}
e^{-\rho t}Q_{t}(f)(x)-\frac{h(x)}{\eta_{\infty}(h)}~ \eta_{\infty}(f)=\sum_{n\geq 2} e^{\rho^h_nt}~\varphi_n(x)~\nu(\varphi_nf) \quad \mbox{with}\quad
\rho^h_n=\rho_n-\rho_1.
\end{equation}

This yields the following result.
\begin{prop} 
For any time horizon $t\geq 0$  and any $f\in \LL_2(\nu)$ 
 we have the exponential estimates \begin{equation}\label{expo-decays-Ka}
\left\Vert e^{-\rho t}Q_{t}(f)-\frac{h}{\eta_{\infty}(h)}~\eta_{\infty}(f) \right\Vert_{\LL_2(\nu)}\leq e^{\rho^h_2t}~\left(\nu(f^2)-\nu(hf)^2
\right)^{1/2}.
\end{equation}
\end{prop}

Whenever $Q_t$ is a positive semigroup of {\em self-adjoint operators} on $\LL_2(\nu)$ the Doob $h$-transform $P^h_t$  is a semigroup of self-adjoint operators on $\LL_2(\eta_{\infty}^h)$ and we  have the following spectral decomposition
\begin{lem} 
For any $t\geq 0$ and $f\in \LL_2(\eta_{\infty}^h)$ we have
\begin{equation}\label{ref-spectral-dec}
P^h_t(x,dy)=\eta_{\infty}^h(dy)+\sum_{n\geq 2}~e^{\lambda_n t}~h_n(x)~h_n(y)~\eta_{\infty}^h(dy)
\end{equation}
with the $\LL_2(\eta_{\infty}^h)$   orthonormal basis $(h_n)_{n\geq 2}$ defined for any $n\geq 2$ by
$$
h_n:=\varphi_n/h\quad \mbox{and}\quad \lambda_n=\rho_n-\rho_1< 0
\quad\mbox{and}\quad
\eta_{\infty}^h=\Psi_{h^2}(\nu).
$$
\end{lem}
Note that the density of the integral operator $P^h_t(x,dy)$ w.r.t. $\eta_{\infty}^h(dy)$ is given by
\begin{equation}\label{ref-p-h}
p^h_t(x,y)=e^{-\rho_1 t}~\frac{q_{t}(x,y)}{h(x)h(y)}=1+\sum_{n\geq 2}~e^{\lambda_n t}~h_n(x)~h_n(y).
\end{equation}

We further assume that  $h\in \Ba_0(E)$ and  $P^h_t$ is ultra contractive, in the sense that
 for any $t>0$ we have
\begin{equation}\label{ref-ultra-contract-sinus}
\vertiii{P^h_t}_{\LL_2(\eta_{\infty}^h)\mapsto \LL_{\infty}(\eta_{\infty}^h)}= e^{-\rho_1 t}~ \sup_{(x,y)\in E^2}\frac{q_t(x,y)}{h(x)h(y)}=\sup_{(x,y)\in E^2}p^h_t(x,y)<\infty.
\end{equation}
\begin{prop}
Assume that $\nu(E)<\infty$ and $h\in \Ba_0(E)$. In addition, for any $t>0$ (\ref{ref-ultra-contract-sinus}) holds  and the mapping
$
x\mapsto \int~p^h_t(x,y)~\nu(dy)
$
is upper semi-continuous and locally lower bounded. In this situation, the function 
$V:=1/h\in \Ba_{\infty}(E)$
 and for any $t>0$ we have $ Q_{t}(V)/V\in \Ba_0(E)$.
 In addition,  for any $t>0$  we have \begin{equation}\label{ex-dirichlet}
 Q_{t}(V)/V
 \leq c_t/V^2 \in \Ba_0(E). 
\end{equation}
\end{prop}

\subsection{Sub-Markov semigroups}\label{sub-Markov-sec-intro}
Sub-Markov semigroups  are prototype-based models of positive integral operators. In time homogeneous settings, these stochastic models are defined in terms of a  stochastic flow $X_{t}(x)$ evolving on some metric Polish space $(\Ea,d)$, some non negative absorption potential function $U$ on some non necessarily bounded Borel subset $E\subset \Ea$. For a given $x\in E$ we denote by
  $T(x)$ the exit time of the flow $X_t(x)$ from $E$.
  
 We associate with these objects, the sub-Markov semigroup $Q^{[U]}_{t}$ defined for any $f\in \Ba_b(E)$ and $x\in E$ by
\begin{equation}\label{Q-U-gen}
Q^{[U]}_{t}(f)(x)=\EE\left(f(X_{t}(x))~1_{T(x)>t}~\exp{\left(-\int_0^t U(X_{s}(x))ds\right)}\right).
\end{equation}
The above model can be interpreted as the distribution of a stochastic flow evolving in an absorbing medium with hard and soft obstacles. Before killing, the flow starts at $x\in E$ and evolves as $X_t(x)$. Then, it is killed at rate $U$ or as soon as it exits the set $E$.
In the case $E=\Ea$, the flow cannot exit the set $E$ and it is only killed at rate $U$. This situation is sometimes referred a sub-Markov semigroup with soft obstacles represented by the absorbing potential function $U$ on $E$. When the flow may exit the set $E\subset\Ea$, the complementary subset $C:=\Ea-E$ is interpreted as an hard obstacle, a.k.a. an infinite energy barrier.

We illustrate the $V$-positive semigroup analysis developed in this article through three typical examples of solvable sub-Markov semigroups arising in physics and applied probability.

\subsubsection{The harmonic oscillator}\label{one-d-harmonic}
Consider the case  $E=\Ea=\RR$, and let $X_t(x)=B_t(x)$ be a Brownian motion starting at $x\in\RR$ and let $U(x)=x^2/2$. In this situation, the semigroup $Q^{[U]}_{t}=Q_t$ defined in  (\ref{Q-U-gen}) coincides with the one dimensional harmonic oscillator. For any $t>0$, the integral operator $Q_t$ has a continuous density w.r.t. the uniform measure $\nu$ on $E$ given by 
\begin{equation}\label{q-harmonic}
q_t(x,y)=\sum_{n\geq 1}~e^{\rho_nt}~\varphi_n(x)\varphi_n(y)
\end{equation}
with the $\LL_2(\nu)$ orthonormal basis eigenstates
$$
\varphi_n(x)=(2^{n-1} (n-1)!\sqrt{\pi})^{-1/2}~e^{-x^2/2}~ \HH_{n-1}(x)
$$
associated with the eigenvalues
$$
 \rho_n=-(n-1/2)\quad \mbox{and the Hermite polynomials}\quad \HH_n(x)=(-1)^n~e^{x^2}~\partial^ne^{-x^2}.
$$
In this context, the eigenstate  associated with the top eigenvalue $\rho=\rho_1=-1/2$ is given by the harmonic function 
\begin{equation}\label{h-harmonic}
h(x)=\varphi_1(x)=\pi^{-1/4}~e^{-x^2/2}.
\end{equation}
 The spectral resolution of  integral operator $P^h_t(x,dy)$ and its density $p^h_t(x,y)$ with respect to the invariant measure
 $$
  \eta_{\infty}^h(dy)=\frac{1}{\sqrt{\pi}}~e^{-y^2}~dy
 $$ are given as in (\ref{ref-spectral-dec}) and (\ref{ref-p-h}) with $\LL_2(\eta^h_{\infty})$ orthonormal basis defined for any $n\geq 2$ by
 $$
h_n=(2^{n-1} (n-1)!)^{-1/2}~ \HH_{n-1}\quad\mbox{\rm and}\quad\rho^h_n=\rho_n-\rho_1=-(n-1).
 $$
In this context, the $h$-process is given by the Ornstein-Uhlenbeck diffusion
\begin{equation}\label{h-proc-harmonic}
dX^h_t(x)=\partial \log{h(X^h_t(x))}~dt+dB_t=-X^h_t(x)~dt+dB_t.
\end{equation}
In the above display,  $B_t=B_t(0)$ stands for the one dimensional Brownian motion starting at the origin. 
The conjugate formula 
\begin{equation}\label{conjugate-formula}
Q_t(hf)/Q_t(h)=P^h_t(f)\Longleftrightarrow
Q_t(f)=e^{\rho t} h~P^h_t(f/h)
\end{equation}
 yields the following proposition.
\begin{prop}\label{solv-harmonic-P}
For any time horizon $t\geq 0$ we have 
$$
Q_t(x,dy)=\frac{1}{\sqrt{\cosh(t)}}~\exp{\left(-\frac{x^2}{2}~p_t\right)}~\frac{1}{\sqrt{2\pi p_t}}~\exp{\left(-\frac{(y-m_t(x))^2}{2p_t}\right)}~dy
$$
with the mean and variance parameters $(m_t(x),p_t)$ defined by
$$
m_t(x)={x}/{\cosh(t)}\quad \mbox{and}\quad
p_t=\tanh(t).
$$
\end{prop}
The proof of the above proposition is a direct consequence of the conjugate formula (\ref{conjugate-formula}), thus it is provided in the appendix, on page~\pageref{solv-harmonic-P-proof}.  

Choosing 
$
V(x)=1+\vert x\vert^n$, for some $n\geq 1$, we readily check that
\begin{equation}\label{harm-lyap-proof}
V\in \Ca_{\infty}(E)\quad \mbox{\rm and}\quad
Q_t(V)/V\leq v_t~Q_t(1)~\in \Ca_0(E)
\end{equation}
where $v_t$ is a constant depending only on $t$.

\subsubsection{The half-harmonic oscillator}\label{half-oscillator-sec}
Consider the case  $E=]0,\infty[\subset \Ea=\RR$, and let $X_t(x)=B_t(x)$ be a Brownian motion starting at $x\in\Ea$ and let $U(x)=x^2/2$. In this situation, the semigroup $Q^{[U]}_{t}=Q_t$ defined in  (\ref{Q-U-gen}) coincides with the harmonic oscillator with an infinite barrier at the origin $\partial E=\{0\}$ (a.k.a. the half-harmonic oscillator).
Using the fact that
$$
 e^{x^2/2}~\frac{1}{2}~\partial^2\, e^{-x^2/2}=U(x)-1/2
$$
by an exponential change of probability measure (cf. for instance Section 18.3 in~\cite{dm-penev-2017}) we have the conjugate formula
  \begin{eqnarray*}
Q_{t}(f)(x)
&=&e^{-t/2}~e^{-x^2/2}~\EE\left(f(Y_t(x))~e^{Y_t(x)^2/2}~1_{T^Y(x)>t}\right)
\end{eqnarray*}
with the Ornstein-Uhlenbeck diffusion
\begin{equation}\label{h-proc-harmonic-demi}
\begin{array}{l}
dY_t(x)=-Y_t(x)~dt+dB_t\\
\\
\mbox{\rm and}
\quad
T^Y(x):=\inf\left\{t\geq 0~:~Y_{t}(x)\in \partial E\right\}\quad \mbox{\rm with}
\quad \partial E=\{0\}.
\end{array}
\end{equation}
Note that the stochastic flow $Y_t(x)$ coincides with the $h$-process of the harmonic oscillator discussed in (\ref{h-proc-harmonic}). 
Thus, by reflection arguments we have
\begin{eqnarray*}
Q^Y_t(f)(x)&:=&\EE\left(f(Y_t(x))~1_{T^Y(x)>t}\right)=
\EE\left(f(B_{\sigma_t}(\epsilon_tx))~1_{T(\epsilon_tx)>t}\right)\\
&=&
\int_0^{\infty}~f(y)~q^Y_t(x,y)~dy\quad\mbox{\rm with} \quad q^Y_t(x,y):=(r_t(x,y)-r_t(x,-y)).
\end{eqnarray*}
In the above display, $(\epsilon_t,\sigma_t)$ stands for the parameters $$
\left(\epsilon_t,\sigma_t\right):=\left(e^{-t},\sqrt{\frac{1-\epsilon_t^2}{2}}\right)
\quad \mbox{\rm and}\quad
r_t(x,y)=\frac{1}{\sqrt{2\pi\sigma^2_t}}~\exp{\left(-\frac{1}{2\sigma^2_t}\left(y-\epsilon_t x\right)^2\right)}.
$$
This yields the following proposition.
\begin{prop}\label{prop-half-harmonic}
For any $t>0$ and $x\in E$ the normalized semigroup $\overline{Q}_t$ defined in (\ref{def-Phi-s-t-intro-0}) is given by
$$
\overline{Q}_t(x,dy)=\frac{\sinh{(y~m_t(x))}}{ \PP\left(0\leq Z\leq m_t(x)/\sqrt{p_t} \right)}~\displaystyle\times \frac{1}{\sqrt{2\pi p_t}}~\exp{\left(-\frac{y^2+m_t(x)^2}{2p_t}\right)}~\nu(dy).
$$
In the above display,  $\nu(dy):=1_{[0,\infty[}(y)~dy$ stands for the trace of the Lebesgue measure on the half-line, $Z$ is a centered Gaussian variable with unit variance and  $(m_t(x),p_t)$ are the mean and variance parameters defined in Proposition~\ref{solv-harmonic-P}. In addition,  the total mass function $Q_t(1)(x)$ is given by the formula
$$
 Q_t(1)(x)= 2~\frac{e^{-\frac{x^2}{2}\,p_t}}{\sqrt{\cosh(t)}}\times \PP\left(0\leq Z\leq m_t(x)/\sqrt{p_t}  \right)\in \Ca_0(E).
$$
\end{prop}
The proof of the above proposition follows the same lines of arguments as the proof of Proposition~\ref{solv-harmonic-P};  it is provided in the appendix, on page~\pageref{prop-half-harmonic-proof}.

Choosing 
$
V(x)=x^n+1/x$, for some $n\geq 1$, we readily check that
\begin{equation}\label{half-harm-lyap}
V\in \Ca_{\infty}(E)\quad \mbox{\rm and}\quad
Q_t(V)/V\leq c_t/V~\in \Ca_0(E).
\end{equation}
The proof of the above estimate follows elementary but lengthly calculations, thus it is provided in the appendix on page~\pageref{half-harm-lyap-proof}.

For any $t>0$, the integral operator $Q_t$ has a continuous density w.r.t. the uniform measure $\nu$ on $E$ given by 
 \begin{eqnarray*}
q_t(x,y)&=&\sum_{n\geq 1}~e^{\rho_{n}t}~\varphi_{n}(x)~\varphi_{n}(y)
\end{eqnarray*}
with the $\LL_2(\nu)$ orthonormal basis eigenstates
$$
\varphi_n(x)=\sqrt{2}~(2^{2n-1} (2n-1)!\sqrt{\pi})^{-1/2}~e^{-x^2/2}~ \HH_{2n-1}(x)
$$
 associated with the eigenvalues
 $$
  \rho_n=-((2n-1)+1/2).
$$
In this context, the eigenstate  associated with the top eigenvalue $\rho=\rho_1=-3/2$ is given for any $x\in ]0,\infty[$ by the harmonic function 
$$
h(x)=\varphi_1(x)=2\pi^{-1/4}~x~e^{-x^2/2}=h_0(x)~\HH_1(x)
$$
with the ground state $h_0$  of the harmonic oscillator discussed in (\ref{h-harmonic}).
Note that $h$ coincides with the restriction on $]0,\infty[$ of the first excited state of the harmonic-oscillator (negative on $]-\infty,0]$ and crossing the origin at $x=0$). 

 The spectral resolution of  integral operator $P^h_t(x,dy)$ and its density $p^h_t(x,y)$ with respect to the invariant measure
 $$
  \eta_{\infty}^h(dy)=\frac{4}{\sqrt{\pi}}~y^2~e^{-y^2}~1_{]0,\infty[}(y)~dy
 $$ are given for any $x,y\in ]0,\infty[$ as in (\ref{ref-spectral-dec}) and (\ref{ref-p-h}) with $\LL_2(\eta^h_{\infty})$ orthonormal basis defined for any $n\geq 2$ and $x\in ]0,\infty[$ by the odd Hermite functions
 $$
h_n(x)=(2^{2n} (2n-1)!)^{-1/2}~\HH_{2n-1}(x)/x
\quad\mbox{\rm and}\quad\rho^h_n=-2(n-1).
 $$
 In this context, the $h$-process is given by the diffusion
\begin{equation}\label{h-diff-positive}
dX^h_t(x)=\partial \log{h(X^h_t(x))}~dt+dB_t=\left(\frac{1}{X^h_t(x)}-X^h_t(x)\right)~dt+dB_t.
\end{equation}

\subsubsection{The Dirichlet heat kernel}\label{dhk-sec}
Let $X_t(x)=B_t(x)$ be a Brownian motion starting at $x\in E:=]0,1[\subset \Ea:=\RR$ and $T(x)$ be the first time $t\geq 0$ the process $B_t(x)\in \partial E:=\{0,1\}$. Choosing $U=0$ in (\ref{Q-U-gen}), the semigroup $Q^{[U]}_{t}=Q_t$ takes the following form
$$
Q_t(f)(x):=\EE(f(B_t(x))~1_{T(x)>t}).
$$
For any $t>0$, the integral operator $Q_t$ has a continuous density w.r.t. the uniform measure $\nu$ on $E$ given by 
the Dirichlet heat kernel 
\begin{equation}\label{intro-example}
q_t(x,y)=\sum_{n\geq 1}~e^{\rho_nt}~\varphi_n(x)\varphi_n(y)
\end{equation}
with the $\LL_2(\nu)$ orthonormal basis eigenstates
$$
\varphi_n(x)=\sqrt{2}~\sin{(n\pi x)}\quad \mbox{\rm associated with the eigenvalues}\quad \rho_n=-(n\pi)^2/2.
$$
 In this context, the eigenstate $h(x)=\varphi_1(x)=\sqrt{2}~\sin{(\pi x)}
$ associated with the top eigenvalue $\rho=\rho_1=-\pi^2/2$ is strictly positive except at the boundary $\{0,1\}$.
By removing the boundary,  
the semigroup $P^h_t$ of the process evolving in the ground state $h(x)$ on the open interval $E$  is a self-adjoint operators on $\LL_2(\eta_{\infty}^h)$ 
with
$$
\eta_{\infty}^h(dx)=h^2(x)~\nu(dx)=2~\sin^2{(\pi x)}~1_{E}(x)~dx.
$$
In addition, we  have the spectral decomposition (\ref{ref-spectral-dec})
with the $\LL_2(\eta^h_{\infty})$ orthonormal basis eigenstates
$$
h_n(x):=\sin{(n\pi x)}/\sin{(\pi x)}
$$
associated with the eigenvalues
$$
 \lambda_n=-\pi^2(n^2-1)/2< 0.
$$
Our next objective is to estimate the density $p^h_t(x,y)$ of the integral operator $P^h_t(x,dy)$ w.r.t. $\eta_{\infty}^h$ defined  in (\ref{ref-p-h}).
Recalling that $\vert \sin{(ny)}\vert\leq n\vert \sin{(y)}\vert$, for any $n\geq 1$ and $y\in\RR$, for any $x\in E$ we have the diagonal estimate
$$
p^h_t(x,x)-1=\sum_{n\geq 2}~e^{\rho_n^h t}~h_n(x)^2
$$
with
$$
h_n(x)^2=\left(\frac{\sin{(n\pi x)}}{\sin{(\pi x)}}\right)^2\leq n^2\quad\mbox{\rm so that condition (\ref{ref-ultra-contract-sinus}) is satisfied.}
$$
Observe that the function 
$$
V~:~x\in E\mapsto V(x):=\sqrt{2}/h(x)\in [1,\infty[
$$
is locally bounded with compact level sets given for any $0<\epsilon\leq 1$ by the formulae
$$
\Ka_{\epsilon}:=\{x\in ]0,1[~:~ V(x)\leq 1/\epsilon\}=\{x~:~\sin{(\pi x)}\geq \epsilon\}\subset E.
$$

In any dimension we can use the intrinsic ultracontractivity to produce a Lyapunov function $V$. Let $E$ be a bounded domain of $\RR^n$ 
for some $n \ge 1$ and assume that it is a $C^{1, \alpha}$ domain for some $\alpha > 0$. Denote by $q_t(x,y)$ the Dirichlet heat kernel on $E$. By \cite{Ouhabaz-Wang} one has
$$ q_t(x,y) \le c_t ~d(x, \partial E) d(y, \partial E)$$
for some constant $c_t$ independent on $x$ and $y$. Here $d(x, \partial E)$ denotes the distance from $x$ to the boundary of $E$. Set $V(x) = \frac{1}{d(x, \partial E)}$. The above  intrinsic ultracontractivity implies
$$ Q_t(V)(x) = \int_E q_t(x,y) V(y) dy \le c_t |E| ~ d(x, \partial E)$$
which in turn gives  ${Q_t(V)}/{V} \le {c_t |E|}/{ V^2} \in \Ba_0(E)$, where $\vert E\vert$ stands for the volume of the bounded set $E$.

\section{Lyapunov design principles}\label{lyapunov-principles-sec}

The aim of this section is to present some general principles to construct Lyapunov functions for positive semigroups. Section~\ref{foster-lyap} provides equivalent 
formulations of the Lyapunov 
 condition in (\ref{ref-V-theta-intro}) encountered in the literature in terms of exhausting sequences of compact level sets. This section also presents simple ways to design Lyapunov functions for sub-Markov semigroups on normed spaces in terms of their generators. Section~\ref{sg-domination} presents some principles to construct Lyapunov functions for positive semigroups dominated by semigroups with known Lyapunov functions. Section~\ref{sg-conjugacy} is dedicated to the design of Lyapunov functions for conjugate semigroups. All the principles discussed in this section are illustrated in Section~\ref{ricc-sec} as well as in Section~\ref{conditional-diff-sec} in the context of conditional diffusions.

\subsection{Foster-Lyapunov conditions}\label{foster-lyap}

For time homogeneous models $Q_{s,s+t}:=Q_t$,  the l.h.s.
 condition in (\ref{ref-V-theta-intro}) takes the form $Q_{\tau}(V)/V\leq \Theta_{\tau}\in \Ba_0(E)$. In terms of the compact sets
$
K_{\epsilon}:=\{\Theta_{\tau}\geq \epsilon\}$,  the l.h.s
Lyapunov condition in (\ref{ref-V-theta-intro})
yields for any $\tau>0$ the estimate
\begin{equation}\label{ref-epsilon-V-K}
Q_{\tau}(V)(x)\leq \epsilon ~V(x)+1_{K_{\epsilon}}(x)~c_{\epsilon}
\end{equation}
for any $\epsilon > 0$ with the parameter $c_{\epsilon}:=
\sup_{K_{\epsilon}}(V\Theta_{\tau})<\infty.$ 
This implies that for any $n\geq 1$ we have
 \begin{equation}\label{ref-epsilon-V-K-n}
Q_{\tau}(V)(x)\leq \epsilon_n~V(x)+1_{K_{\epsilon_n}}(x)~c_{\epsilon_n}
\end{equation}
where $K_{\epsilon_n}\subset E$ stands for some increasing sequence of compacts sets and $c_{\epsilon_n}$ some finite constants, indexed by a decreasing 
sequence of parameters $\epsilon_n\in [0,1]$ such that $\epsilon_n\longrightarrow 0$ as ${n\rightarrow\infty}$. In the reverse angle, assume that $Q_{\tau}(V)/V$ is locally lower bounded and lower semicontinuous. In this situation, condition (\ref{ref-epsilon-V-K-n}) ensures that $Q_{\tau}(V)/V\in \Ba_0(E)$ for any $\tau>0$. Indeed,
 for any $\delta>0$, there exists some $n\geq 1$ such that $\epsilon_n< \delta$ and we have
$$
\{Q_{\tau}(V)/V\geq \delta\}\subset\{Q_{\tau}(V)/V> \epsilon_n\}\subset K_{\epsilon_n}.
$$
Since $\{Q_{\tau}(V)/V\geq \delta\}$ is a closed subset of a compact set it is also compact.

More generally, whenever (\ref{ref-epsilon-V-K-n}) is met for  some exhausting sequence of compact sets $K_{\epsilon_n}$, in the sense that for any compact subset $K\subset E$ there exists some $n\geq 1$ such that $K\subset K_{\epsilon_n}$ we have
 $$
 \inf_K Q_{\tau}(V)/V\geq  \inf_{K_{\epsilon_n}} Q_{\tau}(V)/V \geq \epsilon_n.
 $$
 This ensures that the function $Q_{\tau}(V)/V$ is necessarily locally lower bounded. In this situation,  we have $Q_{\tau}(V)/V\in \Ba_0(E)$ as soon as 
 $Q_{\tau}(V)/V$ is  lower semicontinuous.
 
 Notice that the sub-level set $\Va(r):=\{V\leq r\}$ of the Lyapunov function $V\in \Ba_{\infty}(E)$ and the $\epsilon$-super-level sets $K_{\epsilon}:=\{\Theta_{\tau}\geq \epsilon\}$ of $\Theta_{\tau}\in \Ba_0(E)$ are equivalent compact exhausting sequences, in the sense that for any $r\geq 1$ we have
$$
\Va(r)\subset K_{\epsilon_r}\subset
\Va(r_{\epsilon})\quad \mbox{\rm with}\quad
\epsilon_r:=\inf_{\Va(r)}\Theta_{\tau}\quad \mbox{\rm and}\quad
r_{\epsilon}:=\sup_{K_{\epsilon_r}}V.
$$
 
Whenever $E$ is a locally compact Polish space,  the abstract sequence $C_n:=K_{\epsilon_n}$ in (\ref{ref-epsilon-V-K-n}) is automatically exhausting; that is, we have that $E=\cup_{n\geq 0}C_n$ with $C_n$ is included in the interior $C^{\circ}_{n+1}$ of the compact set $C_{n+1}$. To check this claim, observe that for any $n\geq 1$ there exists some $m_n\geq n$ such that
$$
 C_n\subset \{\Theta_{\tau}\geq \inf_{C_n}\Theta\}\subset C_{m_n}\subset \{\Theta_{\tau}\geq \inf_{C_{m_n}}\Theta\}.
$$
Thus,  the exhausting sequence $C_n$ is equivalent to the one defined by the super-level sets of $\Theta_{\tau}$.

The rather abstract condition (\ref{ref-epsilon-V-K-n}) is often presented in the literature as an initial condition 
to check on a case-by-case basis to analyze the stability property of time homogenous sub-Markov semigroups (see for instance ~\cite{ferre,guillin-nectoux-Wu}, as well as Section 17.5 in~\cite{dm-penev-2017} in the context of Markov semigroups and the references therein).

We end this section with a brief discussion on condition (\ref{ref-epsilon-V-K-n}) in the context of the sub-Markov semigroup  discussed in (\ref{Q-U-gen}). Note that this semigroup can be turned into a Markov semigroup by sending the killed process into a cemetery state, say $\Delta$, at the killing time. In this interpretation, functions on $E$ are extended to $E_{\Delta}=E\cup\{\Delta\}$ by setting $f(\Delta)=0$. More interestingly, whenever $E$ is locally compact  its topology coincides with the weak topology induced by $\Ca_0(E):=\Ba_0(E)\cap \Ca_b(E)$, and inversely (cf. Proposition 2.1 in~\cite{aliabad}).  In this context a continuous function $f$ vanishes at infinity 
if and only if its extension to 
the one point compactification (a.k.a.
Alexandroff compactification)
$E_{\Delta}:=E\cup\{\Delta\}$ (obtained by setting $f(\Delta)=0$) 
is continuous. For locally compact spaces, we also recall that the one point extension
$E_{\Delta}$ is compact. 

Whenever it exists, the generator $L^U$ of these sub-Markov semigroups $Q^{[U]}_t$ are defined on domain of functions $D(L^U)\subset \Ba_0(E)$. As expected, the analysis of this class of models in terms of generators  often requires to develop a sophisticated analysis taking into account the topological structure of the set $E$.
To the best of our knowledge, there is no simple sufficient condition to check (\ref{ref-epsilon-V-K-n}) in terms of these generators.

The situation is greatly simplified for sub-Markov semigroups with soft obstacles. When $E=\Ea$ is a locally compact normed space $(E,\Vert\point\Vert)$ we let $L$  be the generator of the flow $X_t(x)$. In this situation, the generator of the sub-Markov semigroup $Q^{[U]}_t$ is given by $L^U=L-U$. We further assume that $L$ and $L^U$ are defined on some common domain $\Da(L)\subset \Ba(E)$.
\begin{lem}[\cite{ferre}]\label{lem-LUK-proof}
Let $V,V_0\in \Da(L)$ be a couple of functions  such that $V,V_0\geq 1$ and
\begin{equation}\label{ref-W-V-1}
V(x)\longrightarrow_{\Vert x\Vert\rightarrow\infty}\infty\quad\mbox{\rm and}\quad
V(x)/V_0(x)\longrightarrow_{\Vert x\Vert\rightarrow\infty}\infty.
\end{equation}
In this situation, condition  (\ref{ref-epsilon-V-K-n})
 is satisfied as soon as there exists some finite constant $c_0<\infty$ such that
\begin{equation}\label{ref-M-M-LV}
 L^U(V_0)/V_0\leq c_0\quad\mbox{and}\quad
L^U(V)(x)/V(x)\longrightarrow_{\Vert x\Vert\rightarrow\infty}-\infty.
\end{equation}
\end{lem}
Note that in this context, the compact sets in (\ref{ref-epsilon-V-K-n}) are given for some sufficiently large radii $r_{\epsilon}>0$ by the closed balls:
\begin{equation}\label{ref-M-MK}
 \Ka_{\epsilon}=\overline{\BB}(0,r_{\epsilon}):=\{x\in E~: ~\Vert x\Vert\leq r_{\epsilon}\}. 
\end{equation}

\subsection{Semigroup domination}\label{sg-domination}
For a given $p\geq 1$  we clearly  have 
$$
V\in \Ba_{\infty}(E)\Longleftrightarrow V^p\in \Ba_{\infty}(E)\quad \mbox{\rm and}\quad
\Ba_{V^{1/p}}(E)\subset\Ba_{V}(E)\subset \Ba_{V^p}(E). 
$$ We say that a $V$-positive semigroup $Q_{s,t}$ is $p$-dominated by a collection of integral operators $\Qa_{s,t}$ on $\Ba_{V^p}(E)$ and we write $Q\ll_p \Qa$ as soon as for any non negative function $f\in \Ba_{V}(E)$ and any $s\leq t$  we have
$$
Q_{s,t}(f)\leq c_{t-s}(p)~\Qa_{s,t}(f^p)^{1/p}.$$
 To simplify notation, when  $p=1$ we write $Q\ll \Qa$ instead of $Q\ll_1 \Qa$.
Observe that
$$
Q\ll_p \Qa\quad\Longrightarrow\quad \forall s\leq t\qquad
(Q_{s,t}(V)/V)^p\leq c_{t-s}(p)^p~\Qa_{s,t}(V^p)/V^p.
$$
 This yields for any $\tau>0$ and $\theta_{\tau}\in \Ba_0(E)$ the Lyapunov estimate
\begin{equation}\label{ref-dom-il}
 \Qa_{s,s+\tau}(V^p)/V^p\leq \theta^p_{\tau}\Longrightarrow Q_{s,s+\tau}(V)/V\leq  c_{\tau}~\theta_{\tau}.
 \end{equation}
We illustrate the above domination property with 
  the Langevin diffusion flow $\Xa_t^{(a)}(z)=(X_t(z),Y_t(z))\in (\RR^{n}\times\RR^n)$ starting at
$z=(x,y)\in (\RR^{n}\times\RR^n)$ and defined  by the hypo-elliptic  diffusion 
\begin{eqnarray}
dX_t(z)&=&Y_t(z)/m~dt\nonumber\\
dY_t(z)&=&\left(a(X_t(z))-\gamma X_t(z)-\beta Y_t(z)/m\right)~dt+\sigma~dB_t.\label{langevin-intro-ref}
\end{eqnarray}
In the above display, $\sigma,\gamma,\beta,m>0$ stands for some parameters and $a$ some Lipschitz function on $\RR^n$, with $n\geq 1$.
Notice that when $a=0$, the flow $\Xa_t^{(0)}(z)$ resumes to an hypo-elliptic Ornstein-Ulhenbeck on $\RR^{2n}$.
Consider a bounded open connected domain $D\subset\RR^n$ and  set
$$
\forall z\in E:=D\times \RR^n\qquad T^{(a)}(z):=\inf\left\{t\geq 0~:~\Xa_t^{(a)}(z)\in \partial E\right\}.
$$ 
We associate with these objects, the sub-Markov semigroup defined for any $f\in \Ba_b(E)$ and $z=(x,y)\in E$ by
$$
\Qa^{(a)}_t(f)(z):=\EE\left(f(\Xa^{(a)}_t(z))~1_{T^{(a)}(z)>t}\right).
$$
In this situation, we have
\begin{equation}\label{langevin-ref-intro}
\sup_{D}a<\infty\quad \Longrightarrow\quad \forall p>1\qquad
\Qa^{(a)}\ll_p\Qa^{(0)}.
\end{equation}
The proof of the above assertion is a direct consequence of Girsanov's theorem and H\"older's inequality. For the convenience of the reader, a detailed proof is provided in the appendix on page~\pageref{langevin-ref-intro-proof}.

 To emphasize the role of the absorption in sub-Markov semigroups we return to the class of models discussed in (\ref{Q-U-gen}).  We let $P_{t}$ be the free evolution Markov semigroup associated with the stochastic flow $X_{t}(x)$.
Assume that $Q^{[U]}_{t}(1)\in \Ba_{0}(E)$ and 
\begin{equation}\label{Q-U-V-V}
\Vert Q^{[U]}_{t}(V)/V \Vert<\infty\quad \mbox{\rm for some $t>0$ and $V\in\Ba_{\infty}(E)$.}
\end{equation} 
Applying H\" older's inequality and choosing $V_p:=V^{1/p} \in\Ba_{\infty}(E)$ with $p>1$ we readily check the estimate
\begin{equation}\label{Q-U-V-V-2}
Q^{[U]}_{t}(V_p)/V_p\leq c_t(p)~Q^{[U]}_{t}(1)^{1-1/p}\in \Ba_{0}(E).
\end{equation} 
The next lemma provides several practical conditions to check the uniform estimate (\ref{Q-U-V-V}) for sub-Markov semigroups associated with soft obstacles.
\begin{lem}\label{lem-LU}
Consider the sub-Markov semigroup discussed in (\ref{Q-U-gen}) 
when $E=\Ea$ is a locally compact normed space $(E,\Vert\point\Vert)$. Assume that the generators $L$ and $L^U$ of the flows $P_t$ and $Q^{[U]}_t$ are defined on some common domain $\Da(L)\subset \Ba(E)$.
 In this situation, for any $V\in \Ba_{\infty}(E)\cap \Da(L)$ and parameter $a>0$ we have
 \begin{equation}\label{L-U-V-ref}
L^U(V)\leq -a V+c\quad \Longrightarrow\quad \forall t\geq 0\qquad\Vert Q^{[U]}_{t}(V)/V\Vert<\infty.
\end{equation}
Whenever $U\in \Ba_{\infty}(E)\cap \Da(L)$, for any $a_0\geq 0$ and $a_1\in\RR$ we have
 \begin{equation}\label{L-U-ref}
L(U)\leq a_0+a_1 U\quad \Longrightarrow\quad \forall t\geq 0\qquad \Vert Q^{[U]}_{t}(U)\Vert<\infty.
\end{equation}

\end{lem}
The proof of the above lemma follows essentially the same lines of arguments as the proof of Lemma~\ref{lem-LUK-proof}; thus it is provided in the appendix, on page~\pageref{lem-LU-proof}.

 Whenever $E=\Ea$  and the absorption potential function  $U$ is bounded, we have
 $P\ll Q^{[U]}\ll P$. In this context, {\em there is no hope to have that $Q^{[U]}_{t}(1)\in \Ba_{0}(E)$ for some $t>0$}. Nevertheless, for any $V\in \Ba_{\infty}(E)$ and any time horizon $t>0$ we have
$$
 Q^{[U]}_{t}(V)/V\in \Ba_0(E)\Longleftrightarrow  P_{t}(V)/V\in \Ba_0(E).
$$
In this situation, the design of Lyapunov functions $V$ satisfying (\ref{ref-V-theta-intro}) or equivalently Foster-Lyapunov conditions of the form (\ref{ref-epsilon-V-K-n}) is equivalent  to the problem of finding a Lyapunov function for the Markov semigroup $P_t$.

 Whenever $P_t$ is stable, in the sense that it has a Lyapunov $V\in \Ba_{\infty}(E)$ such that $P_t(V)/V\in \Ba_0(E)$ for some $t>0$, then the domination property $Q^{[U]}\ll P$  yields automatically a Lyapunov function for $Q^{[U]}_t$.

 Whenever $P_t$ is not necessarily stable but we have $ \Vert P_t(V)/V\Vert<\infty$  for some $t>0$ and $V\in \Ba_{\infty}(E)$, applying (\ref{Q-U-V-V-2})  the domination property $Q^{[U]}\ll P$  
ensures that for any $p>1$ we have $V_p:=V^{1/p} \in\Ba_{\infty}(E)$ and
 $$
Q_t^{[U]}(1)\in \Ba_0(E)\quad\Longrightarrow\quad 
 Q^{[U]}_{t}(V_p)/V_p\in \Ba_{0}(E).
 $$
 
 Last, but not least, note that the above discussion extends without difficulties to time varying models.

 \subsection{Some conjugacy principles}\label{sg-conjugacy}
For any given $V\in \Ba_{\infty}(E)$, observe that for any positive function $H$,
$$
H\in \Ba_{0,V}(E)\Longleftrightarrow V^H:=V/H\in \Ba_{\infty}(E).
$$
Thus,  $\Qa_t$ is a $V$-positive semigroup on $\Ba_V(E)$ if and only if the $H$-conjugate semigroup $\Qa_t^H(f):=\Qa_t(fH)/H$ is
a $V^H$-positive semigroup on $\Ba_{V^H}(E)$. In this situation, any semigroup  $Q\ll\Qa^H$ dominated by $\Qa^H$ yields for any $s\geq 0$ and $t>0$ the Lyapunov estimate
$$
Q_{s,s+t}(V^H)/V^H\leq c_t~\Qa_t(V)/V\in \Ba_{0}(E). 
$$
To get one step further, observe that
$$
\Qa_t(V)/V=\Qa_t(1)~\overline{\Qa}_t(V)/V.
$$
In this notation, for any $H\in \Ba_{0,V}(E)$  and any $V$-positive semigroup $\Qa_t$ on $\Ba_V(E)$ such 
$\Qa_t(1)\in \Ba_0(E)$ and $\vertiii{\overline{\Qa}_{t}}_V<\infty$ we have
\begin{equation}\label{ref-Q-Qa-1-H}
Q\ll\Qa^H
\Longrightarrow
Q_{s,s+t}(V^H)/V^H\leq c_t~\Qa_t(1)\in \Ba_{0}(E). 
\end{equation}

We illustrate the above comparison principles with an elementary example. Let $\Ea:=\RR$
 and $W\in \Ba_{\infty}(\RR)$ be some non negative function. Consider the stochastic flow  $X^W_t(x)$  of a one-dimensional Langevin diffusion on $\Ea$ with generator 
\begin{equation}\label{ref-langevin-generator}
L(f)=\frac{1}{2}~e^{2W}\partial \left(e^{-2W}\partial f\right).
\end{equation}
We associate with a given  open connected interval  $E\subset \Ea$, the sub-Markov semigroup $Q_t$ on $\Ba_b(E)$ defined by
\begin{equation}\label{ref-langevin-one-D-examp}
Q_t(f)(x):=\EE(f(X^W_t(x))~1_{T_{\partial E}^{W}(x)>t})\quad \mbox{\rm with}\quad
T_{\partial E}^{W}(x):=\inf{\left\{t\geq 0~:~X^W_t(x)\in\partial E\right\}}.
\end{equation}
Observe that
\begin{equation}\label{ref-langevin-one-def-H}
H:=e^{-W}\Longrightarrow U:=H^{-1}\frac{1}{2}~\partial^2H=\frac{1}{2}~\left((\partial W)^2-\partial^2W\right).
\end{equation}
When $W=0$ the flow $X^0_t(x)=B_t(x)$ coincides with the Brownian flow $B_t(x)$ starting at $x$.
Thus, by a change of probability we check that
\begin{equation}\label{ref-Q-U-half-harmonic}
Q_t
=\Qa_t^{H}\quad \mbox{\rm with}\quad\Qa_t(f)(x):=
\EE\left(f(B_t(x))~1_{T^0_{\partial E}(x)>t}~e^{-\int_0^t~U(B_s(x))~ds}\right).
\end{equation}
Whenever $E=]0,1[$ the semigroup $\Qa_t$ is dominated by the Dirichlet heat kernel on $]0,1[$ discussed in Section~\ref{dhk-sec}. When $E=\RR$, respectively $E=]0,\infty[$, and $U(x)\geq c+\varsigma~x^2/2$, for some $c<\infty$ and $\varsigma>0$, the semigroup $\Qa_t$ is dominated by the harmonic oscillator discussed in Section~\ref{one-d-harmonic},  respectively the half-harmonic oscillator discussed in Section~\ref{half-oscillator-sec}. All of these dominating semigroups are completely solvable with $\Qa_t(1)\in \Ba_0(E)$ and known Lyapunov functions.

\section{Boundary problems}\label{boundary-sec}
Let $(\Ea,d)$ be a locally compact Polish space with a distinguished complete metric $d~:(x,y)\in E^2\mapsto d(x,y)\in \RR_+$. We recall that these metric spaces are complete $\sigma$-compact and locally compact metric spaces, thus
they have the Heine-Borel property, that is each closed and bounded subsets in $\Ea$ are compact. 

We also recall that a subspace $E\subset \Ea$ is Polish if and only if it  is the intersection of a countable collection of open subsets. 
The distance from $x\in \Ea$ to a measurable subset $\Aa\subset\Ea$ is denoted by
$$
d(x,\Aa):=\inf{\left\{d(x,y)~:~y\in \Aa\right\}}.
$$ 
We also denote by 
$\partial E:=\overline{E}-E^{\circ}$ the boundary of some domain (open and connected) $E\subset \Ea$, where $\overline{E}$ and $E^{\circ}$ stand for the closure and the interior of a subset $E$. 

In the further development of the article, $\cchi$ stands for some
 decreasing positive function $\cchi$ on $]0,\infty[$ such that for any $0<\alpha<1$
 we have
 $$
\lim_{\alpha\rightarrow 0}  \cchi(\alpha)=+\infty\qquad
 \cchi(\alpha)< 1/\alpha\quad \mbox{\rm and}\quad
 \overline{\cchi}(\alpha):=\int_0^{\alpha}\cchi(u)du<\infty.
$$
\begin{defi}
We associate with $\cchi$  the function $V_{\partial}\in \Ca(E)$ defined by
\begin{equation}\label{def-V-delta}
V_{\tiny \partial}~:~x\in E~\mapsto~V_{\tiny \partial}(x):=\cchi(d(x,\partial E))\in ]0,\infty[.
\end{equation}
\end{defi}
For instance, we can choose
$
\cchi(u)=1/u^{1-\epsilon}
$, for some $\epsilon\in ]0,1[$.  For any $r>0$ the $r$-sub-level sets of $V_{\partial}$ are given by the closed subsets
$$
\Va_{\tiny \partial}(r):=\{x\in E~:~V_{\tiny \partial}(x)\leq r\}=\{x\in E~:~d(x,\partial E)\geq \cchi^{-1}(r)\}.
$$
Note that $V_{\tiny \partial}\in \Ca_{\infty}(E)$ as soon as
$ \overline{E}$ is compact. 
\subsection{Bounded domains}\label{bounded-dom-sec}
  Let $E\subset \Ea:=\RR^n$ be some  bounded  domain with locally Lipschitz boundary $\partial E$, for some $n\geq 1$. Consider a semigroup of integral operators 
\begin{equation}\label{ref-Q-absolute-cont}
Q_{t}(x,dy)=q_{t}(x,y)~dy
\end{equation}
 having for any $t>0$ a bounded density $(x,y)\in E^2\mapsto q_{t}(x,y)\in [0,\infty[$ w.r.t. the trace of the Lebesgue measure $\nu(dy)=dy$ on $E$.   In this situation, we have the following lemma.
 \begin{lem}\label{lem-lip-dom-compact}
 For any $t>0$ we have
 \begin{equation}\label{lip-dom-compact}
 V_{\partial}\in \Ca_{\infty}(E)\quad \mbox{and}\quad
 \Vert Q_t(V_{\partial})\Vert \leq c_t~\int_E \cchi(d(x,\partial E))~dx<\infty.
 \end{equation}
 \end{lem}
The proof of the above lemma follows from an elementary change of variable formulae, thus it is provided in the appendix, on page~\pageref{lip-dom-compact-proof}.

The estimate (\ref{lip-dom-compact}) clearly applies to the class of sub-Markov semigroups $Q^{[U]}_t$  defined in (\ref{Q-U-gen}) for any choice of the absorption 
potential function, as soon as the semigroup $Q^{[U]}\ll Q$ is dominated by a collection of integral operators $Q_t(x,dy)$
having a bounded density $q_t(x,y)$ on $E^2$ w.r.t. the Lebesgue measure on $E$. For instance, when the transition semigroup of the free evolution flow $X_t(x)$ in (\ref{Q-U-gen}) has a density $p_t(x,y)$ for any non negative function $f$ on $E$ and any $x\in E$ we have
$$
Q^{[U]}_t(f)(x)\leq \int~q_t(x,y)~f(y)~dy\quad \mbox{\rm with}\quad q_t(x,y):=p_t(x,y)~1_{E}(y).
$$

We summarize the above discussion with the following proposition.
\begin{prop}\label{prop-domain-bounded}
Assume that $Q^{[U]}\ll Q$ is dominated by a collection of integral operators $Q_t$ satisfying (\ref{ref-Q-absolute-cont}). Then, 
$$
Q^{[U]}_t(V_{\partial})/V_{\partial}\leq c_t/V_{\partial}\in \Ba_0(E).
$$
\end{prop}

The choice of the Lyapunov function $V$ is clearly not unique. For instance, when $E=]0,1[$ instead of $V_{\partial}$ we can choose $V(x):=1/\sqrt{x}+1/\sqrt{1-x}$. For the Dirichlet heat kernel discussed in Section~\ref{dhk-sec} we can also choose $V(x)=1/\sin{(\pi x)}$.

We emphasize that sub-Markov  integral operators on the compact interval $E=[0,1]$
 with {\em a positive} continuous density w.r.t. the Lebesgue measure on $E$
  arise when the free evolution process is reflected at both sides of the interval. 
  In this context  the process is not conditioned by  any type of non absorption  at the boundaries. In this context, the unit function $V=1$ belongs to 
  $ \Ba_{\infty}(E)$. In the same vein, sub-Markov  integral operators with mixed boundary conditions on the left-closed interval $E=[0,1[$, or respectively on the right-closed interval $E=]0,1]$ arise when the free evolution process is reflected at the Neumann boundary $ \partial_{N}E:=\{0\}$ and non absorbed at the Dirichlet boundary $ \partial_{D}E=\{1\}$, or respectively 
  reflected at $ \partial_{N}E:=\{1\}$ and non absorbed at $ \partial_{D}E=\{0\}$. 

  More generally, 
consider a bounded domain $\Omega\subset \RR^n$ with Lipschitz boundary $\partial \Omega=\partial_{D}\Omega\cup \partial_{N}\Omega$ consisting of two disjoint 
  connected components $ \partial_{D}\Omega$ and $ \partial_{N}\Omega$ closed in $\RR^n$,   and set $E:=\Omega\cup \partial_{N}\Omega$. 
In this notation,   the  function $V_{\partial}(x):=\chi\left(d(x,\partial_D E)\right)$ belongs to  $ \Ca_{\infty}(E)$. 
   In addition,  for any bounded density $q_t(x,y)$ on $E^2$ we have the uniform estimate
$$
  \int_E~q_t(x,y)~V_{\partial}(y)~dy \leq c_t~\int_E V_{\partial}(y)~dy<\infty.$$
  The above estimate also holds for the function $V_{\partial}(x)=\chi\left(d(x,\partial E)\right)$. 
  \subsection{Unbounded domains}\label{bo-sec}

  When the domain $E$ is not bounded the function $V_{\tiny \partial}\not\in \Ba_{\infty}(E)$. In this context, one natural way to design a Lyapunov function $V\in \Ba_{\infty}(E)$ is to consider an auxiliary function $V_{\Ea}\in \Ca_{\infty}(\Ea)$ with $V_{\Ea}(x)\geq 1$ for any $x\in E$. In this situation, we have 
  $$V:=V_{\tiny \partial}+V_{\Ea}\in\Ca_{\infty}(E).$$
  To check this claim, observe that the sub-level sets of $V_{\tiny \partial}$ are given by the closed subsets
  $$
 \Va_{\tiny \partial}(r):=  \{V_{\tiny \partial}\leq r\}=\{x\in E~:~d(x,\partial E)\geq \cchi^{-1}(r)\} \subset E 
  $$
and we have the compact inclusion
$$
\Va(r):= \{V\leq r\}\subset \Va_{\tiny \Ea}(r)\cap \Va_{\tiny \partial}(r)\quad\mbox{\rm with}\quad \Va_{\tiny \Ea}(r):=\{x\in \Ea~:~V_{\Ea}(x)\leq r\}.
$$
This yields the following easily checked proposition.
\begin{prop}\label{prop-bb-ub}
For any $t>0$ we  have
$$
\Vert Q_t(V_{\tiny \partial})\Vert\vee
\Vert Q_t(V_{\Ea})\Vert<\infty\quad \Longrightarrow \quad Q_t(V)/V\leq c_t/V\in \Ba_0(E).
$$
When $Q_t(1)\in \Ba_0(E)$ we also have
$$
\begin{array}{l}
\displaystyle 
\Vert \overline{Q}_t(V_{\tiny \partial})\Vert\vee\Vert\overline{Q}_t(V_{\Ea})/V_{\Ea}\Vert<\infty
\Longrightarrow
Q_t(V)/V
\leq c_t~Q_t(1)\in \Ba_0(E).
\end{array}$$
\end{prop}

The design of a function $V_{\Ea}$  is rather flexible. For instance, assume that  $\overline{Q}\ll P$ is dominated by some Markov integral operators $P_t$ on $\Ba_b(\Ea)$ such that $\Vert P_t(V_{\Ea})/V_{\Ea}\Vert<\infty$ for some $V_{\Ea}\in \Ba_{\infty}(\Ea)$.
In this situation, we have $\Vert\overline{Q}_t(V_{\Ea})/V_{\Ea}\Vert<\infty
$ as well as
$$
 \Vert V_{\Ea}~Q_t(1)\Vert<\infty\Longrightarrow
\Vert Q_t(V_{\Ea})\Vert<\infty.
$$
For instance, when $P_t$ satisfies the sub-Gaussian estimates (\ref{def-OU-A-Sigma-lem}) on $\Ea=\RR^n$ we can choose $V_{\Ea}(x):=1+\Vert x\Vert^k$, for some $k\geq 1$, as soon as the function $Q_t(1)(x)\longrightarrow_{\Vert x\Vert\rightarrow\infty}0$ faster than $\Vert x\Vert^{-k}$.

When the domain $E$ and its boundary $\partial E$ are both non necessarily bounded, 
it may happens that $Q_t(1)\in \Ba_0(E)$ but $\overline{Q}_t(V_{\tiny \partial})\not\in \Ba_b(E)$. In this situation, we can use the following proposition.

\begin{prop}
Assume there exists some $V_{\Ea}\in \Ca_{\infty}(\Ea)$ with $V_{\Ea}(x)\geq 1$ for any $x\in E$ and such  that
$$
\Vert\overline{Q}_t(V_{\tiny \partial})/V_{\Ea}\Vert\vee\Vert\overline{Q}_t(V_{\Ea})/V_{\Ea}\Vert\vee \Vert Q_t(1)V_{\Ea}\Vert<\infty.$$
Then we have
$$
Q_t(V)/V\leq c_t/V\in \Ba_0(E).
$$
\end{prop}
\proof
Using the following  decompositions
$$
Q_t(V_{\tiny \partial})=Q_t(1)V_{\Ea}~\overline{Q}_t(V_{\tiny \partial})/V_{\Ea}\quad
\mbox{\rm and}\quad Q_t(V_{\Ea})=Q_t(1)V_{\Ea}~\overline{Q}_t(V_{\Ea})/V_{\Ea}
$$
and applying Proposition~\ref{prop-bb-ub} we  have
$$
\Vert Q_t(V_{\tiny \partial})\Vert\vee
\Vert Q_t(V_{\Ea})\Vert<\infty \quad\mbox{\rm and therefore}\quad \quad Q_t(V)/V\leq c_t/V\in \Ba_0(E).
$$
This ends the proof of the proposition.
\cqfd

The case $Q_t(1)\not\in \Ba_0(E)$ can also be handle whenever the pair $(V_{\partial},V_{\Ea})$ can be chosen so that 
\begin{equation}\label{V-cond}
\forall \delta >0\qquad
V_{\partial}\,V_{\Ea}^{\delta}\in \Ca_{\infty}(E).
\end{equation}
For instance we can choose for some $v>0$ and $\epsilon\in ]0,1[$ the functions
$$
 V_{\Ea}(x):=\exp{(v\Vert x\Vert)}\quad \mbox{\rm and}\quad \cchi(u):=1/u^{1-\epsilon}.
$$
Observe that
$$
d(x,\partial E)\leq \Vert x\Vert+d(0,\partial E)
\quad \mbox{\rm and}\quad
V_{\partial}(x)\geq \cchi(\Vert x\Vert+(1\vee d(0,\partial E)))
$$
and for any $m\geq 0$ and $\delta>0$ we have
$$
V_{\Ea}(x)\geq c_v(m,\delta)~(1+\Vert x\Vert)^{(m+1)/\delta}.
$$
This implies that
$$
V_{\Ea}(x)^{\delta}V_{\partial}(x)\geq c_2\frac{(1+\Vert x\Vert)^{m+1}}{(\Vert x\Vert+(1\vee d(0,\partial E)))^{1-\epsilon}}\geq c~(1+\Vert x\Vert)^{m+\epsilon}.
$$
Using the fact that $V_{\Ea}(x)\geq 1$ for any $x\in E$, this implies that
$$
\{x\in E~:~V_{\Ea}(x)^{\delta}V_{\partial}(x)\leq r\}\subset \{ x\in E~:~c~(1+\Vert x\Vert)^{m+\epsilon}\leq r\}\cap\{x\in E~:~V_{\partial}(x)\leq r\}.
$$
We conclude that $V_{\Ea}^{\delta}V_{\partial}$ has compact level sets and (\ref{V-cond}) is satisfied.

In this context, we have the following proposition.
\begin{prop}\label{cond-gen-unb-prop}
Consider a couple of functions $(V_{\partial},V_{\Ea})$ satisfying (\ref{V-cond}). Assume  there exists some parameters $t>0$,  $\delta_t>0$ and $\epsilon\geq 0$  such that
\begin{equation}\label{cond-gen-unb}
Q_t(V_{\Ea})/V_{\Ea}\leq c_t/V_{\Ea}^{\delta_t} \quad \mbox{and}\quad 
Q_t(V_{\partial})\leq c_t~V_{\Ea}^{\,\epsilon\, \delta_t}.
\end{equation}
In this situation, for any  $p>1+\epsilon$ we have 
 $$
V:=V_{\Ea}^{1-1/p}~V_{\partial}^{1/p}\in\Ca_{\infty}(E)
$$
as well as
$$
Q_t\left(V\right)/V\leq c_t/(V_{\Ea}^{\delta_t \epsilon_p}V_{\partial}^{1/p})\in \Ca_0(E) \quad \mbox{with}\quad
\epsilon_p:=(1-(1+\epsilon)/p)>0.
$$
\end{prop} 
\proof
Observe that
for any $p>1+\epsilon$ we have
$$
V_{\partial} V_{\Ea}^{p-1}\in\Ca_{\infty}(E)\quad \mbox{\rm and therefore}\quad
V:=V_{\partial}^{1/p}\,V_{\Ea}^{1-1/p}\in\Ca_{\infty}(E).
$$
In the same vein, for any $\epsilon\geq 0$ we have
$$
(\ref{V-cond})\Longrightarrow
V_{\partial} V_{\Ea}^{p\delta_t \epsilon_p}\in\Ca_{\infty}(E)\quad \mbox{\rm and therefore}\quad V_{\partial}^{1/p}\,V_{\Ea}^{\delta_t \epsilon_p}\in \Ca_{\infty}(E).
$$
On the other hand, using H\" older's inequality, we have
\begin{eqnarray*}
Q_t\left(V\right)/V&\leq& (Q_t(V_{\Ea})/V_{\Ea})^{1-1/p}~\left(Q_t(V_{\partial})/V_{\partial}\right)^{1/p}\\
&\leq &c_t(1) (Q_t(V_{\partial})/(V_{\Ea}^{\delta_t (p-1)}V_{\partial}))^{1/p}\leq
c_t(2) (1/(V_{\Ea}^{p\delta_t (1-(1+\epsilon)/p)}V_{\partial}))^{1/p}.
\end{eqnarray*}
This ends the proof of the proposition.
\cqfd

The design of a function $V_{\Ea}$ satisfying (\ref{cond-gen-unb}) is rather flexible. For instance, (\ref{cond-gen-unb}) is automatically satisfied when $Q\ll P$ is dominated by some Markov integral operators $P_t$ on $\Ba_b(\Ea)$ such that
$$
P_t(V_{\Ea})(x)/V_{\Ea}(x)\leq c_t(1)/V_{\Ea}(x)^{\delta_t}.
$$
Section~\ref{sec-Markov-lyap} discusses a variety of Lyapunov functions $V_{\Ea}$ satisfying the above condition for Markov diffusion semigroups. These Lyapunov functions can also 
be designed using the domination principles presented in Section~\ref{sg-domination}. For instance, consider the semigroup $Q_t:=\Qa_t^{(a)}$ associated with the  Langevin diffusion flow  on a cylinder discussed in (\ref{langevin-ref-intro}). In this situation, combining (\ref{ref-dom-il}) with Proposition~\ref{ref-lem-OU-A-Sigma}, for any $v\geq 0$ and $t>0$ there exists some finite constant $\delta_t>0$ such that
$$
V_{\Ea}(x):=\exp{\left(v\vert x\vert\right)}\Longrightarrow
Q_{t}(V_{\Ea})/V_{\Ea}\leq c_{t}/V_{\Ea}^{\delta_{t}}.$$

Next, we illustrate the r.h.s. condition in (\ref{cond-gen-unb}) when $q_t$ are sub-Gaussian densities; in the sense 
that for any $x,y\in E$ we have
\begin{equation}\label{ref-sub-gauss-b}
q_{t}(x,y)\leq c_t~g_t(x,y)\quad \mbox{\rm with} \quad g_t(x,y):=\frac{1}{(2\pi \sigma_t^2)^{n/2}}~\exp{\left(-\frac{1}{2\sigma_t^2}~\Vert y-m_t(x)\Vert^2\right)}
\end{equation}
for some parameter $\sigma_t>0$ and some non necessarily bounded 
function $m_t$ on $E$.

\begin{prop}
Let $\varphi$ be a Lipschitz function on $\RR^{n-1}$ with uniformly bounded gradient and set
$$
E:=\{x=(x_i)_{1\leq i\leq n}\in \RR^{n}~:~x_n>\varphi(x_{-n})\}
\quad \mbox{with}\quad
 x_{-n}:=(x_i)_{1\leq i<n}\in \RR^{n-1}.
$$
Then  the r.h.s. condition in (\ref{cond-gen-unb}) is met with $\epsilon=0$ for any positive semigroup satisfying (\ref{ref-sub-gauss-b}). The same property holds when the boundary $\partial E$ can be decomposed as a finite union of graphs of   differentiable functions on $\RR^{n-1}$ with uniformly bounded gradients. 
\end{prop}
\proof
We choose $\alpha>0$ sufficiently small so that for any $$x\in  
D_{\alpha}(E):=\{x\in \overline{E}~:~d(x,\partial E)\leq  \alpha\}$$
there exists a  projection $ \overline{x}\in \partial E$ with
$
d(x,\partial E)=\Vert x- \overline{x}\Vert
$. 
Let $\Ca_{\varpi}(\overline{x})$ be an interior cone with a given base vertex  $\overline{x}=(x_{-n},\varphi(x_{-n}))\in \partial E$ and a given  half-opening angle $\varpi$ around the axis $\Aa(\overline{x}):=\{(x_{-n},x_n)~:~x_n\geq \varphi(x_{-n})\}$. For any $x\in \Aa(\overline{x})$  there exists a  projection $ \widehat{x}\in \partial \Ca_{\varpi}(\overline{x})$ on the boundary $\partial \Ca_{\varpi} (\overline{x})$ with $$
 d(x,\partial \Ca_{\varpi} (\overline{x}))=d(x,\widehat{x})=\cos\left(\frac{\pi}{2}-\varpi\right)~(x_n-\varphi(x_{-n}))\leq 
  d(x, \overline{x})
$$ 
On the other hand, for any $y\in \partial E$ we have
 \begin{eqnarray*}
   z:=(y_{-n},\varphi(x_{-n}))&\Longrightarrow& 0\leq \frac{\pi}{2}-\varpi\leq \widehat{y \overline{x} z}
 \quad\mbox{\rm 
 and}\quad
 \tan(\widehat{y \overline{x}  z})=
 \frac{\vert \varphi(x_{-n})-\varphi(y_{-n})\vert }{\Vert x_{-n}-y_{-n}\Vert}. \end{eqnarray*}
 This yields the estimate
 $$
 \cos\left(\frac{\pi}{2}-\varpi\right)\geq \cos\left(\widehat{y \overline{x} z}\right)=\frac{1}{\sqrt{1+\tan^2(\widehat{y \overline{x}  z})}}
 $$
 from which we conclude that
 $$
 0\leq x_n-\varphi(x_{-n})\leq   \kappa~ d(x, \overline{x})
$$
with
$$
  \kappa:=\sqrt{1+\Vert \nabla \varphi \Vert^2}
  \quad \mbox{\rm and}\quad
  \Vert \nabla \varphi\Vert:=\sup_{y\in \RR^{n-1}}\Vert \nabla \varphi(y)\Vert<\infty.
 $$
  
 This implies that
 $$
 \begin{array}{l}
\displaystyle \int_{D_{\alpha}(E)} \chi\left(d(y,\partial E)\right)
q_t(x,y)~dy\\
 \\
\displaystyle \leq  c_t(1)~\int_{D_{\alpha}(E)} \chi\left(\left(y_n-\varphi(y_{-n})\right)/\kappa)\right)~\exp{\left(-\frac{1}{2\sigma^2_t}( (y_n-\varphi(y_{-n}))+(\varphi(y_{-n})-(m_t(x))_n))^2\right)} \\
 \\
\displaystyle \hskip3cm~\times~\exp{\left(-\frac{1}{2\sigma_t^2}~\Vert y_{-n}-(m_t(x))_{-n}\Vert^2\right)}~ dy_ndy_{-n}.
 \end{array}
 $$
Using the change of variables
 $$
 z:=\left(y_n-\varphi(y_{-n})\right)/\kappa\Longrightarrow
d y_n= \kappa~dz
 $$
we find that
 $$
 \begin{array}{l}
\displaystyle \int_{D_{\alpha}(E)} \chi\left(d(y,\partial E)\right)
q_t(x,y)~dy\\
 \\
\displaystyle \leq  \kappa~ c_t(1)~\overline{\cchi}(\alpha)~\int_{\RR^{n-1}} ~\exp{\left(-\frac{1}{2\sigma_t^2}~\Vert y_{-n}-(m_t(x))_{-n}\Vert^2\right)}~ dy_{-n}\leq c_t(2).
 \end{array}
 $$
 On the other hand, for any $\alpha>0$ we have
 $$
 \int_{E-D_{\alpha}(E)} \chi\left(d(y,\partial E)\right)~q_t(x,y)~dy\leq 
  \chi\left(\alpha\right)~\Vert Q_t(1)\Vert \leq c_t(3).
 $$
 This ends the proof of the proposition.
 \cqfd
\subsection{Smooth boundaries}\label{sec-smooth-boundary}
Next, we illustrate the Lyapunov conditions on $V_{\partial}$  in the context of absolutely continuous sub-Markov semigroup of the form (\ref{ref-Q-absolute-cont}) with a bounded density $q_t(x,y)$ on a non necessarily bounded domain $E\subset \RR^n$ with  smooth non necessarily bounded $\Ca^2$-boundary with uniformly bounded interior curvature. 

We assume that there exists  $\alpha>0$ sufficiently small so that every point of the $\alpha$-offset of $\partial E$ (a.k.a. $\alpha$-tubular neighborhood) defined by
$$
\mbox{\em Tub}_{\alpha}(\partial E):=\{x\in \RR^n~:~d(x,\partial E)\leq \alpha\}
$$
 lies on some normal ray passing through a point on $\partial E$ and no two normal rays passing through different points of $\partial E$ intersect in $\mbox{\em Tub}_{\alpha}(\partial E)$. We let $N(z)$ be the unit normal vector at $z\in \mbox{\em Tub}_{\alpha}(\partial E)$ pointing inward $E$, and
let  $\Da_{r}(E)$ the closed subset defined for any $r\leq \alpha$ by
 $$
\Da_{r}(E):=\{x\in \overline{E}~:~d(x,\partial E)\leq  r\}\quad
\mbox{\rm and}\quad \Da_{-r}(E):=\{x\in \RR^n-E~:~d(x,\partial E)\leq  r\}.
$$ 
In this notation, the inverse of the normal coordinate map
\begin{equation}\label{Fermi-coordinates}
F~:~(z,r)\in \partial E\times [-\alpha,\alpha]\mapsto F(z,r)=z+r\,N(z)\in \mbox{\em Tub}_{\alpha}(\partial E)
\end{equation}
is given for any $x\in  \mbox{\em Tub}_{\alpha}(\partial E)
$ by
$$
F^{-1}(x)=\left(\mbox{\rm proj}_{\partial E}(x),d_{\alpha}(x,\partial E)\right)
$$
where $\mbox{\rm proj}_{\partial E}(x)$ stands for the projection of $x\in  \mbox{\em Tub}_{\alpha}(\partial E)
$ onto $\partial E$ and $d_{\alpha}(x,\partial E)$ stands for the signed distance function
$$
d_{\alpha}(x,\partial E):=d(x,\partial E)~1_{\Da_{\alpha}(E)}(x)-
d(x,\partial E)~1_{\Da_{-\alpha}(E)}(x)\in [-\alpha,\alpha].
$$
In addition, the inward normal $N(x)$ at any $x$ on the  $\Ca^2$ boundary $\partial E$ is given by
$$
\nabla d_{\alpha}(x,\partial E)=N(x).
$$
The Hessian of the signed distance function on the boundary $\partial E$ gives the Weingarten map $\Wa(x)$. With this notation at hand, we have
$$
\int_{\Da_{\alpha}(E)}~f(d(y,\partial E))~q_{t}(x,y)~dy
=\int_0^\alpha f(r)~q^{\partial}_t(x,r)~dr
$$
with the level-set density function
\begin{eqnarray}
q^{\partial}_t(x,r)&:=&\int_{\partial E_r}~q_{t}(x,y)~\sigma_{\partial,r}(dy)\label{level-set-density}\\
&=&\int_{\partial E}~q_t\left(x,z+rN(z)\right)~\left\vert\mbox{\rm det}\left(I-r~\Wa(z)\right)\right\vert~\sigma_{\partial}(dz).\nonumber
\end{eqnarray}
In the above display, $\sigma_{\partial,r}(dz)$ stands for the Riemannian volume measure on the $r$-extended boundary 
$$
\partial E_r:=\{x\in E~:~d(x,\partial E)=r\}.$$

Moreover, 
since $E$ has uniformly bounded interior curvature, for any $r\leq \alpha$ we have
\begin{gather*}
\kappa_{\partial}(\alpha):=\sup\left\vert\mbox{\rm det}\left(I- r~\Wa(z)\right)\right\vert<\infty
\quad \mbox{\rm and}\quad \kappa^-_{\partial}(\alpha):=\sup\left\vert\mbox{\rm det}\left(I+r~\Wa(y)\right)\right\vert<\infty.
\end{gather*}
In the above display, the supremum is taken over all
 $z\in \partial E$,   $y\in \partial E_r$, and  $r\leq \alpha$. 
 Several examples of hypersurface boundaries satisfying the above conditions are discussed in Section~\ref{hyper-surface-bound-sec}  (cf. for instance Proposition~\ref{prop-vol-form-weignarten}).

We denote by $\overline{q}^{\partial}_t\geq q^{\partial}_t$ the function defined as $q^{\partial}_t$ by replacing $q_t$ by $\overline{q}_t$. 
Using the fact that
$$
\overline{Q}_t(V_{\tiny \partial})(x)\leq \cchi(\alpha)+
\int_0^\alpha~\cchi(r)~ \overline{q}^{\partial}_t(x,r)~dr
$$
we readily check the following proposition.
\begin{prop}
For any
 $t>0$ we have
\begin{eqnarray}
\sup_{0\leq r\leq \alpha}\sup_{x\in E} q^{\partial}_t(x,r)<\infty&\Longrightarrow&
Q_t(V_{\tiny \partial})\leq \cchi(\alpha)~Q_t(1)+  c_t(\alpha) ~\overline{\cchi}(\alpha)\nonumber\\
\sup_{0\leq r\leq \alpha}\sup_{x\in E} \overline{q}^{\partial}_t(x,r)<\infty&\Longrightarrow& \overline{Q}_t(V_{\tiny \partial})\leq \cchi(\alpha)+ c_t(\alpha)~ \overline{\cchi}(\alpha).\label{ub-obst-ref}
\end{eqnarray}
When the boundary $\partial E$ is bounded, for any $t>0$ we have the estimate
\begin{equation}\label{ref-bounded-boundary}
 \Vert\overline{Q}_t(V_{\tiny \partial})\Vert\leq  c_t(\alpha)~\left(\cchi(\alpha)+ \overline{\cchi}(\alpha)~\sup_{0\leq r\leq \alpha}~\sigma_{\partial,r}\left(\partial E_r\right)\right).
\end{equation}

\end{prop}

We end this section with some practical tools to estimate the level-set density functions discussed in Section~\ref{sec-smooth-boundary}. Most of our estimates are based on the following technical lemma.

\begin{lem}\label{lem-boundary-f-g}
Consider a couple of non negative functions $f,g$ on $\RR^n$ and some parameter $\alpha>0$ such that
$$
 \sup_{\Vert u\Vert\leq \alpha}f(z+u)\leq \iota(\alpha)~g(z)
\quad\mbox{
for some $\iota(\alpha)<\infty$.}
$$
 In this situation, we have the uniform estimate
$$
\sup_{0\leq r\leq \alpha}\int_{\partial E_{r}}~f(z)~\sigma_{\partial,r}(dz)\leq \iota(\alpha)~\kappa_{\partial}(\alpha)~\int_{\partial E}~g\left(z\right)~\sigma_{\partial}(dz)
$$
as well as the co-area estimate
$$
\int_{\partial E}~f(z)~\sigma_{\partial}(dz)\leq \frac{1}{\alpha}~\iota(\alpha)~\kappa_{\partial}^-(\alpha)\int_{\Da_{\alpha}(E)}~g(z)~dz.
$$
\end{lem}
The proof of the above lemma is provided in the appendix, on page~\pageref{lem-boundary-f-g-proof}.

Note that the level-set density function defined in (\ref{level-set-density}) can be estimated for any $0\leq r\leq \alpha$ by the formula
\begin{eqnarray}
q^{\partial}_t(x,r)&\leq& \kappa_{\partial}(\alpha)~\int_{\partial E}~q_t\left(x,z+rN(z)\right)~~\sigma_{\partial}(dz).\nonumber
\end{eqnarray}
\begin{prop}
Assume that $q_t(x,y)\leq \varpi_t~g_t(x,y)$ is dominated by some probability density $y\mapsto g_t(x,y)$ on $\RR^n$ for some $t>0$ and  some parameter $\varpi_t<\infty$. In addition, we have
\begin{equation}\label{ref-cond-q-g}
 \sup_{\Vert u\Vert\leq \alpha}g_t(x,y+u)\leq \iota_t(\alpha)~g_{\alpha,t}(x,y)
\end{equation}
for some probability density  $y\mapsto g_{\alpha,t}(x,y)$ and some $\iota_t(\alpha)<\infty$. 
In this situation, we have the uniform density estimates
\begin{equation}\label{ref-cond-q-g-est}
\sup_{0\leq r\leq \alpha}\sup_{x\in E} q^{\partial}_t(x,r)\leq \varpi_t~\iota_t(\alpha)\kappa_{\partial}^-(\alpha)\kappa_{\partial}(\alpha)/\alpha.
\end{equation}
\end{prop}
\proof
By (\ref{ref-cond-q-g}) for any $0\leq r\leq \alpha$ we have
$$
q^{\partial}_t(x,r)
\leq \varpi_t~ ~\kappa_{\partial}(\alpha)~\int_{\partial E}~g_t\left(x,z+rN(z) \right)~~\sigma_{\partial}(dz).
$$
On the other hand, we have 
$$\int_{\Da_{\alpha}(E)}~g_{t,\alpha}(x,y)~dy\leq 1.$$ 
The estimate (\ref{ref-cond-q-g-est}) is now a direct consequence of the co-area estimate stated in Lemma~\ref{lem-boundary-f-g}.
This ends the proof of the proposition. 
\cqfd
We illustrate the above condition when $q_t$ are the sub-Gaussian densities discussed in (\ref{ref-sub-gauss-b}). In this situation, using the fact that $2a^{\prime}b\leq \frac{1}{\epsilon}\Vert a\Vert^2+\epsilon \Vert b\Vert^2 $ for any $0<\epsilon<1$ and $\Vert u\Vert\leq \alpha $  we check that
$$
-\frac{1}{2\sigma_t^2}~\Vert (y+u)-m_t(x)\Vert^2\leq -\frac{(1-\epsilon)}{2\sigma_t^2}~\Vert y-m_t(x)\Vert^2+\frac{1}{2\sigma_t^2}\left(\frac{1}{\epsilon}-1\right)~\alpha^2.
$$
In this context, condition (\ref{ref-cond-q-g}) is met with the gaussian density
$$
g_{\alpha,t}(x,y):=\frac{e^{-\frac{1}{2\sigma_t(\epsilon)^2}~\Vert y-m_t(x)\Vert^2}}{(2\pi \sigma_t(\epsilon)^2)^{n/2}}~$$
with
$$
\iota_t(\alpha):=c_t~e^{\frac{\alpha^2/\epsilon}{2\sigma_t(\epsilon)^2}}
\quad \mbox{\rm and}\quad\sigma_{t}(\epsilon)^2:=\sigma_t^2/(1-\epsilon).
$$

\section{Riccati type processes}\label{ricc-sec}

\subsection{Positive diffusions}
Consider the Riccati type diffusion on $E=]0,+\infty[$ defined  for any $x\in E$ by
$$
dX_t(x)=\left(a_0+a_1~X_t(x)-b~X_t(x)^2\right)~dt+ \sigma_1(X_t(x))~dB^1_t+
\sigma_2(X_t(x))~dB^2_t,  \ \  X_0(x) = x
$$
for some Brownian motion $(B^1_t,B^2_t)$ on $\RR^2$, the diffusion functions
$$
\sigma_1(x):=\varsigma_1~\sqrt{x}\qquad
\sigma_2(x):=\varsigma_2~x
$$
and the parameters 
 $$a_1\in\RR\qquad a_0>\varsigma_1^2\qquad b>0\quad \mbox{\rm and}\quad
 \varsigma_1, \varsigma_2\geq 0.
$$
Applying It\^o's formula, we readily check that
$$
\partial_t\EE(X_t(x))\leq \mbox{\rm Ricc}\left(\EE(X_t(x))\right)\quad\mbox{\rm and}\quad
\partial_t\EE(1/X_t(x))\leq \mbox{\rm Ricc}^-\left(\EE(1/X_t(x))\right)
$$
with the Riccati drift functions defined by
\begin{equation}\label{ref-riccati-drift}
\mbox{\rm Ricc}(z):=a_0+a_1z-b z^2
\quad \mbox{and}\quad
\mbox{\rm Ricc}^-(z):=a_0^-+a_1^-z-b^- z^2
\end{equation}
with the parameters
$$
a_0^-:=b\quad a_1^-:=(\varsigma_2^2-a_1) \quad \mbox{\rm and}\quad b^-:=a_0-\varsigma_1^2.
$$
Consider the Lyapunov function $V\in \Ba_{\infty}(E)$ defined by $V(x):=x+1/x$.
By well known properties of Riccati flows, for any $t>0$  we have
$
\Vert P_{t}(V)\Vert<\infty 
$. 
 For a more thorough discussion on this class of one-dimensional Riccati diffusions, we refer to the article~\cite{bishop-19}.

\subsection{Matrix valued diffusions}
Let $\overline{E}$ and $E$ be the space of  $(n\times n)$-positive semi-definite and definite matrices respectively. Also let  $\lambda_1(x)\geq \ldots\geq \lambda_n(x)$ denote the ordered eigenvalues of $x\in \overline{E}$.
Let $\Wa_t$ denotes an $(n\times n)$-matrix with independent Brownian entries. Also let $A$ be an $(n\times n)$-matrix with real entries and let 
$R,S\in E$. We associate with these objects the $E$-valued diffusion
$$
	dX_t \,=\, (AX_t + X_tA^{\prime} + R - X_tSX_t)\,dt \,+\, \frac{\epsilon}{2}\left[X_t^{1/2}d\Wa_t\,R^{1/2} + R^{1/2}\,d\Wa_t^{\prime}\,X_t^{1/2}\right].
$$
Whenever $\epsilon\leq  2/\sqrt{n+1}$, the diffusion $X_t$ has a unique strong solution that never hits the boundary $\partial E=\overline{E}-E$. In addition, 
the transition semigroup $P_t$ of $X_t$ is strongly Feller and admits a smooth density w.r.t. the Lebesgue measure on $E$, thus it is irreducible.
Furthermore, when 
$\epsilon^2(1+n)/2\leq \lambda_{n}(R)/\lambda_{1}(R)$ then the function $V(x)=\tr( x)+\tr (x^{-1})$ is a Lyapunov function with compact level subsets. For a detailed proof of the above assertion for more general classes of Riccati matrix valued diffusions we refer to~\cite{bishop-19-2} (see for instance the stability Theorem 2.4 and Section 5.4 in~\cite{bishop-19-2}).

\subsection{Logistic birth and death process}
Let $X_t(x)$ be the stochastic flow on $E:=\NN-\{0\}$ with generator
$L$ defined for any $f\in \Ba_b(E)$ and $x\geq 2$ by 
$$
L(f)(x)=J(x,x-1)~(f(x-1)-f(x))+J(x,x+1)~(f(x+1)-f(x)) 
$$
and for $x=1$ by
$$
L(f)(1)=J(1,2)(f(2)-f(1)).
$$
In the above display, the birth and death rates  are given by
\begin{equation}\label{ex-2-bd}
J(x,x+1):=\lambda_b~x+\upsilon_b\quad \mbox{\rm and}\quad J(x,x-1):=\lambda_d~x+ \lambda_l~ x(x-1)+\upsilon_d
\end{equation}
for some non negative parameters  $\lambda_d,\lambda_b,\upsilon_b,\upsilon_d\geq 0$ and $
\lambda_l>0
$. Consider the identity function  $V:x\in E\mapsto V(x)=x$.
For any $x\geq 2$ we have
$$
L(V)(x)=
J(x,x+1)-J(x,x-1)= \left(F\circ V\right)(x)
$$
with the concave  function
\begin{equation}\label{drift-riccati}
z\in \RR_+\mapsto F(z):=(\upsilon_b-\upsilon_d)+(\lambda_b+\lambda_l-\lambda_d)~z-\lambda_l~z^2\in \RR.
 \end{equation}
 Observe that
 $$
 L(V)(1)-F(V(1))=J(1,2)-F(1)=J(1,0)=\upsilon_d+\lambda_d.
 $$
This yields the estimate
\begin{eqnarray*}
P_t(L(V))(x)&=&P_t(1_{[2,\infty[}~L(V))(x)+P_t(1_{\{1\}}~L(V))(x)\\
&=&P_t((F\circ V))(x)+P_t(1_{\{1\}})(x)~\left(L(V)(1)-F(V(1))\right)\\
&\leq &F(P_t(V))(x) +J(1,0)
\end{eqnarray*}
from which we check that
\begin{eqnarray*}
\partial_t P_t(V)(x)&\leq 
&\mbox{\rm Ricc}\left(P_t(V)(x)\right)\end{eqnarray*}
with the Riccati drift function defined in (\ref{ref-riccati-drift}) 
with the parameters
$$
a_0:=\upsilon_b+\lambda_d\qquad a_1:=\lambda_b+\lambda_l-\lambda_d \quad \mbox{\rm and}\quad
b:=\lambda_l>0.
$$
By well known properties of Riccati flows, for any $t>0$ we conclude that
$
\Vert P_t(V)\Vert<\infty
$.

\subsection{Multivariate birth and death processes}
We  denote by $e:=\{e_i,~1\leq i\leq n\}$ the collection of column vector $e_i$ on $\{0,1\}^n$ with entries $e_i(j)=1_{i=j}$ and with a slight abuse of notation we denote by $0$ the null state in $\NN^n$.
Let $X_t(x)$ be a stochastic flow on $E=\NN^n-\{0\}$ with generator $L$ defined by
\begin{equation}\label{def-L-birth-death}
L(f)(x):=\sum_{y\in E}J(x,y)~(f(y)-f(x)).
\end{equation}
Let $\lambda,\mu,\upsilon,\varsigma$ be  some column vectors and let $C,D$ some $(d\times d)$-matrices with real entries  such that for any $1\leq i\leq d$ and any
$x\in E$ we have
$$
J(x,x+e_i):=\upsilon_i+x_i~(\lambda_i+(C x)_i)\geq 0\quad \mbox{\rm and}\quad 
J(x,x-e_i):=\varsigma_i+x_i~(\mu_i+(D x)_i)\geq 0.
$$
We also set
$$
J(x,y)=0\quad 
\mbox{as soon as}\quad   |x-y|\geq 2.
$$
We further assume that
$$
\vert \upsilon\vert\geq \vert \varsigma\vert \qquad 
B:=(D-C)\geq b~I>0 \quad \mbox{for some $b>0$. }
$$
and we set
$$
a_0:=\vert \upsilon\vert-\vert \varsigma\vert \geq 0
\qquad
a_1:=\vee_{1\leq i\leq n} (\lambda_i-\mu_i)$$
and for any $x\in\NN^n$
$$
\Vert x\Vert:=\left(\sum_{1\leq i\leq n}x_i^2\right)^{1/2}\geq \vert x\vert:=\sum_{1\leq i\leq n}x_i.
$$

 Consider the Lyapunov function
$$
x\in E\mapsto V(x)=\vert x\vert \in \NN_+.
$$
Note that $V$ is locally bounded with finite  level sets  and for
 any $x\in E-e$ we have
\begin{eqnarray*}
L(V)(x)&=&\sum_{1\leq i\leq n}
\left(\left(\upsilon_i+x_i~(\lambda_i+(C x)_i)\right)-\left(\varsigma_i+x_i~(\mu_i+(D x)_i)\right)\right).
\end{eqnarray*}
In this situation, we have the formula
\begin{equation}\label{hypothesis-birth-death}
L(V)(x)=a_0+(\lambda-\mu)^{\prime} x-x^{\prime}Bx\leq a_0+a_1\vert x\vert-b\Vert x\Vert^2.
\end{equation}
On the other hand, for any $y=e_j$ we have
$$
L(V)(y)=\sum_{1\leq i\leq n} J(y,y+e_i).
$$
This implies that
\begin{eqnarray*}
P_t(L(V))&=&P_t(1_{E-e}~L(V))+P_t(1_{e}~L(V))\\
&=&a_0+a_1~P_t(V)-b~P_t(V^2)\\
&&\hskip.1cm +\sum_{1\leq j\leq n}~P_t(1_{e_j})~\left(
L(V)(e_j)-\left(a_0+(\lambda-\mu)^{\prime} e_j-e_j^{\prime}Be_j\right)\right)
\end{eqnarray*}
from which we readily check that
$$
\partial_t\EE(V(X_t(x)))\leq a^+_0+a_1\EE(V(X_t(x)))-b~(\EE(V(X_t(x))))^2
$$
with
\begin{eqnarray*}
\displaystyle a^+_0&:=&a_0+\sum_{1\leq j\leq d}~\left(
\sum_{1\leq i\leq d} J(e_j,e_j+e_i)-\left(\vert \upsilon\vert-\vert \varsigma\vert+(\lambda-\mu)^{\prime} e_j-e_j^{\prime}(D-C)e_j\right)\right)\\
&=&a_0+\sum_{1\leq j\leq d} \left(\vert \varsigma\vert+\mu^{\prime} e_j+e_j^{\prime}D e_j\right).
\end{eqnarray*}
We conclude that
$\Vert P_t(V)\Vert<\infty$, for any $t>0$. The semigroup analysis discussed above can be extended without difficulties to more general process on countable spaces models satisfying condition (\ref{hypothesis-birth-death}). The extension to time varying models can also be handle using a more refined analysis on time varying Riccati equations.
 
We also mention, that the case $\vert \upsilon\vert=0=\vert \varsigma\vert$ coincides with the competitive and multivariate 
Lotka-Volterra birth and death process discussed in Theorem 1.1 in~\cite{champagnat-6}.

\section{Some conditional diffusions}\label{conditional-diff-sec}

\subsection{Coupled harmonic oscillators}\label{muli-dim-osc-harm}
Consider the $\RR^n$-valued diffusion (\ref{ref-X-diff-stable}) with $(b(x),\sigma(x))=(A x,\Sigma)$, 
for  some  {\em non necessarily stable} drift matrix $A$ and some diffusion matrix $\Sigma$ with appropriate dimensions. We associate with a given semi-definite positive $(n\times n)$ matrix $S\geq 0$ the potential function
\begin{equation}\label{ref-Ua-q}
U(x):=\frac{1}{2}~x^{\prime}Sx\quad \mbox{\rm and we set}\quad R=\Sigma\Sigma^{\prime}.
\end{equation}
We assume that the pairs $(A,R^{1/2})$  and $(A^{\prime},S^{1/2})$ are both controllable.
Let $Q_t=Q_t^{[U]}$ be the sub-Markov semigroup defined in (\ref{Q-U-gen}) on the Euclidean space $\Ea=E=\RR^n$.
As shown in~\cite{dm-h-21}, the leading-triple $(\rho,h,\eta_{\infty})$ discussed in (\ref{intro-rho-h}) is given by 
\begin{eqnarray}
\rho&=&-\tr(R q_{\infty})/2=-\tr(p_{\infty} S)/2\nonumber\\
h(x)&=&\exp{\left(-x^{\prime}q_{\infty}x/2\right)}\quad\mbox{\rm and}\quad
\eta_{\infty}(dx)=\frac{\exp{\left(-x^{\prime}p_{\infty}^{-1}x/2\right)}}{\sqrt{\mbox{\rm det}(2\pi p_{\infty})}}~dx,\label{def-rho-h-2d}
\end{eqnarray}
with the positive fixed points $p_{\infty}$ and $q_{\infty}$ of the dual algebraic Riccati matrix equation
$$
Ap_{\infty}+p_{\infty}A^{\prime}+R-p_{\infty}Sp_{\infty}=0\quad \mbox{\rm and}\quad
A^{\prime}q_{\infty}+q_{\infty}A+S-q_{\infty}R q_{\infty}=0.
$$
In this context,  the $h$-process, denoted $(X_t^h(x))_{t \ge 0}$ and defined by the
stochastic differential equation
\begin{equation}\label{def-X-h-2d}
dX^h_{t}(x)=A^hX^h_{t}(x)~dt+\Sigma~dB_t\quad \mbox{\rm with}\quad A^h:=A-R~q_{\infty}.
\end{equation}
Our controllability conditions ensures that $A^h$ is a stable matrix.
Note that $X^h_t(x)$ is an $\RR^n$-valued Gaussian random variable with mean
$m^h_t(x)$ and covariance matrix $p^h_t \in \RR^{n\times n}$ given for any $t>0$ by
$$
m^h_t(x)=\exp{\left(A^ht\right)}~x\quad \mbox{\rm and}\quad
p^h_t=\int_0^t~\exp{\left(A^hs\right)}~R\exp{\left((A^h)^{\prime}s\right)}~ds >0.
$$
This yields the explicit formula
$$
P^h_t(x,dy)=\frac{1}{\sqrt{\mbox{\rm det}(2\pi p^h_t)}}~\exp{\left(-\frac{1}{2}(y-m_t^h(x))^{\prime}(p^{h}_t)^{-1}
(y-m_t^h(x))\right)}~dy.
$$
Moreover the invariant measure $\eta^h_{\infty}=\eta^h_{\infty} \Pa^h_t$ is unique and given by
$$
\eta^h_{\infty}(dx)=\frac{1}{\sqrt{\mbox{\rm det}(2\pi p^h_{\infty})}}~\exp{\left(-\frac{1}{2}y^{\prime}(p^{h}_{\infty})^{-1}
y\right)}~dy
$$
with the limiting covariance matrix
$$
p^h_{\infty}:=\int_0^{\infty}~\exp{\left(A^hs\right)}\Sigma^2\exp{\left((A^h)^{\prime}s\right)}~ds =(p_{\infty}^{-1}+q_{\infty})^{-1}>0.
$$
For any time horizon $t\geq 0$ and any measurable function $F$ on the set $\Ca([0,t],\RR^n)$ of continuous paths from $[0,t]$ into $\RR^n$ we have the 
path space exponential change of measure Feynman-Kac formula
$$
\EE\left(F(\XX_t(x))\exp{\left(\int_0^t\UU_s(\XX_s(x))~ds\right)}\right)=e^{\rho t}~h(x)~\EE\left(F(\XX^h_t(x))/h(X^h_t(x))\right)
$$
with the historical processes
$$
\XX_t(x):=(X_s(x))_{0\leq s\leq t}, \qquad \XX^h_t(x):=(X^h_s(x))_{0\leq s\leq t }\quad\mbox{\rm and}\quad \UU_s(\XX_s(x)):=U_s(X_s(x)).
$$
This yields the conjugate formulae
$$
Q_t(f)=e^{\rho t}~h~P^h_t(f/h).
$$
We denote by $(m_t(x),p_t)\in (\RR^n\times \RR^{n\times n})$  the mean and covariance parameters satisfying the linear evolution and the Riccati matrix differential equations
\begin{equation}\label{ref-riccati}
\left\{
\begin{array}{rcl}
\partial_t m_t(x)&=&(A-p_tS)~m_t(x)\\
&&\\
\partial_t p_t&=&Ap_t+p_tA^{\prime}+\Sigma^2-p_tSp_t\quad \mbox{with}\quad (m_0(x),p_0)=(x,0).
\end{array}
\right.
\end{equation}

The next proposition provides an explicit description of these semigroups.
\begin{prop}[\cite{dm-h-21}]\label{ref-coupled-harmonic-explicit}
For any time horizon $t>0 $ we have $p_t>0$ and the normalized semigroup $\overline{Q}_t$ defined in (\ref{def-Phi-s-t-intro-0}) is given by
\begin{equation}\label{solv-harmonic-P-2d}
\overline{Q}_t(x,dy)=\frac{1}{\sqrt{\mbox{\rm det}(2\pi p_t)}}~\exp{\left(-\frac{1}{2}(y-m_t(x))^{\prime}p^{-1}_t
(y-m_t(x))\right)}~dy
\end{equation}
as well as
$$
-2\log{Q_t(1)(x)}=x^{\prime}\left(\int_0^t~F^{\prime}_sSF_s~ds\right) x+\int_0^t~\tr(Sp_s)~ds
$$
with the fundamental matrix semigroup $F_t$ starting at $F_0=I$ given by
$$
\partial_t F_t=(A-p_tS)~F_t.
$$
\end{prop}

Observe that the normalized Markov operator $\overline{Q}_t$ satisfies (\ref{def-OU-A-Sigma-lem}) and (\ref{def-OU-A-Sigma-lem-epsilon}) with  the parameters
\begin{equation}\label{sub-gauss-harmonic-ref}
c_{t}=\frac{1}{\sqrt{\mbox{\rm det}(2\pi p_{t})}}, \, 
\sigma^2_{t}=\lambda_{\tiny max}(p_t)\quad \mbox{and}\quad
\epsilon_{\tau}=\vert e^{\tau(A-p_{\infty}S)}\vert\longrightarrow 0 \,\,  \mbox{\rm as} \,\,  {\tau\rightarrow\infty}
\end{equation}
for some matrix norm $\vert \point\vert$.
The r.h.s. assertion is a direct consequence of the Floquet representation theorem presented in~\cite{bdp-21} (cf. (1.3) and Theorem 1.1) and  the fact that $(A-p_{\infty}S)$ is a stable matrix. Applying Lemma~\ref{ref-lem-OU-A-Sigma}
 for any $v\geq 0$ and $t>0$ there also exists some finite constant $\delta_t>0$ such that
$$
V(x):=\exp{\left(v\vert x\vert\right)}\Longrightarrow
\overline{Q}_{t}(V)/V\leq c_{t}/V^{\delta_{t}}. 
$$
Using Proposition~\ref{ref-coupled-harmonic-explicit}, for any $k\geq 0$ and $t\geq 0$ it is also readily checked that
$$
V(x):=(1+\Vert x\Vert)^k\Longrightarrow
\Vert\overline{Q}_{t}(V)/V\Vert<\infty\quad\mbox{\rm and}\quad
\Vert Q_{t}(V)\Vert<\infty.
$$

\subsection{Half-harmonic linear diffusions}\label{ref-sec-lin-half}
For one dimensional models,  the coupled harmonic oscillator discussed in Section~\ref{muli-dim-osc-harm} resumes to one dimensional  linear diffusion
\begin{equation}\label{ref-lin-diff}
dX_t(x)=a\,X_t(x)\,dt+dB_t\quad\mbox{and the potential}\quad U(x)=\varsigma x^2/2
\end{equation}
for some parameters $\varsigma>0$ and $a\in\RR$. We set
$\beta:=a+\sqrt{a^2+\varsigma}$. In this notation, the leading pair $(\rho,h)=(\rho_1,\varphi_1)$ is given by 
\begin{equation}\label{ref-beta-h}
\rho=-\beta/2 \quad\mbox{\rm and}\quad h(x)=\left((\beta-a)/\pi\right)^{1/4}\exp{\left(-\beta x^2/2\right)}.
\end{equation}
The quasi-invariant measure is therefore given by
$$
\eta_{\infty}(dx)=\sqrt{\frac{\varsigma}{2\pi\beta}}~
\exp{\left(- \varsigma x^2/(2\beta)\right)}~dx.
$$
Therefore,  the $h$-process resumes to the Ornstein-Uhlenbeck diffusion
\begin{equation}\label{ref-h-proc-lin}
dX^h_t(x)=-b~
X^h_t(x)~dt+dB_t\end{equation}
with the invariant measure
$$
 \eta^h_{\infty}(dx):=\sqrt{\frac{b}{\pi}}~
\exp{\left(-b~ x^2\right)}~dx\quad \mbox{with}\quad b:=(\beta-a)=\sqrt{a^2+\varsigma}>0.
$$
Note that any Ornstein-Uhlenbeck process can be seen as the $h$-process associated with a non absorbed (possibly transient) linear diffusion evolving in some quadratic potential  well.

In this context, Proposition~\ref{ref-coupled-harmonic-explicit} is also satisfied 
with the mean and variance parameters
\begin{equation}\label{ref-m-sigma-lin}
\left\{
\begin{array}{rcl}
\partial_t m_t(x)&=&(a-p_t\varsigma)~m_t(x)\\
&&\\
\partial_t p_t&=&2ap_t+1-\varsigma p_t^2\quad \mbox{with}\quad (m_0(x),p_0)=(x,0).
\end{array}\right.
\end{equation}
We also have
\begin{equation}\label{ref-Q1-proc-lin}
\begin{array}{l}
\displaystyle-\frac{2}{\varsigma}\log{\Qa_t(1)(x)}=\overline{p}_t+\chi_t~x^2
\end{array}\end{equation}
with
$$
\chi_t:=\int_0^t~\exp{\left(-2\int_0^s (a-p_u\varsigma)du\right)}~ds\quad\mbox{and}\quad \overline{p}_t:=
\int_0^tp_sds.
$$

The half-harmonic semigroup  associated with the flow $X_t(x)$ is defined for any $x\in E:=]0,\infty[$ and $f\in \Ba_b(E)$ by the formulae
  \begin{eqnarray}
Q_{t}(f)(x)&:=&\EE\left(f(X_{t}(x))~ 1_{T(x)>t}~\exp{\left\{-\int_0^t U(X_s(x))~ds\right\}}\right).\label{ref-half-lin-diff-harmonic}\end{eqnarray}
In the above display, $T(x)$ stands for the hitting time of the origin. In terms of the $h$-process of the flow in the harmonic potential (\ref{ref-h-proc-lin}) we 
 also have the conjugate formula
  \begin{eqnarray}\label{lin-diff-half}
Q_{t}(f)(x)
&=&e^{t\rho}~e^{-\beta x^2/2}~\EE\left(f(Y_t(x))~e^{\beta Y_t(x)^2/2}~1_{T^Y(x)>t}\right)
\end{eqnarray}
with the parameters $(\rho,\beta)$ defined in (\ref{ref-beta-h}) and the Ornstein-Uhlenbeck diffusion flow defined by
$$
dY_t(x)=-b~
Y_t(x)~dt+dB_t\quad \mbox{\rm with}\quad b:=(\beta-a)>0.
$$
In the above display,  $T^Y(x)$ stands for the hitting time of the origin by the flow $Y_t(x)$ starting at $x>0$.
Arguing as in Section~\ref{half-oscillator-sec} we check that
$$
\begin{array}{l}
\displaystyle
Q_t(x,dy)\\
\\
=\sinh{(y~m_t(x))}~\exp{\left(-\frac{\varsigma}{2}\left(\chi_t\,x^2+\overline{p}_t\right)\right)}
~\displaystyle\times \sqrt{\frac{2}{\pi p_t}}~\exp{\left(-\frac{y^2+m_t(x)^2}{2p_t}\right)}~1_{]0,\infty[}(y)~dy
\end{array}
$$
with the parameters $(m_t(x),p_t)$ and $(\chi_t,\overline{p}_t)$ defined in (\ref{ref-m-sigma-lin}) and (\ref{ref-Q1-proc-lin}). 

Arguing as in (\ref{half-harm-lyap}), choosing the Lyapunov function
$
V(x)=x^n+1/x$, for some $n\geq 1$, we readily check that
\begin{equation}\label{half-harm-lyap-l}
V\in \Ca_{\infty}(E)\quad \mbox{\rm and}\quad
Q_t(V)/V\leq c_t/V~\in \Ca_0(E).
\end{equation}

\subsection{Linear diffusions in some domains}\label{non-abs-OU}
Consider the one-dimensional stochastic flow $Y_t(x)$ of an Ornstein-Uhlenbeck
$$
dY_t(x)=-b~
Y_t(x)~dt+dB_t\quad \mbox{for some}\quad b>0.
$$
In the above display, $B_t$ is a one-dimensional Brownian motion starting at the origin.
For a given $x\in E:=]0,\infty[$, we let $T^Y(x)$ be the hitting time of the origin by the flow $Y_t(x)$ starting at $x>0$. Consider the semigroup
$$
Q^Y_{t}(f)(x):=\EE(f(Y_{t}(x))~ 1_{T^Y(x)>t}).
$$
Choosing $(a,\varsigma,\beta,\rho)=(0,b^2,b,-b/2)$ in  (\ref{ref-beta-h}), formula (\ref{lin-diff-half}) takes the form
  $$
  Q^Y_t(f)(x)=
e^{-\rho t}~H(x)^{-1}Q_{t}(fH)(x)\quad \mbox{\rm with}\quad
H(x)=\exp{\left(-b x^2/2\right)}
$$
with the semigroup $Q_t$ defined in (\ref{ref-half-lin-diff-harmonic}) with $U(x)=b^2 x/2$.

For any given $n\geq 1$ we have
$$
V(x):=x^n+{1}/{x}\Longrightarrow V\in \Ca_{\infty}(E)\quad \mbox{\rm and}\quad
 V^{H}:=V/H\in  \Ca_{\infty}(E).
$$
Using (\ref{half-harm-lyap-l}) we conclude that
$$
V^{H}\in \Ca_{\infty}(E)\quad \mbox{\rm and}\quad
 Q^Y_t(V^{H})/V^{H}=e^{-\rho t}~Q_{t}(V)/V\leq c_t/V\in \Ca_0(E).
$$

The long time behavior of the positive semigroup $Q^Y_t$  is also studied in~\cite{lladser}, and more recently in~\cite{ocafrain} in terms of the tangent of the $h$-process.

More generally, consider the $\RR^n$-valued diffusion flow $X_t(x)$ defined in (\ref{ref-X-diff-stable}) with $(b(x),\sigma(x)=(A x,\Sigma)$, 
for some matrices $(A,\Sigma)$ with appropriate dimensions. Assume that
 $R:=\Sigma\Sigma^{\prime}$ is positive semi-definite 
and the pair of matrices $(A,R^{1/2})$ are controllable. In this situation, the Markov semigroup $P_t$ of the stochastic flow $X_t(x)$ satisfies the sub-Gaussian estimate (\ref{ref-sub-gauss-b}) for some parameters $(\sigma_t,m_t(x))$. 

Consider a  domain $E\subset \RR^n$ with  $\Ca^2$-boundary with uniformly bounded interior curvature.
For any given $x\in E$, let $Q_t$ be the sub-Markov semigroup
\begin{equation}\label{sub-Markov-OU}
Q_{t}(f)(x):=\EE(f(X_{t}(x))~ 1_{T(x)>t})\quad \mbox{\rm with}\quad
T(x):=\inf\left\{t\geq 0~:~X_{t}(x)\in \partial E\right\}.
\end{equation}
We clearly have
$
Q_t(V_{\partial})\leq P_t(V_{\partial})
$, with the function $V_{\partial}$ defined in (\ref{def-V-delta}). When $E$ is non necessarily bounded  but its boundary $\partial E$ is bounded we known from (\ref{ref-bounded-boundary}) that  $\Vert Q_t(V_{\tiny \partial})\Vert<\infty$. For non necessarily bounded boundaries the sub-gaussian property 
(\ref{ref-sub-gauss-b}) ensures that $\Vert Q_t(V_{\tiny \partial})\Vert<\infty$.
 
When $E$ is bounded, applying Lemma~\ref{lem-lip-dom-compact} (see also Proposition~\ref{prop-domain-bounded}) we have
$$
V_{\partial}\in \Ca_{\infty}(E)\quad \mbox{and}\quad
 Q_t(V_{\partial})/V_{\partial}\leq c_t/V_{\partial}\in \Ca_{0}(E).
$$
For unbounded domains  we need to ensure that $A$ is stable so that (\ref{stable-exp}) is satisfied for some norm $\vert \point\vert$ on $\RR^n$. In this situation, applying Proposition~\ref{prop-ref-lang} for any $t>0$ there exists some $\delta_t>0$ such that
$$
V_{\Ea}(x):=\exp{\left(v \vert x\vert\right)}\Longrightarrow
 Q_t(V_{\Ea})/V_{\Ea}\leq c_t/V_{\Ea}^{\delta_t}.
$$
Applying Proposition~\ref{cond-gen-unb-prop} with $\epsilon=0$, for any $p>1$ we conclude that
\begin{equation}\label{ref-OU-Lyap-dom}
V_p:=V_{\Ea}^{1-1/p}~V_{\partial}^{1/p}\in\Ca_{\infty}(E)
\quad\mbox{\rm
and}\quad
Q_t\left(V_p\right)/V_p\leq c_t ~\Theta_{p,t}
 \end{equation}
 with the function
 $$
 \Theta_{p,t}:=1/(V_{\Ea}^{\delta_t (1-1/p)}V_{\partial}^{1/p})\in \Ca_0(E). 
 $$

\subsection{Langevin diffusions in some domains}
Consider the semigroup $Q_t$ of the  one-dimensional Langevin diffusion defined in (\ref{ref-langevin-one-D-examp}) with $E=]0,\infty[$ and a quadratic confinement potential
$$
W(x)=x^2/2\Longrightarrow H(x):=e^{-W(x)}=e^{-x^2/2}\quad \mbox{and}\quad
U:=\frac{1}{2}~\left(x^2+1\right). 
$$
In this case, the semigroup $\Qa_t$ defined in (\ref{ref-Q-U-half-harmonic}) coincides with the semigroup of the half-harmonic oscillator discussed in Section~\ref{half-oscillator-sec}. By (\ref{half-harm-lyap})   for any $n\geq 1$ we have
$$
V(x):=x^n+{1}/{x}\Longrightarrow \Qa_t(V)/V\leq c_t/V~\in \Ca_0(E).
$$
Notice that
\begin{equation}\label{ref-U2}
V^H(x):=V(x)/H(x)=x^n~e^{x^2/2}+\frac{e^{x^2/2}}{x}.
\end{equation}
Using (\ref{ref-Q-U-half-harmonic}) we conclude that
$$
V^H\in \Ca_{\infty}(E)\quad \mbox{\rm and}\quad
Q_t(V^H)/V^H=\Qa_t(V)/V\leq c_t/V\in \Ca_0(E).
$$
More generally, consider the case $E=]0,\infty[$ with at least a quadratic confinement potential $U$, in the sense that
$$
U(x)=\frac{1}{2}~\left((\partial W)^2-\partial^2W\right)(x)
\geq ~U_2(x):=c+\varsigma ~x^2/2\quad \mbox{\rm for some $\varsigma>0$.}
$$
In this situation,
$
\Qa_t\ll\Qa^{[U_2]}
$ is dominated by
 the semigroup $\Qa^{[U_2]}$ of the half-harmonic oscillator discussed in Section~\ref{half-oscillator-sec}. Arguing as in (\ref{ref-U2}) we have
$$
H:=e^{-W}\qquad
V^H:=V/H\in \Ca_{\infty}(E)\quad \mbox{\rm and}\quad
Q_t(V^H)/V^H\leq c_t/V\in \Ca_0(E).
$$
For instance, whenever the confinement potential $W$ is chosen so that
$$
 W(x)\geq \epsilon_0 \log{x}+W_1(x)\quad \mbox{ for some $0\leq \epsilon_0<1$}
 $$
 and some function $W_1\geq 1$ such that $W_1(x)\longrightarrow_{x\rightarrow\infty}\infty$
we have
$$
H=e^{-W}\Longrightarrow
V^H(x):=V(x)/H(x)=x^n~e^{W(x)}+\frac{e^{W(x)}}{x}\geq x^n~e^{W(x)}+\frac{e^{W_1(x)}}{x^{1-\epsilon_0}}.
$$
Using (\ref{ref-Q-U-half-harmonic}) we conclude that
$$
V^H\in \Ca_{\infty}(E)\quad \mbox{\rm and}\quad
Q_t(V^H)/V^H=\Qa_t(V)/V\leq c_t/V\in \Ca_0(E).
$$
We illustrate the above result, with the logistic diffusion discussed in~\cite{cattiaux-09}.
Consider the generalized Feller diffusion
$$
dY_t(x):= \left(2a~Y_t(x)-(8b/\sigma^2)~Y_t(x)^2\right)~dt+\sigma~\sqrt{Y_t(x)}~dB_t
$$
starting at $x\in E:=]0,\infty[$. In the above display, $B_t$ is a one dimensional Brownian motion starting at the origin and $a,b,\sigma>0$ some parameters. Observe that
$$
X_t(x):=(2/\sigma)\sqrt{Y_t(x)}\Longrightarrow  dX_t(x)=-\partial W\left(X_t(x)\right)~dt+dB_t
$$
with the potential function
$$
\partial W(x)=\frac{1}{2x}-a~x+b~x^3\quad \mbox{\rm with}\quad a,b>0.
$$
Thus, choosing
$$
W(x)=\frac{1}{2}~\log{x}+b\frac{x^4}{4}-a\frac{x^2}{2}$$
we readily check that
$$V^H(x):=e^{W(x)}(x^n+1/x)=(x^{n-1/2}+1/\sqrt{x})~e^{b\frac{x^4}{4}-a\frac{x^2}{2}}\Longrightarrow V^H\in \Ca_{\infty}(E).
$$

More generally, consider the  Langevin diffusion flow $$\Xa_t(z)=(X_t(z),Y_t(z))\in (\RR^{n}\times\RR^n)$$ starting at
$z=(x,y)\in (\RR^{n}\times\RR^n)$ and defined  by the hypo-elliptic  diffusion (\ref{langevin-intro-ref}). We further assume that $\sup_{D}a<\infty$ for some bounded open connected domain $D\subset\RR^n$ with  $\Ca^2$-boundary, and for any  $z\in E:=D\times \RR^n$ and $f\in \Ba_b(E)$ we set
$$
\Qa_t(f)(z):=\EE\left(f(\Xa_t(z))~1_{T(z)>t}\right)
\quad\mbox{\rm with}\quad
 T(z):=\inf\left\{t\geq 0~:~X_t(z)\in \partial D\right\}.
$$ 
We know from (\ref{langevin-ref-intro}) that for any $q>1$ we have
$
\Qa\ll_qQ
$ is $q$-dominated by the sub-Markov semigroup $Q_t$ associated with the Ornstein-Uhlenbeck diffusion on $E$ defined in (\ref{sub-Markov-OU}), with the matrices $(A,\Sigma)$ defined in (\ref{ref-A-Sigma}). In terms of the functions $(V_p,\Theta_{p,t})$ defined in (\ref{ref-OU-Lyap-dom}), combining (\ref{ref-dom-il}) with (\ref{ref-OU-Lyap-dom})
for any $p,q>1$ we conclude that
$$
\Qa_t(V_{p,q})/V_{p,q}\leq c_t(p,q)~~\Theta_{p,q,t}
$$
with the collection of Lyapunov functions
$$
V_{p,q}:=V_p^{1/q}\in \Ca_{\infty}(E)\quad\mbox{\rm and the function}\quad
\Theta_{p,q,t}:=\Theta_{p,t}^{1/q}\in \Ca_0(E).
$$

\subsection{Coupled oscillators in some domains}\label{one-wall-sec}
Consider the $\RR^n$-valued diffusion $X_t(x)$ and the quadratic potential function $U$
discussed in Section~\ref{muli-dim-osc-harm}, for some $n\geq 2$ and set $E:=]0,\infty[\times\RR^{n-1}$.
Let $Q_t$ be the semigroup  defined for any $f\in \Ba_b(E)$ and $x\in E$ by the formulae
  \begin{eqnarray}
Q_{t}(f)(x)&:=&\EE\left(f(X_{t}(x))~ 1_{T(x)>t}~\exp{\left(-\int_0^t U(X_s(x))ds\right)}\right)\end{eqnarray}
with the quadratic function $U$ in (\ref{ref-Ua-q}) and the exit time $T(x)$ given by
$$
T(x):=\inf\left\{t\geq 0~:~X_{t}(x)\in \partial E\right\}\quad \mbox{\rm with}\quad \partial E=\{0\}\times\RR^{n-1}.
$$
 In terms of the $h$-process $Y_t(x):=X^h_t(x)$ associated with the leading pair $(\rho,h)$ defined in  (\ref{def-X-h-2d})  we 
 also have the conjugate formula
$$
Q_{t}(f)
=e^{\rho t}~h~Q^{Y}_t(f/h)
\quad \mbox{with}\quad
Q^{Y}_t(f)(x):=\EE\left(f(Y_t(x))~1_{T^Y(x)>t}\right).
$$
In the above display, $T^Y(x)$ stands for the boundary hitting time
$$
T^{Y}(x):=\inf\left\{t\geq 0~:~Y_{t}(x)\in \partial E\right\}.
$$

When $n=2$, the linear diffusion $X_t(x)$ associated to  the matrices $A_{1,2}=A_{2,1}=A_{2,2}=0$ and $A_{1,2}=1$ and $\Sigma_{1,1}=\Sigma_{1,2}=\Sigma_{2,1}=0$ and $\Sigma_{2,2}=1$ coincides with the integrated Wiener process model discussed in~\cite{goldman,lachal,mcKean-63}. In a seminal article~\cite{mcKean-63}, McKean obtained the joint distribution of the pair $(T(x),X^2_{T(x)})$ in the absence of soft absorption, that is when $U=0$. To the best of our knowledge, an explicit description of the distribution of this pair and the corresponding sub-Markov semigroup is unknown in  
more general situations.

Observe that for any $x\in E$ and any non negative function $f\in \Ba_b(\RR^n)$ we have
$$
Q_{t}(f)(x)\leq  \Qa_t(f)(x):=e^{\rho t}~h(x)~\EE\left(f(Y_t(x))/h(Y_t(x))\right).
$$
The semigroup $\Qa_t$ defined above coincides with the semigroup of the coupled harmonic oscillator discussed in  Section~\ref{muli-dim-osc-harm}. We know from (\ref{sub-gauss-harmonic-ref}) that $\overline{\Qa}_t$ satisfies the sub-Gaussian estimates (\ref{def-OU-A-Sigma-lem}) with
$$
c_{t}=\frac{1}{\sqrt{\mbox{\rm det}(2\pi p_{t})}}\quad\mbox{\rm and}\quad
\sigma^2_{t}=\lambda_{\tiny max}(p_t)
$$
with the solution $p_t$ of the Riccati-matrix equation (\ref{ref-riccati}).
Using Proposition~\ref{ref-coupled-harmonic-explicit}  for any $k\geq 1$ we have
$$
V_{\Ea}(x):=1+\Vert x\Vert^k\Longrightarrow
\Vert\overline{\Qa}_t(V_{\Ea})/V_{\Ea}\Vert<\infty\Longrightarrow  Q_t(V_{\Ea})\leq c_t~\Qa_t(1) V_{\Ea}.
$$
Recalling that $\Qa_t(1)(x)$ tends to $0$ exponentially fast as $\Vert x\Vert\rightarrow\infty$, this implies that 
$$\forall t>0\qquad\Vert Q_t(V_{\Ea}) \Vert<\infty.
$$

On the other hand for any $y=(y_1,y_{-1})\in E=(]0,\infty[\times\RR^{n-1})$, the distance to the boundary is given by 
$d(y,\partial E)=y_1$. In terms of the function $V_{\partial}$ defined in (\ref{def-V-delta}) his implies that
$$
Q_t(V_{\partial})(x)\leq \int~\Qa_t(x,dy)~1_{]0,1[}(y_1)~\cchi(y_1)+\cchi(1)~\Qa_t(1)(x)$$
from which we check that $\Vert Q_t(V_{\partial})\Vert<\infty$. Applying Proposition~\ref{prop-bb-ub} we conclude that
$$
V:=V_{\partial}+V_{\Ea}\in \Ca_{\infty}(E)\quad \mbox{\rm and}
 \quad Q_t(V)/V\leq c_t/V\in \Ca_0(E).
$$
The same analysis applies by replacing the half line $E_1$ by the unit interval $E_1:=]0,1[$. In this context, the boundary is given by the two infinite potential walls
$$
\partial E=(\{0\}\times\RR^{n-1})\cup (\{1\}\times\RR^{n-1})
\quad\mbox{and}\quad
d(x,\partial E)=x_1\wedge (1-x_1).
$$
More generally, consider a  domain $E\subset \RR^n$ with  $\Ca^2$-boundary with uniformly bounded interior curvature.
In this situation, the sub-Gaussian property (\ref{ref-sub-gauss-b}) ensures that
$
\Vert\overline{\Qa}_t(V_{\partial})\Vert<\infty
$ and therefore
$$\Vert Q_t(V_{\partial})\Vert\leq \Vert \Qa_t(V_{\partial})\Vert= \Vert \Qa_t(1)\overline{\Qa}_t(V_{\partial})\Vert<\infty.
$$
 Applying Proposition~\ref{prop-bb-ub}, we conclude that 
 $$
 V:=V_{\partial}+V_{\Ea}\in \Ca_{\infty}(E)\quad \mbox{\rm and}
\quad \mbox{\rm and}
 \quad Q_t(V)/V\leq c_t/V\in \Ca_0(E).
$$

\section{Some hypersurface boundaries}\label{hyper-surface-bound-sec}
\subsection{Defining functions and charts}
Consider a smooth function $y\in \RR^{n-1}\mapsto \varphi(y)\in \RR$ with non empty and connected level set, for some $n\geq 2$.  Consider a domain $E$ in $\RR^n$ with a smooth boundary $\partial E=\overline{\varphi}^{-1}(\{0\})$ defined as the null level set  of the function
$$
x=(x_i)_{1\leq i\leq n}\in \RR^{n}\mapsto
 \overline{\varphi}(x):=\varphi(x_{-n})-x_n\quad \mbox{\rm with}\quad
 x_{-n}:=(x_i)_{1\leq i<n}\in \RR^{n-1}.
$$
Consider the column vectors $\nabla \varphi(x_{-n}):=\left(\partial_{x_i}\varphi(x_{-n})\right)_{1\leq i<n}$.
In this notation, the unit normal vector $N(x)$ at $x\in \partial E$ is given by the column vectors
$$
N(x)=\frac{\nabla\overline{\varphi}(x)}{\Vert \nabla\overline{\varphi}(x)\Vert} = \frac{1}{\sqrt{1+\Vert \nabla \varphi(x_{-n})\Vert^2}}~\left(
\begin{array}{c}
\nabla \varphi(x_{-n})\\
-1
\end{array}
\right).$$
Observe that the vector $N(x)$ is the outward-pointing normal direction to $E$ as soon as $E=\overline{\varphi}^{-1}\left(]-\infty,0[\right)$
and the inward-pointing normal direction to $E$ when $E=\overline{\varphi}^{-1}\left(]0,+\infty[\right)$.

Consider the column vectors  $e_i:=(1_i(j))_{1\leq i<n}$, with $1\leq i<n$.
In this notation, the $(n-1)$ tangential column vectors  $T_{i}(x)$ 
 at $x\in\partial E$ are given for any $1\leq i<n$ by the column vectors
$$
 T_i(x):=\left(
\begin{array}{c}
e_i\\
\partial_{x_i}\varphi(x_{-n})
\end{array}
\right).$$

The inner product $g(x)$ on the tangent space $\TT_x(\partial E)$  (a.k.a. the first fundamental form on $\partial E$) is given by the Gramian matrix
$$
g(x)=\left(T_i(x)^{\prime}T_j(x)\right)_{1\leq i,j<n}=T(x)T(x)^{\prime}\quad\mbox{\rm
with}
\quad
T(x)^{\prime}:=\left(T_1(x),\ldots,T_{n-1}(x)\right).
$$
This yields the matrix formula
$$
g(x)=\left(\begin{array}{c}
I,
 \nabla\varphi(x_{-n})\end{array}\right)\left(\begin{array}{c}
 I\\
 \nabla\varphi(x_{-n})^{\prime}\end{array}\right)=I+\nabla\varphi(x_{-n})\nabla\varphi(x_{-n})^{\prime}.
$$
In this notation, the projection $\mbox{\rm proj}_{\TT_x(\partial E)}$ on the tangent space $\TT_x(\partial E)$ is defined for any column vector $V=(V^i)_{1\leq i\leq n}\in \RR^n$ by
$$
\mbox{\rm proj}_{\TT_x(\partial E)}(V):=\left(T_1(x),\ldots,T_{n-1}(x)\right)g(x)^{-1}\left(
\begin{array}{c}
T_1(x)^{\prime}\\
\vdots\\
T_{n-1}(x)^{\prime}\end{array}
\right)\left(
\begin{array}{c}
V^1\\
\vdots\\
V^n\end{array}
\right).
$$
In matrix notation, the projection of $m$ column vectors $V_i\in\RR^n$, with $i\in \{1,\ldots,m\}$ and any $m\geq 1$ takes the synthetic form
\begin{eqnarray*}
\mbox{\rm proj}_{\TT_x(\partial E)}(V_1,\ldots,V_m)&=&
\left(T(x)^{\prime}g(x)^{-1}T(x)\right)(V_1,\ldots,V_m)\\
&=&\left(
\mbox{\rm proj}_{\TT_x(\partial E)}(V_1),\ldots,\mbox{\rm proj}_{\TT_x(\partial E)}(V_m)\right).
\end{eqnarray*}
Equivalently, if $g(x)^{i,j}$ denotes the $(i,j)$-entry of the inverse matrix
$g(x)^{-1}$, the projection of a column vector $V\in\RR^n$ onto $\TT_x(\partial E)$ is defined by
$$
\mbox{\rm proj}_{\TT_x(\partial E)}(V)=\sum_{1\leq i,j< n}~g(x)^{i,j}~\left(T_j(x)^{\prime}V\right) ~T_i(x).
$$
\subsection{The shape matrix}\label{shape-sec}
 Consider the Monge parametrization 
\begin{equation}\label{monge-p}
\psi~:~\theta=(\theta_i)_{1\leq i<n}\in \Sa:=\RR^{n-1}\mapsto  \psi(\theta)=\left(
\begin{array}{c}
\theta\\
\varphi(\theta)
\end{array}
\right)\in \partial E\subset\RR^n.
\end{equation}
In this chart, the tangent vectors and the normal unit vector at  $x=\psi(\theta)$ are given for any $1\leq i<n$ by
$$
T_i^{\psi}(\theta):=\partial_{\theta_i}\psi(\theta)=T_i\left(\psi(\theta)\right)\in \TT_{x}(\partial E)\quad \mbox{\rm and}\quad  N^{\psi}(\theta):=N(\psi(\theta))\in \TT^{\perp}_{x}(\partial E).
$$
For any $1\leq i,j<n$  we have
$$
\begin{array}{l}
\left(\partial_{\theta_i}\psi(\theta)\right)^{\prime} N^{\psi}(\theta)=0\\
\\
\Longrightarrow
\Omega(\psi(\theta))_{i,j}:=
\left(\partial_{\theta_i,\theta_j}\psi(\theta)\right)^{\prime} N^{\psi}(\theta)=-\left(\partial_{\theta_i}\psi(\theta)\right)^{\prime} \partial_{\theta_j}N^{\psi}(\theta).
\end{array}$$
Observe that for  $x=\psi(\theta)$, 
$$
\partial_{\theta_i} N^{\psi}(\theta)=\sum_{1\leq k\leq n}\left(\partial_{x_k}N\right)
(x)~\partial_{\theta_i}\psi^k(\theta)=\left(\nabla N(x)\right)^{\prime}~\partial_{\theta_i}\psi(\theta)
$$
from which we check that for any $1\leq i,j<n$ the coefficients of the second fundamental form can  be computed as follows:
$$
\Omega(x)_{i,j}=-
\left(\partial_{\theta_i}\psi(\theta)\right)^{\prime} \left(\nabla N(x)\right)^{\prime}\partial_{\theta_i}\psi(\theta).
$$
We set
$$
\left(\partial N^{\psi}(\theta)\right)^{\prime}:=\left(\partial_{\theta_1} N^{\psi}(\theta),\ldots,\partial_{\theta_{n-1}} N^{\psi}(\theta)\right)\in \left(\TT_{\psi(\theta)}(\partial E)\right)^{n-1}.$$
In this notation, for any $x=\psi(\theta)$ we have the matrix formulation 
\begin{eqnarray*}
\Omega(x)&:=&-\partial\psi(\theta)~\left(\partial N^{\psi}(\theta)\right)^{\prime}=\left(\left(\partial_{\theta_i,\theta_j}\psi(\theta)\right)^{\prime} N(x)\right)_{1\leq i,j<n}\\
&=&-\frac{\nabla^2 \varphi(\theta)}{\sqrt{1+\Vert\nabla \varphi(\theta)\Vert^2}}\quad\mbox{with}\quad
\nabla^2 \varphi(\theta):=\left(\partial_{\theta_i,\theta_j}\varphi(\theta)\right)_{1\leq i,j<n}.
\end{eqnarray*}
We also readily check the matrix formulation of the Weingarten's equations
\begin{eqnarray*}
\left(\partial N^{\psi}(\theta)\right)^{\prime}
&=&\left(\left(\partial \psi(\theta)\right)^{\prime}g(\psi(\theta))^{-1}\right)\left(
\partial \psi(\theta)\right)\left(\partial N^{\psi}(\theta)\right)^{\prime}=-\left(\partial \psi(\theta)\right)^{\prime}\Wa \left(x\right).
\end{eqnarray*}
In the above display, $\Wa(x)$ stands for  the shape matrix  (a.k.a. the Weingarten map or the mixed second fundamental form) defined by
\begin{eqnarray*}
\Wa(x)&:=&g(x)^{-1}\Omega(x)\\
&=&-\frac{1}{\sqrt{1+\Vert\nabla \varphi(x_{-n})\Vert^2}}~\left(I+\nabla\varphi(x_{-n})\nabla\varphi(x_{-n})^{\prime}
\right)^{-1}\nabla^2 \varphi(x_{-n}).
\end{eqnarray*}
We summarize the above discussion in  the following proposition.
\begin{prop}\label{prop-weingarten}
For any $1\leq i<n$  we have the Weingarten's equations
$$
\partial_{\theta_i} N^{\psi}(\theta)=-\sum_{1\leq k<n}\Wa \left(\psi(\theta)\right)
_{k,i}~\partial_{\theta_k}\psi(\theta)\in \TT_{\psi(\theta)}(\partial E).$$
\end{prop}

\begin{examp}\label{examp-2d}
For $n=2$ we have $x\in\RR\mapsto \overline{\varphi}(x)=\varphi(x)-x$, so that the boundary $\partial E=\overline{\varphi}^{-1}(\{0\})$ coincides with the graph of the function $\varphi$. In this context, the metric and Weingarten map at $x\in \partial E=\{x=(x_1,x_2)\in \RR^2~:~x_2=\varphi(x_1)\}$ take the form
$$
g(x)=1+\Vert\partial\varphi(x_1)\Vert^2\quad\mbox{and}\quad
\Wa(x)=-\frac{1}{\left(1+\Vert\partial\varphi(x_1)\Vert^2\right)^{3/2}}~\partial^2 \varphi(
x_{1}).
$$

\end{examp}

\begin{examp}\label{examp-2d-3d}
For $n=3$,  the boundary $\partial E$ is given by the surface in $\RR^3$ defined
$$
\partial E:=\{x=(x_i)_{1\leq i\leq 3}\in\RR^3~:~x_3=\varphi(x_1,x_2)\}.
$$ 
The Monge parametrization is given by
$$
\psi~:~\theta=(\theta_1,\theta_2)\in \RR^{2}\mapsto  \psi(\theta)=\left(
\begin{array}{c}
\theta_1\\
\theta_2\\
\varphi(\theta_1,\theta_2)
\end{array}
\right)\in \partial E\subset\RR^3.
$$
In this situation, the tangent vectors at $x\in \partial E$ are given by
$$
T_1(x)=\left(
\begin{array}{c}
1\\
0\\
\partial_{x_1}\varphi(x)
\end{array}
\right)\quad \mbox{and}\quad
T_2(x)=\left(
\begin{array}{c}
0\\
1\\
\partial_{x_2}\varphi(x)
\end{array}
\right).
$$
In the same vein, whenever $E=\{x\in\RR^3~:~ \varphi(x_1,x_2)\leq x_3\}
$ the outward pointing unit normal at $x\in \partial E$ is given by
$$
N(x)=\frac{1}{\sqrt{1+(\partial_{x_1}\varphi(x))^2+(\partial_{x_2}\varphi(x))^2}}\left(
\begin{array}{c}
\partial_{x_1}\varphi(x)\\
\partial_{x_2}\varphi(x)\\
-1
\end{array}
\right).
$$
The inner product $g(x)$ is easily computed and given by
$$
g(x)=\left(
\begin{array}{cc}
1+(\partial_{x_1}\varphi(x))^2&(\partial_{x_1}\varphi(x))(\partial_{x_2}\varphi(x))\\
(\partial_{x_1}\varphi(x))(\partial_{x_2}\varphi(x))
&1+(\partial_{x_2}\varphi(x))^2
\end{array}
\right).
$$
The inverse metric is given by
\begin{gather*}
g(x)^{-1}=\frac{1}{\mbox{\rm det}(g(x))}~\left(
\begin{array}{cc}
1+(\partial_{x_2}\varphi(x))^2&-(\partial_{x_1}\varphi(x))(\partial_{x_2}\varphi(x))\\
-(\partial_{x_1}\varphi(x))(\partial_{x_2}\varphi(x))
&1+(\partial_{x_1}\varphi(x))^2
\end{array}
\right)
\end{gather*}
with
$$
\mbox{\rm det}(g(x))=1+(\partial_{x_1}\varphi(x))^2+(\partial_{x_2}\varphi(x))^2=1+\Vert\nabla \varphi(x)\Vert^2.
$$
The second fundamental form is also given by
\begin{gather*}
\Omega(x)=-\frac{1}{\sqrt{1+\Vert\nabla \varphi(x)\Vert^2}}~\left(
\begin{array}{cc}
\partial^2_{x_1}\varphi(x)&\partial_{x_1,x_2}\varphi(x)\\
\partial_{x_1,x_2}\varphi(x)
&\partial^2_{x_2}\varphi(x)
\end{array}
\right)
\end{gather*}
and the Weingarten map is defined by
\begin{gather*}
\begin{array}{l}
\displaystyle\Wa(x)=-\frac{1}{(1+\Vert\nabla \varphi(x)\Vert^2)^{3/2}}~\\
\\
\times\left(
\begin{array}{cc}
(1+(\partial_{x_2}\varphi(x))^2)\partial^2_{x_1}\varphi(x)-(\partial_{x_1}\varphi(x))(\partial_{x_2}\varphi(x))\partial_{x_1,x_2}\varphi(x)&(1+(\partial_{x_2}\varphi(x))^2)\partial_{x_1,x_2}\varphi(x)--(\partial_{x_1}\varphi(x))(\partial_{x_2}\varphi(x))\partial^2_{x_2}\varphi(x)\\
&\\
-(\partial_{x_1}\varphi(x))(\partial_{x_2}\varphi(x))\partial^2_{x_1}\varphi(x)+(1+(\partial_{x_1}\varphi(x))^2)\partial_{x_1,x_2}\varphi(x)&-(\partial_{x_1}\varphi(x))(\partial_{x_2}\varphi(x))\partial_{x_1,x_2}\varphi(x)+(1+(\partial_{x_1}\varphi(x))^2)\partial^2_{x_2}\varphi(x)
\end{array}
\right).
\end{array}
\end{gather*}
\end{examp}

\subsection{Surface and volume forms}
The surface form $\sigma_{\partial}$ on the boundary $\partial E$ expressed in the chart $\psi$ introduced in (\ref{monge-p}) is given by
$$
\left(\sigma_{\partial}\circ \psi^{-1}\right)(d\theta)=\sqrt{\mbox{\rm det}\left(g(\psi(\theta))\right)}~d\theta$$
 with the Gramian of the coordinate chart \begin{eqnarray*}
g(\psi(\theta))&:=&\mbox{\rm Gram}\left(\partial_{\theta_1}\psi(\theta),\ldots,\partial_{\theta_{n-1}}\psi(\theta)\right)\\
&:=&\left(\partial \psi(\theta)\right)\left(\partial \psi(\theta)\right)^{\prime}=I+\nabla\varphi(\theta)(\nabla\varphi(\theta))^{\prime}
\end{eqnarray*}
with the coordinates tangent vectors
$
\partial\psi(\theta)=T^{\psi}(\theta):=T(\psi(\theta))$. To check this claim recall that the surface area spaced by the column vectors
$$
\partial \psi(\theta)^{\prime}:=\left(
\partial_{\theta_1}\psi(\theta),\ldots,\partial_{\theta_{n-1}}\psi(\theta)\right)
$$
is equal to the volume of the parallelepided generated by the column vectors
$$
\left(\partial \psi(\theta)^{\prime},N(\psi(\theta))\right):=\left(
\partial_{\theta_1}\psi(\theta),\ldots,\partial_{\theta_{n-1}}\psi(\theta),N(\psi(\theta))\right)
$$
which is given by the determinant of the column vectors, so that
$$
\left(\sigma_{\partial}\circ \psi^{-1}\right)(d\theta)=\vert\mbox{\rm det}\left(\partial \psi(\theta)^{\prime},N(\psi(\theta))\right)\vert~d\theta.
$$
On the other hand, we have
\begin{gather*}
\begin{array}{l}
\left(
\begin{array}{c}
\partial \psi(\theta)\\
N(\psi(\theta))^{\prime}
\end{array}
\right)\left(
\begin{array}{c}
\partial \psi(\theta)^{\prime},
N(\psi(\theta))
\end{array}
\right)\\
\\
=
\left(
\begin{array}{c}
(\partial_{\theta_1}\psi(\theta))^{\prime}\\
\vdots\\
(\partial_{\theta_1}\psi(\theta))^{\prime}\\
N(\psi(\theta))^{\prime}
\end{array}
\right)\left(
\begin{array}{c}
\partial_{\theta_1}\psi(\theta),
\ldots,
\partial_{\theta_1}\psi(\theta),
N(\psi(\theta))
\end{array}
\right)=\left(\begin{array}{cc}
\left(\partial \psi(\theta)\right)\left(\partial \psi(\theta)\right)^{\prime}&0_{n-1,1}\\
0_{1,n-1}&1
\end{array}\right).
\end{array}
\end{gather*}
This implies that
\begin{eqnarray*}
\vert\mbox{\rm det}\left(\partial \psi(\theta)^{\prime},N(\psi(\theta))\right)\vert&=&\sqrt{\vert\mbox{\rm det}\left(\left(\partial \psi(\theta)^{\prime},N(\psi(\theta))\right)^{\prime}\left(\partial \psi(\theta)^{\prime},N(\psi(\theta))\right)\right)\vert}\\
&=&\sqrt{\mbox{\rm det}\left(\left(\partial \psi(\theta)\right)\left(\partial \psi(\theta)\right)^{\prime}\right)}.
\end{eqnarray*}

Using the determinant perturbation formula w.r.t. rank-one matrices $\mbox{\rm det}\left(I+uv^{\prime}\right)=1+v^{\prime}u$ which is valid for any column vectors $u,v\in\RR^n$ we check that $$
\mbox{\rm det}\left(I+\nabla\varphi(\theta)\nabla\varphi(\theta)^{\prime}
\right)=1+\Vert\nabla \varphi(\theta)\Vert^2.
$$
This yields the formula
$$
\left(\sigma_{\partial}\circ \psi^{-1}\right)(d\theta)=\sqrt{1+\Vert\nabla \varphi(\theta)\Vert^2}~d\theta.$$

The mapping $F$ defined in (\ref{Fermi-coordinates}) can also be rewritten as a chart $\overline{\psi}$ on $\Da_{r}(E)$ defined for any $(\theta,u)\in \left(\Sa\times [0,r]\right)$ defined  by
$$
\overline{\psi}(\theta,u):=
F(\psi(\theta),u)=\psi(\theta)+ u~N(\psi(\theta))\in \Da_{r}(E).
$$
The Jacobian matrix of $\overline{\psi}$ is given by
$$
\mbox{\rm Jac}(\overline{\psi})(\theta,u)=\left(\partial_{\theta_1}\overline{\psi}(\theta,u),\ldots,\partial_{\theta_{n-1}}\overline{\psi}(\theta,u),N(\psi(\theta))\right).
$$
By Proposition~\ref{prop-weingarten} we have
\begin{eqnarray*}
\partial_{\theta_i}\overline{\psi}(\theta,u)&=&\partial_{\theta_i}\psi(\theta)+u~\partial_{\theta_i}N^{\psi}(\theta)\\
&=&\partial_{\theta_i}\psi(\theta)-u~\sum_{1\leq k<n}\partial_{\theta_k}\psi(\theta)~\Wa \left(\psi(\theta)\right)
_{k,i}.
\end{eqnarray*}
This yields the formula
$$
\left(\partial_{\theta_1}\overline{\psi}(\theta,u),\ldots,\partial_{\theta_{n-1}}\overline{\psi}(\theta,u)\right)
=\left(\partial_{\theta_1}\psi(\theta,u),\ldots,\partial_{\theta_{n-1}}\psi(\theta)\right)~\left(I-u~\Wa(\psi(\theta))\right)
$$
from which we check that
$$
\vert \mbox{\rm det}\left(\mbox{\rm Jac}(\overline{\psi})(\theta,u)\right)\vert =\sqrt{\mbox{\rm det}(g(\psi(\theta))}~\left\vert\mbox{\rm det}\left(I-u~\Wa(\psi(\theta))\right)\right\vert.
$$
Note that $\overline{\psi}(\theta,0)=\psi(\theta)$, and for any given $u<r$, the mapping $\theta\mapsto\overline{\psi}(\theta,u)$ is a chart on $\partial E_u$. 
This yields the following proposition. For the convenience of the reader, a more detailed proof of the next proposition is provided in the appendix on page~\pageref{prop-vol-form-weignarten-proof}. 
\begin{prop}\label{prop-vol-form-weignarten}
For any $u\leq r$, the surface form  $\sigma_{\partial,u}$ on the boundary $\partial E_{u}$ expressed in the chart  $\theta\in \Sa\mapsto \overline{\psi}(\theta,u):=F(\psi(\theta),u)$ is given by the formula 
$$
\left(\sigma_{\partial,u}\circ \overline{\psi}(\point,u)^{-1}\right)(d\theta)=
\left\vert\mbox{\rm det}\left(I-u~\Wa(\psi(\theta))\right)\right\vert ~\left(\sigma_{\partial}\circ \psi^{-1}\right)(d\theta)
$$
with
$$
\begin{array}{l}
\left\vert\mbox{\rm det}\left(I-u~\Wa(\psi(\theta))\right)\right\vert
\\
\\
\displaystyle=\vert\mbox{\rm det}\left(I+\frac{u}{\sqrt{1+\Vert \nabla \varphi(\theta)\Vert^2}}~\left(I+\nabla\varphi(\theta)\nabla\varphi(\theta)^{\prime}
\right)^{-1}\nabla^2\varphi(\theta)\right)\vert.
\end{array}$$
 In addition, the volume form $\sigma_{\Da_r(E)}$ on $\Da_{r}(E)$ expressed in the chart 
 $\overline{\psi}$ is given by
$$
\left(\sigma_{\Da_{r}(E)}\circ \overline{\psi}^{-1}\right)(d(\theta,u))=\left\vert\mbox{\rm det}\left(I-u~\Wa(\psi(\theta))\right)\right\vert~ \left(\sigma_{\partial}\circ \psi^{-1}\right)(d\theta)~du.
$$

\end{prop}

Using Jacobi's formula for the derivative of determinants, we also have
$$
\displaystyle\partial_u\log{\mbox{\rm det}\left(I-u~\Wa(x)\right)}=-\tr\left(\left(I-u~\Wa(x)\right)^{-1}\Wa(x)\right).
$$

The level-set density function defined in (\ref{level-set-density})
 expressed in the chart $\psi$ is given by the formula
$$
\begin{array}{l}
\displaystyle q^{\partial}_t(x,r)\\
\\
\displaystyle =\int_{\Sa}~q_t\left(x,\psi(\theta)+ r~N(\psi(\theta))\right)~\left\vert\mbox{\rm det}\left(I-r~\Wa(\psi(\theta))\right)\right\vert~\sqrt{\mbox{\rm det}(g(\psi(\theta)))}~d\theta.
\end{array}$$
\subsection{Boundary decompositions}
 For some given coordinate index $k\in \{1,\ldots,n\}$ and $x=(x_i)_{1\leq i\leq n}\in \RR^{n}$ we set
 $$
 x_{-k}:=(x_i)_{i\in \Ia}\quad \mbox{\rm with}\quad
 \Ia:=\{1,\ldots,n\}-\{k\}
$$
We further assume that $$\partial E=\left\{x\in \RR^{n}~:~x_{-k}\in  \Sa\quad\mbox{\rm and}\quad\varphi(x_{-k})=x_{k}\right\}=\overline{\varphi}^{-1}(\{0\})$$ is defined as the null level set  of some global defining function of the form
$$
 \overline{\varphi}~:~x\in \{(x_i)_{1\leq i\leq n}\in \RR^{n}~:~x_{-k}\in \Sa\}\mapsto
 \overline{\varphi}(x):=\varphi(x_{-k})-x_k\in\RR
 $$
 for some open domain $\Sa\subset\RR^{n-1}$. 
 \begin{examp}[Cylindrical boundaries]\label{ref-ex-cylindrical-boundary}
Let $1\leq k\leq n_1$ and $n=n_1+n_2$ for some $n_1> 1$ and $n_2\geq 1$.
Consider a  domain $\Sa$ of the form $\Sa=(\widehat{\Sa}\times\RR^{n_2})$ with $\widehat{\Sa}\subset\RR^{n_1-1}$ and assume that
$$
\forall y\in \RR^{n_1}~s.t.~y_{-k}\in \widehat{\Sa} \quad\mbox{and}\quad \forall z\in \RR^{n_2}\quad\mbox{we have}\quad\varphi(y_{-k},z):=\widehat{\varphi}(y_{-k}).
$$ 
In this situation, the set $\partial E$ is a cylindrical boundary  given by the formula
$$
\partial E=\partial\widehat{E}\times\RR^{n_2}\quad \mbox{with}
\quad \partial\widehat{E}:=\left\{y\in \RR^{n_1}~:~y_{-k}\in  \widehat{\Sa}\quad\mbox{and}\quad \widehat{\varphi}(y_{-k})=y_{k}\right\}.
$$
\end{examp}
In this context, the coordinates of the outward normal  by
$$
N^j(x)=\frac{\epsilon}{\sqrt{1+\Vert \nabla\varphi(x_{-k})\Vert^2}}~\left(1_{\Ia}(j)~\partial_{x_{j}}\varphi(x_{-k})+1_{k}(j)~(-1)\right)$$
with the orientation parameter $\epsilon=1$ when $E=\overline{\varphi}^{-1}(]-\infty,0[)$; and $\epsilon=-1$ when  
$E=\overline{\varphi}^{-1}(]0,+\infty[)$.
In the same vein, the entries $T^{j}_i(x)$ of the tangent vectors $T_i(x)$ indexed by $i\in \Ia$ are given for any $1\leq j\leq n$ by
$$
T^{j}_i(x)=1_{i}(j)~+1_{k}(j)~\partial_{x_i}\varphi(x_{-k}).
$$
Consider the $(n\times (n-1))$-matrix
\begin{gather*}
 T(x)^{\prime}:=\left(T_1(x),\ldots,T_{k-1}(x),T_{k+1}(x),\ldots,T_{n}(x)\right).
\end{gather*}
In this notation, the inner product $g(x)$ on the tangent space $\TT_x(\partial E)$ is given by the $((n-1)\times (n-1))$-square Gramian matrix
\begin{eqnarray*}
g(x)&=&T(x)T(x)^{\prime}=I+\nabla\varphi(x_{-k})\nabla\varphi(x_{-k})^{\prime}
\end{eqnarray*}
with the gradient column vector
\begin{gather*}
\nabla\varphi(x_{-k}):=\left(\partial_{x_i}\varphi\left(x_{-k}\right)\right)_{i\in \Ia}=\left(
\begin{array}{c}
\partial_{x_1}\varphi\left(x_{-k}\right)\\
\vdots\\
\partial_{x_{k-1}}\varphi\left(x_{-k}\right)\\
\partial_{x_{k+1}}\varphi\left(x_{-k}\right)\\
\vdots\\
\partial_{x_{n}}\varphi\left(x_{-k}\right)
\end{array}
\right)\in \RR^{n-1}.
\end{gather*}
We check this claim using the fact that for any $i_1,i_2\in \Ia$ we have
\begin{eqnarray*}
T_{i_1}(x)^{\prime}T_{i_2}(x)&=&\sum_{1\leq j\leq n}\left(1_{i_1}(j)~+1_{k}(j)~\partial_{x_{i_1}}\varphi(x_{-k})\right)
\left(1_{i_2}(j)~+1_{k}(j)~\partial_{x_{i_2}}\varphi(x_{-k})\right)\\
&=&1_{i_1=i_2}+\partial_{x_{i_1}}\varphi(x_{-k})~\partial_{x_{i_2}}\varphi(x_{-k}).
\end{eqnarray*}

The parametrization of the hyper surface $\partial E$ is now given by the chart function
$$
\psi~:~\theta=(\theta_i)_{i\in \Ia}\in \Sa
\mapsto
\psi(\theta)\in \partial E$$
with 
$$  \forall 1\leq j\leq n\qquad\psi^{j}(\theta):=1_{\Ia}(j)~\theta_j~+1_{k}(j)~\varphi(\theta).
$$
For any $1\leq j\leq n$ and $ i_1,i_2\in \Ia$ observe that
 $$
 \partial_{\theta_{i_1}}\psi(\theta)=T^{\psi}_{i_1}(\theta):=T_{i_1}(\psi(\theta))\quad\mbox{\rm and}
 \quad
  \partial_{\theta_{i_1,\theta_{i_2}}}\psi^j(\theta)=1_{k_{\iota}}(j)~  \partial_{\theta_{i_1},\theta_{i_2}}\varphi(\theta).
  $$
  This implies that
  $$
\left(  \partial_{\theta_{i_1,\theta_{i_2}}}\psi(\theta)\right)^{\prime}  N(\psi(\theta))=-\epsilon~\frac{ \left(\nabla^2\varphi(\theta)\right)_{i_1,i_2}}{\sqrt{1+\Vert \nabla\varphi(x_{-k})\Vert^2}}\quad\mbox{\rm with}\quad\nabla^2\varphi(\theta):=\left(\partial_{\theta_{i_1},\theta_{i_2}}\varphi(\theta)\right)_{(i_1,i_2)\in \Ia^2}.
  $$
We set $ (\partial\psi(\theta))^{\prime}:= T(\psi(\theta))^{\prime}$ and $N^{\psi}(\theta):=N(\psi(\theta))$.
In this notation,  we also have
$$
\begin{array}{l}
\displaystyle\left(\partial N^{\psi}(\theta)\right)^{\prime}\\
\\
\displaystyle:=\left(\partial_{\theta_1} N^{\psi}(\theta),\ldots,\partial_{\theta_{k_{\iota}-1}} N^{\psi}(\theta),\partial_{\theta_{k_{\iota}+1}} N^{\psi}(\theta),\ldots,\partial_{\theta_{n}} N^{\psi}(\theta)\right)\\
\\
\displaystyle=-\Wa(\psi(\theta)):=-g(\psi(\theta))^{-1}\Omega(\psi(\theta))
\quad \mbox{\rm
with}\quad
\Omega(\psi(\theta))=-\epsilon~\frac{\nabla^2 \varphi(\theta)}{\sqrt{1+\Vert\nabla \varphi(\theta)\Vert^2}}.
\end{array}$$

 \begin{examp}
 For the cylindrical  boundary discussed in Example~\ref{ref-ex-cylindrical-boundary}, the inner product and the Weingarten map on the boundary $\partial \widehat{E}$ are given for any
 $y\in \partial \widehat{E}$ by the matrices
 $$
  \widehat{g}(y)= I_{(n_1-1,n_1-1)}+\nabla\widehat{\varphi}(y_{-k})\nabla\widehat{\varphi}(y_{-k})^{\prime}\quad\mbox{and}
 \quad \widehat{\Wa}(y):=\epsilon~\widehat{g}(y)^{-1}\frac{\nabla^2 \widehat{\varphi}(y_{-k})}{\sqrt{1+\Vert\nabla \widehat{\varphi}(y_{-k})\Vert^2}}
 $$
 with the gradient column vector  and the Hessian matrix given by
\begin{eqnarray*}
 \nabla\widehat{\varphi}(y_{-k})&:=&\left(
\partial_{y_i}\widehat{\varphi}(y_{-k})
\right)_{i\in \widehat{\Ia}}\\
\nabla^2\widehat{\varphi}(y_{-k})&:=&\left(
\partial_{y_{i_1},y_{i-2}}\widehat{\varphi}(y_{-k})
 \right)_{i_1,i_2\in \widehat{\Ia}}\quad \mbox{with}\quad
 \widehat{\Ia}:=\{1,\ldots,n_1\}-\{k\}.
\end{eqnarray*}
Observe that
$$
\mbox{\rm det}(  \widehat{g}(y))=1+\Vert  \nabla\widehat{\varphi}(y_{-k})\Vert^2
\quad \mbox{and}\quad\nabla^2\varphi(y_{-k},z)=\left(
 \begin{array}{ll}
\nabla^2\widehat{\varphi}(y_{-k})&0_{(n_1-1,n_2)}\\
0_{(n_2,n_1-1)} &0_{(n_2,n_2)}
 \end{array}
 \right).
$$
In this case, the inner product and the Weingarten map on the boundary $\partial E$ are given for any
 point $x=(y,z)\in (\partial \widehat{E}\times\RR^{n_2})$  by the matrices
 \begin{gather*}
 g(x)
 =\left(
 \begin{array}{ll}
  \widehat{g}(y)&0_{(n_1-1,n_2)}\\
0_{(n_2,n_1-1)} &I_{(n_2,n_2)}
 \end{array}
 \right)\quad \mbox{and}\quad \Wa(x)
 =\left(
 \begin{array}{ll}
\widehat{\Wa}(y)&0_{(n_1-1,n_2)}\\
0_{(n_2,n_1-1)} &0_{(n_2,n_2)}
 \end{array}
 \right).
 \end{gather*}
Observe that the above matrices are bounded (w.r.t. any matrix norm) as soon as $\partial \widehat{E}$ is bounded.
 \end{examp}

More generally, assume that  the boundary $\partial E\subset\cup_{\iota\in \Ja}\Oa(\iota)\subset\RR^n$ admits a finite covering by open connected subsets $\Oa(\iota)\subset\RR^n$ indexed by some finite set $\Ja$. In addition, there exists some local defining smooth functions $\overline{\varphi}_{\iota}$ with non vanishing gradients on $\Oa(\iota)$ such that
$$
\partial E(\iota):=\partial E\cap \Oa(\iota)=\overline{\varphi}_{\iota}^{-1}\left(\{0\}\right)
\quad\mbox{\rm
and}\quad
E(\iota):=E\cap \Oa(\iota)=\overline{\varphi}_{\iota}^{-1}\left(]0,\infty[\right).
$$
 Up to  
shrinking the set $\Oa(\iota)$, by the implicit function theorem there is no loss of generality to assume that the defining functions are given by
$$
\overline{\varphi}_{\iota}~:~x=(x_i)_{1\leq i\leq n}\in \Oa(\iota)\mapsto
\overline{\varphi}_{\iota}(x)=\varphi_{\iota}(x_{-k_{\iota}})-x_{k_{\iota}}
$$
for some parameter $1\leq k_{\iota}\leq n$ and some smooth function $\varphi_{\iota}$ on some ball $\Sa(\iota)\subset \RR^{n-1}$. We set $\Ia_{\iota}:=\{1,\ldots,n\}-\{k_{\iota}\}$. In this notation, the parametrization of the hyper surface $\partial E(\iota)$ is now given by the smooth homeomorphism
\begin{equation}\label{def-atlas}
\psi_{\iota}~:~\theta=(\theta_i)_{i\in \Ia_{\iota}}\in \Sa(\iota)
\mapsto
\psi_{\iota}(\theta)\in \partial E(\iota)\quad \mbox{\rm with}\quad  \psi^{j}_{\iota}(\theta):=1_{\Ia_{\iota}}(j)~\theta_j~+1_{k_{\iota}}(j)~\varphi_{\iota}(\theta).
\end{equation}
The first and second fundamental forms on $\TT_x\left(\partial E(\iota)\right)$ as well as the Weingarten map at $x\in \partial E(\iota)$ are given by
\begin{eqnarray*}
g_{\iota}(x)&=&I+\nabla\varphi_{\iota}(x_{-k_{\iota}})\nabla\varphi(x_{-k_{\iota}})^{\prime}\\\
\Omega_{\iota}(x)&=&-\frac{\nabla^2 \varphi_{\iota}(x_{-k_{\iota}})}{\sqrt{1+\Vert\nabla \varphi_{\iota}(x_{-k_{\iota}})\Vert^2}}
\quad \mbox{\rm and}\quad \Wa_{\iota}(x):=g_{\iota}(x)^{-1}\Omega_{\iota}(x).
\end{eqnarray*}
The atlas $\Aa=(\psi_{\iota},\Sa_{\iota})_{\iota\in \Ja}$ represents a collection of local  coordinate systems of  the boundary $\partial E=\cup_{\iota\in \Ja}\partial E(\iota)$.
In this situation, the surface form on $\partial E$ and  the volume form $\sigma_{\Da_{r}(E)}$ on $\Da_{r}(E)$  expressed in the atlas $\Aa$ are
defined by the formulae
\begin{eqnarray*}
\sigma^{\Aa}_{\partial}(d\theta)&:=&\sum_{\iota\in \Ja}~\pi_{\iota}\left(\psi_{\iota}(\theta)\right)~1_{\Sa(\iota)}(\theta)~
\sqrt{1+\Vert\nabla \varphi_{\iota}(\theta)\Vert^2}~d\theta
\\
\sigma_{\Da_{r}(E)}^{\Aa}(d(\theta,u))&:=&
\sum_{\iota\in \Ja}~\pi_{\iota}\left(\psi_{\iota}(\theta)\right)~1_{\Sa(\iota)}(\theta)\left\vert\mbox{\rm det}\left(I-u~\Wa_{\iota}(\psi_{\iota}(\theta))\right)\right\vert~ \sqrt{1+\Vert\nabla \varphi_{\iota}(\theta)\Vert^2}~du.
\end{eqnarray*}
In the above display, $\pi_{\iota}~:~\partial E\mapsto [0,1]$ stands for some partition of unity subordinate to the open cover of the boundary induced by the atlas.

\begin{examp}
Observe that the metric in the graph model discussed in Example~\ref{examp-2d} is not necessarily bounded.
In this context, we can also use for any $a<a_+<b_-<b$ a covering of the form  
$$\Oa(0)=]a,b[\times\RR\qquad \Oa(-1)=]b_-,+\infty[\times \RR\quad \mbox{\rm and}\quad
\Oa(1)=]-\infty,a_+[\times\RR.$$ 
For instance when $\varphi(z)=z^2$ and $(a,a_+,b_-,b)=(-2,-1,1,2)$ we have
\begin{eqnarray*}
\partial E(0)
&=&
\{(x_1,x_2)\in ]-2,2[\times ]4,\infty[~:~x_2=\varphi_0(x_1)\}\\
\partial E(1)&=&
\{(x_1,x_2)\in ]-\infty,-1[\times ]1,+\infty[~:~x_1=\varphi_{1}(x_2)\}\\
\partial E(-1)&=&\{(x_1,x_2)\in ]1,\infty[\times  ]1,+\infty[~:~x_1=\varphi_{-1}(x_2)\}
\end{eqnarray*}
with the functions
$$
\varphi_0(z)=z^2\quad\mbox{and}\quad \forall \epsilon\in\{-1,1\}\qquad \varphi_{\epsilon}(z)=-\epsilon \sqrt{z}.
$$
Whenever $E$ is the sub-graph of $\varphi$, the parameter $\epsilon\in\{-1,1\}$ plays the role of the orientation and the outward pointing unit normal 
vector at $x\in \partial E(0)$ and $y\in \partial E(\epsilon)$ are given by
$$
N_0(x)=\frac{1}{\sqrt{1+4x_1^2}}\left(
\begin{array}{c}
2x_1\\
-1
\end{array}
\right)\quad\mbox{and}\quad
N_{\epsilon}(y)=\frac{\epsilon}{\sqrt{1+1/(4y_2)}}\left(
\begin{array}{c}
-1\\
-\epsilon/\sqrt{4y_2}
\end{array}
\right).
$$
The tangent vectors at $x\in \partial E(0)$ and at $y\in \partial E(\epsilon)$ are defined by
$$
T_{0}(x)=\left(
\begin{array}{c}
1\\
2x_1
\end{array}
\right)\quad \mbox{\rm and}\quad
T_{\epsilon}(y)=\left(
\begin{array}{c}
-\epsilon/(\sqrt{4y_2})\\
1
\end{array}
\right).
$$
The above sub-graphs can be described with $3$ charts $\{\psi_0,\psi_{+1},\psi_{-1}\}$ defined for  any $\epsilon\in \{-1,1\}$ by
$$
\psi_0~:~\theta\in ]-2,2[\mapsto \psi_0(\theta)=\left(
\begin{array}{c}
\theta\\
\theta^2
\end{array}
\right)
\quad
\mbox{and}
\quad
\psi_{\epsilon}~:~\theta\in ]1,\infty[\mapsto \psi_{\epsilon}(\theta):=\left(
\begin{array}{c}
-\epsilon\sqrt{\theta}\\
\theta
\end{array}
\right).$$
In this situation, the tangent vectors are given by
$$
\partial_{\theta}\psi_0(\theta)=T_0(\psi_0(\theta))=\left(
\begin{array}{c}
1\\
2\theta
\end{array}
\right)\quad \mbox{and}\quad
\partial_{\theta}\psi_{\epsilon}(\theta)=T_{\epsilon}\left(\psi_{\epsilon}(\theta)\right)=\left(
\begin{array}{c}
-\epsilon/(\sqrt{4\theta})\\
1
\end{array}
\right).
$$
In this context, for any $\theta\in ]-2,2[$ we have
$$
g(\psi_0(\theta))=1+4\theta^2\quad \mbox{and}\quad \Wa(\psi_0(\theta))=-2~\left(1+4\theta^2\right)^{-3/2}.
$$
In addition, 
for any $\theta\in ]1,\infty[$ we have
$$
g(\psi_{\epsilon}(\theta))=1+1/(4\theta)\quad \mbox{and}\quad 
\Wa(\psi_{\epsilon}(\theta))=-2\epsilon~\left(1+4\theta\right)^{-3/2}.
$$
Observe that the metric expressed in the chart  $\{\psi_0,\psi_{+1},\psi_{-1}\}$  is 
defined in terms of bounded functions.

\end{examp}
\begin{examp}
Consider the hyperbolic paraboloid boundary 
\begin{eqnarray*}
\partial E&=&\{(y_1,y_2,y_3)\in\RR^3~:~y_3=y_1^2+y_2^2\}\\
&=&\partial E(0)\cup
\partial E(1,1)\cup \partial E(1,-1)\cup \partial E(2,1)\cup \partial E(2,-1).
\end{eqnarray*}
In the above display, $\partial E(0)$ and $\partial E(i,\epsilon)$ with $i\in \{1,2\}$ and $\epsilon\in \{-1,1\}$ stands for the partition defined for any $\epsilon\in\{-1,1\}$ by
\begin{eqnarray*}
\partial E(0)&:=&\{y\in \RR^3~:~(y_1,y_2)\in\Sa_0\quad y_3=\varphi_0(y_1,y_2):=y_1^2+y_2^2 \}\\
\partial E(1,\epsilon)&:=&
\{y\in \RR^3~:~(y_1,y_3)\in\Sa\quad   y_2=\varphi_{1,\epsilon}(y_1,y_3):=\epsilon\sqrt{y_3-y_1^2} \}\\
\partial E(2,\epsilon)&:=&\{y\in \RR^3~:~(y_2,y_3)\in\Sa\quad  y_1=\varphi_{2,\epsilon}(y_1,y_2):=\epsilon\sqrt{y_3-y_2^2} \}
\end{eqnarray*}
with the sets
\begin{eqnarray*}
\Sa_0&:=&\{(y_1,y_2)\in\RR^2~:~y_1^2+y_2^2<2\}\\
\Sa&:=& \{(y_2,y_3)\in \RR^2~:~y_3>1\quad\&\quad \vert y_2\vert <\sqrt{3y_3/4} \}.
\end{eqnarray*}
On the truncated boundary $\partial E(0)$ we use a single chart defined by
$$
\psi_0~:~\theta=(\theta_1,\theta_2)\in \Sa_0\mapsto
\psi_0(\theta)=\left(
\begin{array}{c}
\theta_1\\
\theta_2\\
\theta_1^2+\theta_2^2
\end{array}\right)\in\partial E(0).
$$
On $\partial E(1,\epsilon)$ we use the chart defined  by
$$
\begin{array}{l}
\displaystyle\psi_{1,\epsilon}~:~\theta=(\theta_1,\theta_3)\in \Sa\mapsto
\psi_{1,\epsilon}(\theta)=\left(
\begin{array}{c}
\theta_1\\
\epsilon\sqrt{\theta_3-\theta_1^2}\\
\theta_3
\end{array}\right)\in \partial E(1,\epsilon).
\end{array}$$
Finally,  on $\partial E(2)$ we use the chart defined  by
$$
\begin{array}{l}
\displaystyle\psi_{2,\epsilon}~:~\theta=(\theta_2,\theta_3)\in \Sa\mapsto
\psi_{2,\epsilon}(\theta)=\left(
\begin{array}{c}
\epsilon\sqrt{\theta_3-\theta_2^2}\\
\theta_2\\
\theta_3
\end{array}\right)\in \partial E(2,\epsilon).\end{array}$$
For any $\theta=(\theta_1,\theta_2)\in \Sa_0$ we have
\begin{gather*}
\partial_{\theta_1}\psi_{0}(\theta)=\left(
\begin{array}{c}
1\\
0\\
2\theta_1
\end{array}
\right)\quad\mbox{and}\quad 
\partial_{\theta_2}\psi_{0}(\theta)=\left(
\begin{array}{c}
0\\
1\\
2\theta_2
\end{array}
\right).
\end{gather*}
In this chart, the metric is given by
\begin{gather*}
g(\psi_0(\theta))=\left(
\begin{array}{cc}
1+4\theta_1^2&4\theta_1\theta_2\\
4\theta_1\theta_2
&1+4\theta_2^2
\end{array}
\right)\quad \mbox{\rm and}\quad g(\psi_0(\theta))^{-1}=\frac{1}{1+4(\theta_1^2+\theta_2^2)}~\left(
\begin{array}{cc}
1+4\theta_2^2&-4\theta_1\theta_2\\
-4\theta_1\theta_2
&1+4\theta_1^2
\end{array}
\right).
\end{gather*}
In addition, the outward pointing unit normal at $\psi_0(\theta)\in \partial E(0)$ is given by
\begin{gather*}
N_0\left(\psi_0(\theta)\right)=\frac{1}{\sqrt{1+4(\theta_1^2+\theta_2^2)}}\left(
\begin{array}{c}
2\theta_1\\
2\theta_2\\
-1
\end{array}
\right)\quad\mbox{and}\quad\Omega_0\left(\psi_0(\theta)\right)=
\frac{1}{\sqrt{1+4(\theta_1^2+\theta_2^2)}}\left(
\begin{array}{cc}
-2&0\\
0&-2
\end{array}
\right).
\end{gather*}
For any $\theta=(\theta_1,\theta_3)\in \Sa$ we have
\begin{gather*}
\partial_{\theta_1}\psi_{1,\epsilon}(\theta)=\left(
\begin{array}{c}
1\\
\frac{-\epsilon\theta_1}{\sqrt{\theta_3-\theta_1^2}}\\
0
\end{array}
\right)\quad\mbox{and}\quad
\partial_{\theta_3}\psi_{1,\epsilon}(\theta)=\left(
\begin{array}{c}
0\\
\frac{\epsilon}{2\sqrt{\theta_3-\theta_1^2}}\\
1
\end{array}
\right).\end{gather*}
In this chart, the metric is given by
\begin{gather*}
g(\psi_{1,\epsilon}(\theta))=\left(
\begin{array}{cc}
1+\frac{\theta_1^2}{\theta_3-\theta_1^2}&-\frac{\theta_1}{2(\theta_3-\theta_1^2)}\\
-\frac{\theta_1}{2(\theta_3-\theta_1^2)}
&1+\frac{1}{4(\theta_3-\theta_1^2)}
\end{array}
\right)\quad\mbox{and}\quad
g(\psi_{1,\epsilon}(\theta))^{-1}=\frac{1}{1+\frac{\theta_1^2}{\theta_3-\theta_1^2}+\frac{1}{4(\theta_3-\theta_1^2)}}\left(
\begin{array}{cc}
1+\frac{1}{4(\theta_3-\theta_1^2)}&\frac{\theta_1}{2(\theta_3-\theta_1^2)}\\
\frac{\theta_1}{2(\theta_3-\theta_1^2)}
&1+\frac{\theta_1^2}{\theta_3-\theta_1^2}
\end{array}
\right).
\end{gather*}
In addition, the outward pointing unit normal at $\psi_{1,\epsilon}(\theta)\in \partial E(1,\epsilon)$ is given by
\begin{gather*}
N_{1,\epsilon}\left(\psi_{1,\epsilon}(\theta)\right)=\frac{-\epsilon}{\sqrt{1+\frac{\theta_1^2}{\theta_3-\theta_1^2}+\frac{1}{4(\theta_3-\theta_1^2)}}}\left(
\begin{array}{c}
\frac{-\epsilon\theta_1}{\sqrt{\theta_3-\theta_1^2}}\\
-1\\
\frac{\epsilon}{2\sqrt{\theta_3-\theta_1^2}}
\end{array}
\right)\quad\mbox{and}\quad\Omega_{1,\epsilon}\left(\psi_{1,\epsilon}(\theta)\right)=
\frac{-\epsilon}{\sqrt{1+\frac{\theta_1^2}{\theta_3-\theta_1^2}+\frac{1}{4(\theta_3-\theta_1^2)}}}\left(
\begin{array}{cc}
\frac{\epsilon\theta_3}{(\theta_3-\theta_1^2)^{3/2}}&\frac{-\epsilon\theta_1}{(\theta_3-\theta_1^2)^{3/2}}\\
\frac{-\epsilon\theta_1}{(\theta_3-\theta_1^2)^{3/2}}&-\frac{\epsilon}{4(\theta_3-\theta_1^2)^{3/2}}
\end{array}
\right).
\end{gather*}
Finally, for any $\theta=(\theta_2,\theta_3)\in \Sa$ we have
\begin{gather*}
\partial_{\theta_2}\psi_{2,\epsilon}(\theta)=\left(
\begin{array}{c}
\frac{-\epsilon\theta_2}{\sqrt{\theta_3-\theta_2^2}}\\
1\\
0
\end{array}
\right)\quad\mbox{and}\quad
\partial_{\theta_3}\psi_{2,\epsilon}(\theta)=\left(
\begin{array}{c}
\frac{\epsilon }{2\sqrt{\theta_3-\theta_2^2}}\\
0\\
1
\end{array}
\right).
\end{gather*}
In this chart, the metric and the outward pointing unit normal at $\psi_{2,\epsilon}(\theta)\in \partial E(2,\epsilon)$ are given by
\begin{gather*}
g(\psi_{2,\epsilon}(\theta))=\left(
\begin{array}{cc}
1+\frac{\theta_2^2}{\theta_3-\theta_2^2}&-\frac{\theta_2}{2(\theta_3-\theta_2^2)}\\
-\frac{\theta_2}{2(\theta_3-\theta_2^2)}
&1+\frac{1}{4(\theta_3-\theta_2^2)}
\end{array}
\right)\quad\mbox{and}\quad
N_{2,\epsilon}\left(\psi_{2,\epsilon}(\theta)\right)=-\frac{\epsilon}{\sqrt{1+\frac{\theta_2^2}{\theta_3-\theta_2^2}+\frac{1}{4(\theta_3-\theta_2^2)}}}\left(
\begin{array}{c}
-1\\
\frac{-\epsilon\theta_2}{\sqrt{\theta_3-\theta_2^2}}\\
\frac{\epsilon}{2\sqrt{\theta_3-\theta_2^2}}
\end{array}
\right).
\end{gather*}
\end{examp}

  \section*{Appendix}

      \subsection*{Proof of (\ref{V-L-ex})}\label{V-L-ex-proof}
We have
$$
P_{s,s+t}(V)\leq V+\int_s^{s+t} \left(-aP_{s,u}(V)+c\right)~du=V+ct- a\int_0^{t} P_{s,s+u}(V)~du
$$
and
\begin{eqnarray*}
\int_0^{t} P_{s,s+u}(V)~du&=&tP_{s,s+t}(V)-\int_0^{t}~u P_{s,s+u}(L_{s+u}(V))~du\geq tP_{s,s+t}(V)-ct^2/2.
\end{eqnarray*}
Combining the above estimates, we readily check that
$$
P_{s,s+t}(V)\leq (1+at)^{-1}V+ct~\frac{1+at/2}{1+at}\leq (1+at)^{-1}V+ct.
$$
This ends the proof of (\ref{V-L-ex}).\cqfd
    \subsection*{Proof of Lemma~\ref{lem-211}}\label{lem-211-proof}
 We have
\begin{equation}\label{sub-ref-append}
1\leq P_{s,s+t}(V)\leq
V-\int_0^{t}\left(P_{s,s+u}(\varphi(V))-c\right)~du 
\end{equation}
By Jensen's inequality
\begin{eqnarray*}
(P_{s,s+u}(\varphi(V))-\varphi(V))/u&\leq&
u^{-1}(\varphi(P_{s,s+u}(V))-\varphi(V))\\
&\leq& u^{-1}\left(\varphi\left(V+c~u-\int_0^{u}~(P_{s,s+v}(\varphi(V))~dv\right)-\varphi(V)\right)
\end{eqnarray*}
Letting $u\rightarrow 0$ we conclude that
$$
L_s\left(\varphi(V)\right)\leq (\partial \varphi)(V)\left(c-\varphi(V)\right)\leq c~
 (\partial \varphi)(V)
$$
For any $s\leq u\leq t$, this implies that
$$
P_{s,t}(\varphi(V))\leq P_{s,u}\left(\varphi(V)\right)+c\int_u^{t} P_{s,v}((\partial \varphi)(V))~dv
$$
Integrating $u\in [s,t]$ and using (\ref{sub-ref-append}) we conclude that
\begin{eqnarray*}
(t-s)P_{s,t}(\varphi(V))
&\leq& \int_s^t~P_{s,u}\left(\varphi(V)\right)du+c~\Vert \partial \varphi\Vert
~\frac{(t-s)^2}{2}\\
&\leq &V-P_{s,t}(V)+\left(c~(t-s)+c~\Vert (\partial \varphi)\Vert
~\frac{(t-s)^2}{2}\right)
\end{eqnarray*}
We conclude that
$$
P_{s,s+\tau}(V_{\tau} )=
P_{s,s+\tau}(V+\tau \varphi(V) )\leq V+c_{\tau}=V_{\tau}-\tau \varphi(V)+c_{\tau}$$
with the parameter $c_{\tau}$ defined in (\ref{ref-c-tau-sub}).

This ends the proof of Lemma~\ref{lem-211}.
\cqfd
    \subsection*{Proof of Lemma~\ref{ref-beta-V-est}}\label{betaV-est-proof}
Set 
 $V_{\rho}:=1/2+\rho V$ with $\rho \in ]0,1[$. We associate with these parameters the function
\begin{eqnarray*}
\Delta_{\rho}(x,y)&:=&\frac{\Vert (\delta_x-\delta_y)P\Vert_{V_{\rho}}}{\Vert \delta_x-\delta_y\Vert_{V_{\rho}}}\\
&\leq& \frac{\Vert (\delta_x-\delta_y)M\Vert_{\tiny tv}}{1+\rho(V(x)+V(y))}
+\frac{\rho (P(V)(x)+P(V)(y))}{1+\rho(V(x)+V(y))}.
\end{eqnarray*}
By (\ref{PVC}) (with $c=1/2$), whenever $V(x)+V(y)\geq r> r_0$ we have
\begin{eqnarray*}
\Delta_{\rho}(x,y)
&\leq& \frac{1}{1+\rho(V(x)+V(y))}
+\frac{\rho (V(x)+V(y))}{1+\rho(V(x)+V(y))}~\left(\epsilon+\frac{1}{r}\right)\\
&=& 1-\left(1-\frac{1}{1+\rho(V(x)+V(y))}\right)\left(1-\left(\epsilon+\frac{1}{r}\right)\right).
\end{eqnarray*}
This yields for any $r>r_0\vee r_{\epsilon}$ the estimate
$$
\Delta_{\rho}(x,y)\leq 1-d^1_{\rho}(r)
$$
with
\[
d^1_{\rho}(r):=(1-\epsilon)~\left(1-\frac{1}{1+ \rho r}\right)\left(1-\frac{r_{\epsilon}}{r}\right).
\]
Recalling that $V\geq 1/2$ we readily check that
$
P(V)/V\leq (1+\epsilon)
$  and 
$$
\begin{array}{l}
V(x)+V(y)\leq r\\
\\
\displaystyle\Longrightarrow
\frac{\rho V(x) (P(V)(x)/V(x))+\rho V(y)~(P(V)(y)/V(y))}{1+\rho(V(x)+V(y))}\leq(1+\epsilon)~\frac{\rho r}{1+\rho r},
\end{array}$$
This yields for any $(x,y)$ s.t. $V(x)+V(y)\leq r$  the estimate
\begin{eqnarray*}
\Delta_{\rho}(x,y)
&\leq&1-d^2_{\rho}(r):= \frac{1-\alpha(r)}{1+\rho} +(1+\epsilon)~\frac{\rho r}{1+\rho r}.
\end{eqnarray*}
Choosing $$
\rho=\rho(r):=\frac{\alpha(r)}{2(1+\epsilon)r}\Longrightarrow
\frac{\rho r}{1+\rho r}=\frac{\alpha(r)}{2}~\frac{1}{(1+\epsilon)+\frac{\alpha(r)}{2}}
$$
 we have
$$
d^1_{\rho}(r)=\frac{\alpha(r)}{2}~\frac{(1-\epsilon)}{(1+\epsilon)+ \frac{\alpha(r)}{2}}\left(1-\frac{r_{\epsilon}}{r}\right)<1,
$$
and
$$
1-d^2_{\rho}(r)= \frac{1-\alpha(r)/2}{1+\rho(r)} \leq 1-\alpha(r)/2.
$$
We conclude that
$$
\beta_{V_{\epsilon,r}}(M)=\sup_{x,y\in E}\Delta_{\rho}(x,y)\leq (1-d^1_{\rho}(r))\vee (1-d^2_{\rho}(r))\leq  1-\alpha_{\epsilon}(r).
$$
This ends the proof of Lemma~\ref{ref-beta-V-est}.\cqfd

    \subsection*{Proof of Proposition~\ref{prop-ref-tangent}}\label{prop-ref-tangent-proof}

We have the following almost sure estimate
$\Vert\nabla X_{\tau}(x)\Vert_2\leq e^{-\lambda\tau}
$, where $\Vert A\Vert_2$ stands for the spectral norm of a matrix $A$.
This yields for any $x,y\in\RR^n$ the almost sure estimate
 \begin{equation}\label{ref-nablax-estimate-0-ae-again}
 \Vert X_{\tau}(x)-X_{\tau}(y)\Vert\leq  e^{-\lambda\tau}~\Vert x-y\Vert.
\end{equation}
Applying the above to $y=0$ we find that
$$
P_{\tau}(V)(x)\leq P_{\tau}(V)(0)~V(x)^{1-\delta}\quad \mbox{with}\quad  \delta=1-e^{-\lambda\tau}.$$
Next, we check that $P^X_{\tau}(V)(0)<\infty$. We have
$$
X_{u}(0)=\int_0^{u} \left(b(0)~ds+\sigma~dB_s\right)+\int_0^{u} \int_0^1~\nabla b(\epsilon X_{s}(0))^{\prime}~X_{s}(0)~d\epsilon~ds.
$$
This implies that
$$
\Vert X_{u}(0)\Vert\leq \beta~\left(u +\Vert B_{u}\Vert\right)+\beta
\int_0^{u}~\Vert X_{s}(0)\Vert~ds$$
with
$\beta:= \sigma\vee \Vert b(0)\Vert\vee \Vert \nabla b\Vert$.
Applying Gr\"onwall lemma we check that
$$
\Vert X_{u}(0)\Vert\stackrel{\tiny law}{\leq } 
\beta~\left(u + \Vert B_u\Vert\right)+\beta^2\int_0^{u}\left(s +\Vert B_{s}\Vert\right)e^{\beta (u-s)}~ds.
$$
On the other hand, we have
$$
\int_0^{u}\Vert B_{s}\Vert~ds=u~\int_0^{1}\Vert B_{us}\Vert~ds\stackrel{\tiny law}{=} 
u^{3/2}~\int_0^{1}\Vert B_{s}\Vert~ds.
$$
This yields the rather crude estimate
$$
\Vert X_{u}(0)\Vert\stackrel{\tiny law}{\leq } 
\beta~\left(u +~u^{1/2}~\Vert B_1\Vert\right)+\beta^2~u^2/2 + \beta^2 e^{\beta u}~u^{3/2}~\int_0^{1}\Vert B_{s}\Vert~ds.
$$
For any $a\geq 0$ by Jensen's inequality
$$
\EE\left(e^{a\int_0^{1}\Vert B_{s}\Vert~ds}\right)\leq 
\int_0^{1}~\EE\left(e^{a\Vert B_{s}\Vert}\right)~ds\leq
e^{a^2r/2}. 
$$
It is now an elementary exercise to check that $\EE(e^{v\Vert X_{\tau}(0)\Vert})<\infty$. 
This ends the proof of the proposition.
\cqfd
  
\subsection*{Proof of Proposition~\ref{condition-ks-prop-0}}\label{condition-ks-prop-proof}

Consider the function
$$
\begin{array}{l}
\displaystyle f_t(x):=\exp{\left(2\epsilon~\left(e^{-\alpha t}~W(x)-\beta~\frac{1-e^{-\alpha t}}{\alpha}\right)\right)}\\
\\
\Longrightarrow -\partial_t\log f_t(x)=2\epsilon~e^{-\alpha t}~(\alpha W(x) +\beta).
\end{array}$$
In the same vein, we check that
\begin{eqnarray*}
\partial_{x_i}f_t(x)/f_t(x)&=&2\epsilon~e^{-\alpha t}~\partial_{x_i} W\\
\partial_{x_i,x_j}f_t(x)/f_t(x)&=&2\epsilon~e^{-\alpha t}~\left(2\epsilon~e^{-\alpha t}~
\partial_{x_i} W \partial_{x_j}W+\partial_{x_i,x_j}W\right).
\end{eqnarray*}
This implies that
$$
\begin{array}{l}
\displaystyle\left(L(f_t)-\partial_t f_t\right)/f_t\\
\\
\displaystyle=2\epsilon~e^{-\alpha t}~\left((\alpha W +\beta)+L(W)+\epsilon~e^{-\alpha t}~\Gamma_{L}(W,W)\right).\end{array}
$$
Combining the above with (\ref{condition-ks-prop}) we find that
$$
\begin{array}{l}
\displaystyle\left(L(f_t)(x)-\partial_t f_t(x)\right)\leq -2~\epsilon^2~e^{-\alpha t}~\left(1-e^{-\alpha t}\right)
\Gamma_{L}(W,W)(x)~f_t(x)\leq 0.
\end{array}
$$
This yields the interpolation formula
$$
\EE\left(f_0(X_t(x))\right)-f_t(x)=\int_0^t \EE\left(\partial_s (f_{t-s}(X_s(x)))\right)~ds\leq 0.
$$
We check (\ref{expo-estimate-LW-Gamma}) after some elementary manipulations, thus there are skipped. This ends the proof of the proposition.\cqfd

\subsection*{Proof of Proposition~\ref{solv-harmonic-P}}\label{solv-harmonic-P-proof}  
Notice that
$$
X^h_t(x)\stackrel{\rm law}{=}\epsilon_t~ x+\sigma_t~Z=B_{\sigma_t}\left(\epsilon_t~ x\right)\quad \mbox{\rm with}\quad\epsilon_t:=e^{-t}\quad \mbox{\rm and}\quad
\sigma_t:=\sqrt{\frac{1-\epsilon_t^2}{2}}
$$
and some centered Gaussian random variable $Z$ with unit variance.
The conjugate formula (\ref{conjugate-formula})
yields the integral operator equation
$$
Q_t(x,dy)=e^{- t/2}~e^{-x^2/2}~\frac{1}{\sqrt{2\pi\sigma_t^2}}~\exp{\left(-\frac{(y-\epsilon_tx)^2}{2\sigma^2_t}+\frac{y^2}{2}\right)}~dy.
$$
Observe that
$$
-\frac{(y-\epsilon_tx)^2}{\sigma^2_t}+y^2=-\frac{1}{p_t}~\left(y-\frac{\epsilon_t}{1-\sigma_t^2}~x\right)^2+x^2~\frac{\epsilon^2_t }{1-\sigma_t^2}
$$
with
\begin{equation}\label{def-p-t}
p_t:=\frac{1-\epsilon^2_t}{1+\epsilon_t^2}=\tanh(t)\quad\Longleftrightarrow\quad\partial_t p_t=1-p_t^2\quad \mbox{\rm with}\quad p_0=0.
\end{equation}
We check this claim using the fact that
$$
\frac{1}{\sigma^2_t}=\frac{2}{1-\epsilon^2_t}=1+\frac{1+\epsilon_t^2}{1-\epsilon^2_t}=1+\frac{1}{p_t}.
$$
On the other hand, we have
$$
\frac{1-\sigma^2_t}{\epsilon_t}=\cosh(t)\quad \mbox{\rm and}\quad
\partial_t\log{\cosh(t)}=p_t=\tanh(t).
$$
This implies that
$$
\int_0^t p_s~ds=\log{\cosh(t)}
\quad \mbox{\rm and}\quad
\frac{\epsilon_t}{1-\sigma^{2}_t}=\frac{1}{\cosh(t)}=\exp{\left(-\int_0^t p_sds\right)}.
$$
We also have
$$
1-\sigma^2_t=1-\frac{1-\epsilon^2_t}{2}=\frac{1+\epsilon^2_t}{2}\Longrightarrow
1-\frac{\epsilon^2_t }{1-\sigma_t^2}=\frac{1-\epsilon^2_t}{1+\epsilon^2_t}=p_t.
$$
This implies that
$$
Q_t(1)(x)=e^{-t/2}~\frac{\sqrt{p_t}}{\sigma_t}~\exp{\left(-\frac{x^2}{2}~p_t\right)}=e^{-t/2}~h(x)~P^h_t(1/h)(x).
$$
Notice that
$$
e^{-t/2}~\frac{\sqrt{p_t}}{\sigma_t}=\sqrt{\epsilon_t~\frac{1-\epsilon^2_t}{1+\epsilon_t^2}\frac{2}{1-\epsilon^2_t}}=\sqrt{\frac{2}{1/\epsilon_t+\epsilon_t}}=\frac{1}{\sqrt{\cosh(t)}}
$$
and
$$
\partial_t m_t(x)=-p_t~m_t(x)\quad\mbox{and}\quad 
\partial_t p_t=1-p_t^2\quad \mbox{with}\quad (m_0(x),p_0)=(x,0).
$$
This ends the proof of the proposition.\cqfd

\subsection*{Proof of Proposition~\ref{prop-half-harmonic}}\label{prop-half-harmonic-proof}
Notice that
$$
\begin{array}{l}
e^{t/2}~e^{x^2/2}~Q_t(x,dy)\\
\\
\displaystyle=\frac{1}{\sqrt{2\pi\sigma_t^2}}~\left(
\exp{\left(-\frac{(y-\epsilon_tx)^2}{2\sigma^2_t}+\frac{y^2}{2}\right)}-\exp{\left(-\frac{(y+\epsilon_tx)^2}{2\sigma^2_t}+\frac{y^2}{2}\right)}\right)~1_{[0,\infty[}(y)~dy.
\end{array}$$
This implies that
$$
\begin{array}{l}
\displaystyle Q_t(1)(x)= \frac{e^{-\frac{x^2}{2}~\tanh(t)}}{\sqrt{\cosh(t)}}~\\
\\
\hskip3cm\displaystyle\times \int_0^{\infty}\frac{1}{\sqrt{2\pi p_t}}~\left(
\exp{\left(-\frac{(y-m_t(x))^2}{2p_t}\right)}-\exp{\left(-\frac{(y+m_t(x))^2}{2p_t}\right)}\right)~dy.
\end{array}$$
We conclude that
$$
\begin{array}{l}
\displaystyle Q_t(1)(x)
\displaystyle= \frac{e^{-\frac{x^2}{2}~\tanh(t)}}{\sqrt{\cosh(t)}}\times \PP\left(-m_t(x)/\sqrt{p_t}\leq Z\leq m_t(x)/\sqrt{p_t} \right)\\
\\
\displaystyle= 2~\frac{e^{-\frac{x^2}{2}~\tanh(t)}}{\sqrt{\cosh(t)}}\times \PP\left(0\leq Z\leq \frac{x}{\sqrt{\sinh(t)\cosh(t)}} \right)\\
\\
\displaystyle\longrightarrow 0\quad \mbox{\rm as $x\rightarrow\infty$ or  $x\rightarrow 0$ or as $t\rightarrow\infty$.}
\end{array}$$
In the above display, $Z$  stands for some centered Gaussian random variable with unit variance.
Note that we have used the fact that
$$
m_t(x)/\sqrt{p_t}=\frac{x}{\cosh(t)~\sqrt{\tanh(t)}}=\frac{x}{\sqrt{\sinh(t)\cosh(t)}}.
$$
In addition, we have 
$$
\begin{array}{l}
\displaystyle
\overline{Q}_t(x,dy)=\displaystyle\frac{1}{ \PP\left(0\leq Z\leq {x}/{\sqrt{\sinh(t)\cosh(t)}} \right)}~\\
\\
\displaystyle\times \frac{1}{2\sqrt{2\pi p_t}}~\left(
\exp{\left(-\frac{(y-m_t(x))^2}{2p_t}\right)}-\exp{\left(-\frac{(y+m_t(x))^2}{2p_t}\right)}\right)~1_{[0,\infty[}(y)~dy.
\end{array}$$
This ends the proof of the proposition.\cqfd

\subsection*{Proof of (\ref{langevin-Lac})}\label{langevin-Lac-proof}

The generator of the process (\ref{Langevin-2d}) is defined by
$$
L(f)(q,p)=\beta~\frac{p}{m}~\frac{\partial f}{\partial q}~-\beta~\left(\frac{\partial W}{\partial q}+\frac{\sigma^2}{2}~\frac{p}{m}\right)
\frac{\partial f}{\partial p}+\frac{\sigma^2}{2}~\frac{\partial^2 f}{\partial p^2}.
$$

Recalling that $2pq\leq p^2+q^2$, we prove that
\begin{eqnarray*}
V(q,p)&\leq& \frac{1}{2}\left(\frac{1}{m}+\epsilon\right)~p^2+\frac{\epsilon}{2}~\left(\frac{\sigma^2}{2}+1\right)~q^2+W(q)\\
&\leq &C^{\star}(\epsilon)~\left(1+p^2+q^2+W(q)\right)
\end{eqnarray*}
with $$C^{\star}(\epsilon):=\max{\left\{\frac{1}{2}\left(\frac{1}{m}+\epsilon\right),\frac{\epsilon}{2}~
\left(\frac{\sigma^2}{2}+1\right),1\right\}}.$$
On the other hand, we have
\begin{eqnarray*}
L(V)&=&\beta\frac{p}{m}\left(\frac{\partial W}{\partial q}+\epsilon~\frac{\sigma^2}{2}~q+\epsilon~p\right)
\\
&&\hskip1cm -\beta~\left(\frac{\partial W}{\partial q}+\frac{\sigma^2}{2}~\frac{p}{m}\right)
~\left(\frac{p}{m}+\epsilon~q\right)
+\frac{\sigma^2}{2m}\\
&=&-\beta~\left[\frac{1}{m}~\left(\frac{\sigma^2}{2m}-\epsilon~\right)~p^2+\epsilon~q~\frac{\partial W}{\partial q}\right]+\frac{\sigma^2}{2m}.
\end{eqnarray*}
Under our assumptions, this implies that for any $\vert q\vert\geq r$ we have
\begin{eqnarray*}
L(V)
&\leq &-\beta~\left[\frac{1}{m}~\left(\frac{\sigma^2}{2m}-\epsilon~\right)~p^2+\epsilon~\delta~\left(W(q)+q^2\right)\right]+\frac{\sigma^2}{2m}\\
&\leq &-C_{\star}(\epsilon,\delta)~\left(1+p^2+q^2+W(q)\right)+c_m(\epsilon,\delta)\end{eqnarray*}
with
$$
C_{\star}(\epsilon,\delta):=\beta~ \min{\left\{
\left(
\frac{1}{m}~\left(\frac{\sigma^2}{2m}-\epsilon~\right), \epsilon~\delta
\right)
\right\}}\quad \mbox{and}\quad c_m(\epsilon,\delta):=C_{\star}(\epsilon,\delta)+\frac{\sigma^2}{2m}.
$$
We conclude that for any $\vert q\vert> r$, 

\begin{eqnarray*}
(V^{-1}L(V))(q,p)&
\leq& -\frac{C_{\star}(\epsilon,\delta)~\left(1+p^2+q^2+W(q)\right)-c_m(\epsilon,\delta)}{V(q,p)}\\
&\leq& -\frac{C_{\star}(\epsilon,\delta)~\left(1+p^2+q^2+W(q)\right)-c_m(\epsilon,\delta)}{C^{\star}(\epsilon)~\left(1+p^2+q^2+W(q)\right)}\\
&=&-\frac{C_{\star}(\epsilon,\delta)}{C^{\star}(\epsilon)}+\frac{c_m(\epsilon,\delta)}{C^{\star}(\epsilon)}~\frac{1}{1+p^2+q^2+W(q)}\\
&\leq& -\left[\frac{C_{\star}(\epsilon,\delta)}{C^{\star}(\epsilon)}-\frac{c_m(\epsilon,\delta)}{C^{\star}(\epsilon)}~\frac{1}{1+p^2+q^2}\right].
\end{eqnarray*}
We choose $r$ sufficiently large to satisfy
$$
\begin{array}{l}
\vert p\vert> r\quad\mbox{or}\quad
\vert q\vert > r\\
\\
\displaystyle\Rightarrow
\frac{C_{\star}(\epsilon,\delta)}{C^{\star}(\epsilon)}-\frac{c_m(\epsilon,\delta)}{C^{\star}(\epsilon)}~\frac{1}{p^2+q^2}\geq \frac{C_{\star}(\epsilon,\delta)}{C^{\star}(\epsilon)}-\frac{c_m(\epsilon,\delta)}{C^{\star}(\epsilon)}~\frac{1}{r^2}\geq a:=
\frac{C_{\star}(\epsilon,\delta)}{2C^{\star}(\epsilon)}> 0,
\end{array}$$
and we set
$$
K_r:=\{(q,p)\in \RR^2~:~\vert p\vert\vee\vert q\vert\leq  r\}.
$$
In this notation, we have
$$
L(V)\leq -a V~1_{E-K_r}+\sup_{K_r}L(V)\leq -aV~+c
\quad\mbox{
with}\quad
c:=\sup_{K_r}L(V)+a\sup_{K_r}V.
$$

\subsection*{Proof of (\ref{half-harm-lyap})}\label{half-harm-lyap-proof}
Observe that for any $0<y\leq 1$ and $z\in E=]0,\infty[$ we have
$$
\sinh{(yz)}\leq y~\sinh{(z)}\quad \mbox{and}\quad
\sinh{(z)}\leq \frac{1}{2}~e^{z}.
$$
This implies that
$$
\int_0^{\infty}~\Qa_t(x,dy)~\frac{1}{y}~
\leq \sinh{(m_t(x))}~\frac{e^{-\frac{x^2}{2}(p_t+e^{-2t}/p_t)}}{\sqrt{\cosh(t)}}~\displaystyle\times \int_0^{\infty}~\sqrt{\frac{2}{\pi p_t}}~\exp{\left(-\frac{y^2}{2p_t}\right)}~dy.
$$
from which we check that
$$
\int_0^{\infty}~\Qa_t(x,dy)~\frac{1}{y}~1_{]0,1]}(y)\leq~\frac{\exp{\left(-\left(
\frac{x^2}{2}(p_t+\frac{e^{-2t}}{p_t})-e^{-t}x\right)\right)}}{2\sqrt{\cosh(t)}}.$$

On the other hand, for any $n\geq 1$ we have
$$
\begin{array}{l}
\displaystyle
\int_{0}^{\infty}~\Qa_t(x,dy)~y^n\\
\\
\displaystyle\leq  \frac{1}{2}~~\frac{1}{\sqrt{\cosh(t)}}~\sqrt{\frac{2}{\pi p_t}}~
\exp{\left(-\frac{\epsilon_t^2 x^2}{2p_t}-\frac{x^2}{2}~p_t\right)}
\int_0^{\infty}~y^n~\exp{\left(y\epsilon_t x-\frac{y^2}{2p_t}\right)}~dy.
\end{array}$$
Notice that
$$
y\epsilon_t x-\frac{y^2}{2p_t}=-\frac{1}{2p_t}\left(y-\epsilon_t x p_t\right)^2+
\frac{x^2}{2}~\epsilon^2_t~p_t
$$
so that
$$
\begin{array}{l}
\displaystyle
\int_0^{\infty}~\Qa_t(x,dy)~y^n~\leq  \frac{1}{2}~~\frac{1}{\sqrt{\cosh(t)}}~\sqrt{\frac{2}{\pi p_t}}~\\
\\
\hskip1cm\displaystyle\times
\exp{\left(-\frac{x^2}{2}~\left((1-\epsilon^2_t)~p_t+
\frac{\epsilon_t^2}{p_t}
\right)\right)}
\int_0^{\infty}~y^n~\exp{\left(-\frac{1}{2p_t}\left(y-\epsilon_t x p_t\right)^2\right)}~dy.
\end{array}$$
For any $n\geq 1$, we conclude that 
$$
V(x):=x^n+{1}/{x}\Longrightarrow V\in \Ca_{\infty}(E)\quad \mbox{and}\quad \Vert Q_t(V)\Vert<\infty.
$$
This ends the proof of (\ref{half-harm-lyap}).\cqfd

\subsection*{Proof of Lemma~\ref{lem-LU}}\label{lem-LU-proof}
To simplify notation, we write $Q_{t}$ instead of $Q_{t}^{[U]}$.  
For any $V\in \Ba_{\infty}(E)\cap \Da(L)$ we have
\begin{eqnarray*}
Q_{t}(V)&=&V+\int_0^t Q_s(L(V)-UV)~ds\\
&\leq&V+\int_0^t \left[-a~Q_s(V)+c~Q_s(1)\right]~ds=V+c\int_0^tQ_s(1)~ds-a\int_0^t Q_s(V)ds.
\end{eqnarray*}
On the other hand, through integration by parts we have
\begin{eqnarray*}
\int_0^t Q_s(V)ds&=&\left[s~Q_s(V)\right]_0^t-\int_0^t~s~\frac{d}{ds}Q_s(W)~ds\\
&=&t~Q_t(V)-\int_0^t~s~Q_s(\underbrace{L(V)-UV}_{\leq c})~ds\geq t~Q_t(V)-c\int_0^t~s~Q_s(1)ds.
\end{eqnarray*}
This implies that
\begin{eqnarray*}
Q_{t}(V)
&\leq &V+c\int_0^tQ_s(1)~ds-a\left(t~Q_t(V)-c\int_0^t~s~Q_s(1)ds\right)
\end{eqnarray*}
from which we conclude that
\begin{eqnarray*}
Q_{t}(V)
&\leq &\frac{V}{1+at}+c\int_0^t~Q_s(1)ds \Longrightarrow Q_{t}(V)
\leq \frac{V}{1+at}+ct.
\end{eqnarray*}
This ends the proof of (\ref{L-U-V-ref}). Now, we come to the proof of (\ref{L-U-ref}).
We have the forward evolution equation given for any $f\in \Da(L)$ by
$$
\partial_tQ_{t}(f)=Q_{t}(L^U(f)).
$$
Applying the above to $f=U$ we readily check that
$$
\partial_tQ_{t}(U)\leq a_0+a_1~Q_{t}(U)-Q_{t}(U^2)\leq a_0+a_1~Q_{t}(U)-(Q_{t}(U))^2/Q_{t}(1)
$$
from which we find the Riccati estimates
$$
\partial_tQ_{t}(U)\leq  a_0+a_1~Q_{t}(U)-(Q_{t}(U))^2\Longrightarrow \forall t>0\quad \Vert Q_{t}(U)\Vert<\infty.
$$

This ends the proof of the lemma.
\cqfd

\subsection*{Proof of (\ref{langevin-ref-intro})}\label{langevin-ref-intro-proof}

By Girsanov theorem we have
$$
\Qa_t^{(a)}(f)(z)=\EE\left(f(\Xa^0_t(z))~Z_{t}(z)~1_{T^0(z)>t}\right)
$$
 with the exponential martingale
$$
Z_{t}(z)=\exp{\left(\frac{1}{\sigma}\int_0^{t} a(X_{s}(z))^{\prime}~dB_s-\frac{1}{2\sigma^2}~\int_0^{t} \Vert a(X_{s}(z))\Vert^2~ds\right)}.
$$
By H\"older's inequality, for any non negative function $f$ on $E$, any $z\in E$ and any 
conjugate parameters $p,q> 1$ with $1/p+1/q=1$ we have
$$
\Qa_t^{(a)}(f)(z)\leq
\EE\left(Z_{t}(z)^q~1_{T^0(z)>t}\right)^{1/q}~ \Qa_t^{(a)}(f^p)(z)^{1/p}. 
$$
On the other hand, we have
\begin{eqnarray*}
\EE\left(Z_{t}(z)^q~1_{T^0(z)>t}\right)&=&
\EE\left(\overline{Z}_{t}(z)~\exp{\left(\frac{q(q-1)}{2\sigma^2}~\int_0^{t} \Vert a(X_{s}(z))\Vert^2~ds\right)}~1_{T^0(z)>t}\right)\\
&\leq& c_t(p):=\exp{\left(\frac{p t}{2((p-1)\sigma)^2}~\sup_Da\right)}
\end{eqnarray*}
with the exponential martingale
$$
\overline{Z}_{t}(z)=\exp{\left(\frac{q}{\sigma}~\int_0^{t} a(X_{s}(z))^{\prime}~dB_s-\frac{q^2}{2\sigma^2}~\int_0^{t} \Vert a(X_{s}(z))\Vert^2~ds\right)}.
$$
This ends the proof of roof of (\ref{langevin-ref-intro}).\cqfd

  \subsection*{Proof of Lemma~\ref{lem-lip-dom-compact}}\label{lip-dom-compact-proof}
For any $z\in\partial E$ there exists some open ball $\BB(z,r)\subset \RR^n$ with $r>0$ and some $\Ca^1$-mapping $g$ from $\RR^{n-1}$ into $\RR$ such that
\begin{eqnarray*}
E\cap \BB(z,r)&=&\{ x\in \BB(z,r)~:~x_n< g(x_{-n}) \}\\
 \partial E\cap \BB(z,r)&=&\{ x\in \BB(z,r)~:~x_n=g(x_{-n}) \}\quad \mbox{\rm with}\quad x_{-n}:=(x_1,\ldots,x_{n-1}).
\end{eqnarray*}
We make the change of variables
$$
\begin{array}{l}
\Ea(z,r):=E\cap \BB(z,r)\\
\\
\mapsto\varsigma(x):=
(x_{-n},x_n-g(x_{-n}))\in \Oa(z,r):=\varsigma(\Ea(z,r))\subset (\RR^{n-1}\times\RR_+)
\end{array}
$$
with Jacobian
$$
\nabla \varsigma(x)=\left(
\begin{array}{cc}
I_{(n-1)\times (n-1)}&-\nabla g(x_{-n})\\
0&1
\end{array}\right). 
$$
Observe that
$$
\begin{array}{l}
\varsigma~:~x\in  \Ea_0(z,r):=\left(\partial E\cap \BB(z,r)\right)\\
\\
\Longrightarrow
\varsigma(x)=(x_{-n},0)\in \Oa_0(z,r):=\varsigma( \Ea_0(z,r))\subset (\RR^{n-1}\times\{0\}).
\end{array}
$$
 The inverse is given by
$$
\begin{array}{l}
y\in \Oa(z,r)\mapsto \varsigma^{-1}(y)=
(y_{-n},y_n+g(y_{-n}))\in \Ea(z,r)\\
\\
\Longrightarrow\nabla \varsigma^{-1}(y)=\left(
\begin{array}{cc}
I_{(n-1)\times (n-1)}&\nabla g(y_{-n})\\
0&1
\end{array}\right).
\end{array}
$$
On the other hand we have
\begin{eqnarray*}
\Vert \varsigma(x)-\varsigma(\overline{x})\Vert&=&\left(\Vert x_{-n}-\overline{x}_{-n}\Vert^2+\left(\vert x_n-\overline{x}_n\vert+\vert g(x_{-n})-g(\overline{x}_{-n})\vert \right)^2\right)^{1/2}\\
&\leq&\left(\Vert x_{-n}-\overline{x}_{-n}\Vert^2+2\vert x_n-\overline{x}_n\vert^2+2\Vert \nabla g\Vert^2\Vert x_{-n}-\overline{x}_{-n}\Vert^2 \right)^{1/2}\\
&\leq& c(g)~\Vert x-\overline{x}\Vert\quad \mbox{\rm with}\quad c(g):=\sqrt{2\vee(1+2\Vert \nabla g\Vert^2)}\geq 1.
\end{eqnarray*}
In the same vein, we have
$$
\Vert \varsigma^{-1}(y)-\varsigma^{-1}(\overline{y})\Vert\leq  c(g)~\Vert y-\overline{y}\Vert
\quad\mbox{\rm so that}\quad
\frac{1}{c(g)}~\Vert y-\overline{y}\Vert\leq \Vert \varsigma^{-1}(y)-\varsigma^{-1}(\overline{y})\Vert.
$$
For  any $x\in \Ea(z,r)$ and $\overline{x}\in \Ea_0(z,r)$
we have $\varsigma(\overline{x})\in \Oa_0(z,r)$ and
$$
\Vert x-\overline{x}\Vert=\Vert \varsigma^{-1}(\varsigma(x))-\varsigma^{-1}(\varsigma(\overline{x}))\Vert\geq \frac{1}{c(g)}~\Vert \varsigma(x)-\varsigma(\overline{x})\Vert \geq \frac{1}{c(g)}~\vert \varsigma(x)_n\vert.
$$
Taking the infimum of all $\overline{x}\in \Ea_0(z,r)$ this implies that
$$
d(x,\Ea_0(z,r))\geq \frac{1}{c(g)}~\vert \varsigma(x)_n\vert\quad \mbox{\rm and}\quad
d(\varsigma^{-1}(y),\Ea_0(z,r))\geq \frac{1}{c(g)}~\vert y_n\vert
$$
for any $x\in \Ea(z,r)$ and $y\in \Oa(z,r)$. We conclude that
$$
\begin{array}{l}
\displaystyle\int_{\Ea(z,r)}~\cchi\left(d(x,\Ea_0(z,r))\right)~dx\\
\\
\displaystyle=
\int_{\Oa(z,r)}~\cchi\left(d(\varsigma^{-1}(y),\Ea_0(z,r))\right)~\vert\mbox{\rm det}\left(\varsigma^{-1}(y)\right)\vert~dy\\
\\
\displaystyle\leq \frac{1}{c(g)}~\sup_{y\in \Oa(z,r)}\vert\mbox{\rm det}\left(\varsigma^{-1}(y)\right)\vert~\int_{\Oa(z,r)}~\cchi( y_n)~dy<\infty.
\end{array}
$$
We end the proof of the lemma by covering $\partial E$ by finitely many boundary 
coordinates patches $(\Ea(z_i,r_i),g_i)_{1\leq i\leq n}$, for some $z_i\in \partial E$, $r_i>0$ and some local defining functions $g_i$.
 \cqfd

\subsection*{Proof of Lemma~\ref{lem-boundary-f-g}}\label{lem-boundary-f-g-proof}

Using the change of variable formulae
$$
\int_{\partial E_{r}}~f(z)~\sigma_{\partial,r}(dz)=\int_{\partial E}~f\left(z+rN(z)\right)~\vert \mbox{\rm det}\left(I-r~\Wa(z)\right)\vert~\sigma_{\partial}(dz)
$$
and
$$
\int_{\partial E}~f(z)~\sigma_{\partial}(dz)=\int_{\partial E_r}~f\left(z-rN(z)\right)~\vert \mbox{\rm det}\left(I+r~\Wa(z)\right)\vert~\sigma_{\partial,r}(dz)
$$
we check that
$$
\int_{\partial E_{r}}~f(z)~\sigma_{\partial,r}(dz)\leq \kappa_{\partial}(\alpha)~
\int_{\partial E}~f\left(z+rN(z))\right)~\sigma_{\partial}(dz)
$$
and
$$
\int_{\partial E}~f(z)~\sigma_{\partial}(dz)\leq  \kappa_{\partial}^-(\alpha)\int_{\partial E_{r}}~f\left(z-rN(z)\right)~~\sigma_{\partial,r}(dz).
$$
This yields the estimate
$$
\int_{\partial E_{r}}~f(z)~\sigma_{\partial,r}(dz)\leq ~\iota(\alpha)~\kappa_{\partial}(\alpha)~
\int_{\partial E}~g(z)~\sigma_{\partial}(dz).
$$
In the same vein, we have 
$$
\int_{\partial E}~f(z)~\sigma_{\partial}(dz)\leq \iota(\alpha) \kappa_{\partial}^-(\alpha)\int_{\partial E_{r}}~g\left(z\right)~~\sigma_{\partial,r}(dz).
$$
Integrating w.r.t. the parameter $r\in [0,\alpha]$ we check the co-area estimate
\begin{eqnarray*}
\alpha~\int_{\partial E}~f(z)~\sigma_{\partial}(dz)&\leq&  \iota(\alpha) \kappa_{\partial}^-(\alpha)
\int_0^{\alpha}~dr~\int_{\partial E_r}~g(z)~\sigma_{\partial,r}(dz)~
\\
&=& \iota(\alpha)~\kappa_{\partial}^-(\alpha)\int_{\Da_{\alpha}(E)}~g(z)~dz.
\end{eqnarray*}
This ends the proof of the lemma.
\cqfd

\subsection*{Proof of Proposition~\ref{prop-vol-form-weignarten}}\label{prop-vol-form-weignarten-proof}
For any given  $\overline{\theta}:=(\overline{\theta}_1,\ldots,\overline{\theta}_n)\in \left(\RR^{n-1}\times [0,r]\right)$ we set $\overline{\theta}_{-n}:=(\theta_1,\ldots,\theta_{n-1})$. In this notation, we have
$$
\begin{array}{rcl}
\overline{\psi}~:~\overline{\theta}\in \left(\RR^{n-1}\times [0,r]\right)\mapsto
\overline{\psi}(\overline{\theta})&:=&
F(\psi(\overline{\theta}_{-n}),\overline{\theta}_n)\\
&&\\
&=&\psi(\overline{\theta}_{-n})+ \overline{\theta}_n~N(\psi(\overline{\theta}_{-n}))\in \Da_{r}(E).
\end{array}
$$

The volume form $\sigma_{\Da_{r}(E)}$ on $\Da_{r}(E)$ expressed in the chart $\overline{\psi}$ is given by
$$
\left(\sigma_{\Da_{r}(E)}\circ \overline{\psi}^{-1}\right)(d\overline{\theta})=\vert \mbox{\rm det}\left(\mbox{\rm Jac}(\overline{\psi})(\theta)\right)\vert =\sqrt{\mbox{\rm det}\left(\left(\partial \overline{\psi}(\overline{\theta})\right)\left(\partial \overline{\psi}(\overline{\theta})\right)^{\prime}\right)}~
d\overline{\theta}.
$$

Arguing as above, we have
$$
\left(\partial \overline{\psi}(\overline{\theta})\right)^{\prime}=\left(\left(\partial_{\overline{\theta}_{-n}}\overline{\psi}(\overline{\theta})\right)^{\prime},\partial_{\overline{\theta}_n}\overline{\psi}(\overline{\theta})\right)\in \left(\TT_{\overline{\psi}(\overline{\theta})}(\Da_{r}(E))\right)^{n}
$$
with the tangent vectors
$$
\left(\partial_{\overline{\theta}_{-n}}\overline{\psi}(\overline{\theta})\right)^{\prime}:=\left(
\partial_{\overline{\theta}_1}\overline{\psi}(\overline{\theta}),\ldots,\partial_{\overline{\theta}_{n-1}}\overline{\psi}(\overline{\theta})\right)
\quad \mbox{\rm and}\quad
\partial_{\overline{\theta}_n}\overline{\psi}(\overline{\theta})=N(\psi(\overline{\theta}_{-n})).
$$
In addition, we have
\begin{eqnarray*}
\left(\partial_{\overline{\theta}_{-n}}\overline{\psi}(\overline{\theta})\right)^{\prime}&=&
\left(\partial \psi(\overline{\theta}_{-n})\right)^{\prime}+\overline{\theta}_n~\left(\partial (N(\psi(\overline{\theta}_{-n})))\right)^{\prime}=\left(\partial \psi(\overline{\theta}_{-n})\right)^{\prime}\left(I-\overline{\theta}_n~\Wa(\psi(\overline{\theta}_{-n}))\right).
\end{eqnarray*}
This yields the formula
$$
\left(\partial_{\overline{\theta}_{-n}}\overline{\psi}(\overline{\theta})\right)\left(\partial_{\overline{\theta}_{-n}}\overline{\psi}(\overline{\theta})\right)^{\prime}
=g(\psi(\overline{\theta}_{-n}))~\left(I-\overline{\theta}_n~\Wa(\psi(\overline{\theta}_{-n}))\right)^2
$$
from which we check that
$$
\left(\partial \overline{\psi}(\overline{\theta})\right)\left(\partial \overline{\psi}(\overline{\theta})\right)^{\prime}=\left(
\begin{array}{cc}
g(\psi(\overline{\theta}_{-n})\left(I-\overline{\theta}_n~\Wa(\psi(\overline{\theta}_{-n}))\right)^2&0_{n-1}\\
0&1
\end{array}
\right).
$$
We conclude that
$$
\sqrt{\mbox{\rm det}\left(\left(\partial \overline{\psi}(\overline{\theta})\right)\left(\partial \overline{\psi}(\overline{\theta})\right)^{\prime}\right)}=\sqrt{\mbox{\rm det}(g(\psi(\overline{\theta}_{-n}))}~\left\vert\mbox{\rm det}\left(I-\overline{\theta}_n~\Wa(\psi(\overline{\theta}_{-n}))\right)\right\vert
$$
and therefore
$$
\left(\sigma_{\Da_{r}(E)}\circ \overline{\psi}^{-1}\right)(d\overline{\theta})=
\left\vert\mbox{\rm det}\left(I-\overline{\theta}_n~\Wa(\psi(\overline{\theta}_{-n}))\right)\right\vert
~d\overline{\theta}_n~\left(\sigma_{\partial,0}\circ \psi^{-1}\right)(d\overline{\theta}_{-n}).
$$
For any given $\overline{\theta}_n=u\in [0,r]$, the volume form $\sigma_{\partial,u}$ on the boundary $\partial E_{u}$ expressed in  the boundary chart
$$
\overline{\psi}(\point,u)~:~
\theta\in\RR^{n-1}\mapsto \overline{\psi}(\theta,u)=F(\psi(\theta),u)\in \partial E_{u}
$$
  is given by
$$
\begin{array}{l}
\left(\sigma_{\partial,u}\circ \overline{\psi}(\point,u)^{-1}\right)(d\theta)
=\left\vert\mbox{\rm det}\left(I-u~\Wa(\psi(\theta))\right)\right\vert~\left(\sigma_{\partial,0}\circ \psi^{-1}\right)(d\theta).
\end{array}$$
This ends the proof of the  proposition.
\cqfd

\end{document}